\documentclass[leqno,12pt, twoside]{article}
\usepackage{ifthen}
\usepackage{soul}
\usepackage{latexsym}
\usepackage{amsthm,amssymb}
\usepackage{amsbsy,amsfonts,amsmath}

\newcommand{\bysame}{\leavevmode\hbox to 3em{\hrulefill}\,}

\usepackage{a4}

\numberwithin{equation}{section}

\date{}

\newtheorem{theorem}{Theorem}[section]
\newtheorem{proposition}[theorem]{Proposition}
\newtheorem{corollary}[theorem]{Corollary}
\newtheorem{lemma}[theorem]{Lemma}

\newtheorem{definition}[theorem]{Definition}

\theoremstyle{remark}
\newtheorem{remark}{Remark}
\newtheorem{notation}[remark]{Notation}

\newcommand{\al}{\alpha}
\newcommand{\stbsis}{e}
\newcommand{\stProj}{P}

\newcommand{\bC}{\mathbb{C}}

\newcommand{\bK}{\mathbb{K}}
\newcommand{\bKt}{\bK^\times}
\newcommand{\bKtinf}{{\bKt_\infty}}

\newcommand{\bR}{{\mathbb{R}}}
\newcommand{\bN}{{\mathbb{N}}}
\newcommand{\bQ}{{\mathbb{Q}}}
\newcommand{\bZ}{{\mathbb{Z}}}
\newcommand{\bZgeqo}{{\bZ_{\geq 0}}}

\newcommand{\fkJ}{J}

\newcommand{\rmGL}{{\rm {GL}}}
\newcommand{\rmid}{{\rm {id}}}
\newcommand{\rmrk}{{\rm {rank}}}

\newcommand{\funcI}{{\mathrm{Map}}^\fkI}

\newcommand{\rmSpan}{{\rm {Span}}}
\newcommand{\rmMax}{{\rm {Max}}}
\newcommand{\rmHom}{{\rm {Hom}}}
\newcommand{\rmCh}{{\rm {Ch}}}

\newcommand{\bZPi}{{\bZ\Pi}}
\newcommand{\bZgeqoPi}{{\bZgeqo\Pi}}

\newcommand{\tU}{{\tilde U}}

\newcommand{\rkN}{N}
\newcommand{\fkI}{I}

\newcommand{\bV}{{\mathbb{V}}}
\newcommand{\orderedPi}{{[\Pi]}}

\newcommand{\mcR}{{\mathcal{R}}}

\newcommand{\rmOb}{{\mathrm{Ob}}}

\newcommand{\mcH}{{\mathcal{H}}}

\newcommand{\mcW}{{\mathcal{W}}}
\newcommand{\mcWmcC}{{\mcW({\mathcal{C}})}}
\newcommand{\mcWmcR}{{\mcW({\mathcal{R}})}}

\newcommand{\preatilell}{{\tilde {\mathcal{L}}}_a}

\newcommand{\tilell}{{\tilde {\ell}}}
\newcommand{\atilell}{{\tilell_a}}

\newcommand{\hatN}{{\hat{N}}}
\newcommand{\hateta}{{\hat{\eta}}}
\newcommand{\hatal}{{\hat{\al}}}

\newcommand{\brre}{\stbsis}

\newcommand{\br}{{\hat r}}
\newcommand{\bal}{{\bar \alpha}}
\newcommand{\barN}{\hatN}
\newcommand{\barR}{{\bar R}}
\newcommand{\breta}{{\bar \eta}}
\newcommand{\barPi}{{\bar \Pi}}
\newcommand{\orderbarPi}{{(\barPi)}}
\newcommand{\parity}{\theta}
\newcommand{\barparity}{{\bar{\parity}}}

\newcommand{\perm}{\wp}

\newcommand{\mfkg}{{\mathfrak{g}}}
\newcommand{\barH}{{\bar{H}}}
\newcommand{\barE}{{\bar{E}}}
\newcommand{\barF}{{\bar{F}}}

\newcommand{\bbeta}{{\bar{\beta}}}

\newcommand{\superAnewtr}{{\breve \newtr}}
\newcommand{\superAc}{{\breve c}}
\newcommand{\superAC}{{\breve C}}
\newcommand{\superAR}{{\breve R}}
\newcommand{\superAmcC}{{\breve{\mcC}}}
\newcommand{\superAmcR}{{\breve{\mcR}}}

\newcommand{\superARp}{\superAR^+}

\newcommand{\superBnewtr}{{\dot \newtr}}
\newcommand{\superBc}{{\dot c}}
\newcommand{\superBC}{{\dot C}}
\newcommand{\superBR}{{\dot R}}
\newcommand{\superBmcC}{{\dot{\mcC}}}
\newcommand{\superBmcR}{{\dot{\mcR}}}
\newcommand{\superBf}{{\dot f}}
\newcommand{\superBp}{{\dot p}}

\newcommand{\superBRp}{\superBR^+}

\newcommand{\truebRrkN}{{\bR^\rkN}}

\newcommand{\ospefkI}{{\ddot  \fkI}}

\newcommand{\ospemfD}{{\ddot {\mathcal{A}}}}

\newcommand{\ospep}{{\ddot p}}

\newcommand{\ospeR}{{\ddot R}}

\newcommand{\ospeal}{{\ddot \al}}

\newcommand{\ospexi}{{\ddot \xi}}

\newcommand{\ospeeta}{{\ddot \eta}}

\newcommand{\ospec}{{\ddot c}}
\newcommand{\ospeC}{{\ddot C}}

\newcommand{\ospenewtr}{{\ddot \newtr}}

\newcommand{\ospemcC}{{\ddot \mcC}}
\newcommand{\ospemcR}{{\ddot \mcR}}

\newcommand{\ospef}{{\ddot f}}

\newcommand{\superordDbase}{{({\ddot \Pi})}}

\newcommand{\ospepCN}{{{\acute p}_\rkN}}

\newcommand{\superordCbase}{{({\acute \Pi})}}

\newcommand{\ospefCN}{{{\acute f}_\rkN}}

\newcommand{\bhm}{\chi}

\newcommand{\bigchiPi}{{\mathcal X}_\rkN}

\newcommand{\bq}{{\hat q}}

\newcommand{\tK}{{\tilde K}}
\newcommand{\tL}{{\tilde L}}
\newcommand{\tE}{{\tilde E}}
\newcommand{\tF}{{\tilde F}}

\newcommand{\tOm}{{\tilde \Omega}}

\newcommand{\bhmop}{{\bhm^{\rm{op}}}}

\newcommand{\tUp}{{\tilde \Upsilon}}

\newcommand{\tEp}{\tE}
\newcommand{\tFp}{\tF}

\newcommand{\tY}{{\tilde Y}}

\newcommand{\cI}{\mathcal{I}}
\newcommand{\tcI}{{\tilde \cI}}
\newcommand{\cJ}{\mathcal{J}}
\newcommand{\tcJ}{{\tilde \cJ}}

\newcommand{\bigchiPiirr}{\bigchiPi^{{\rm{irr}}}}

\newcommand{\mclM}{{\mathcal{M}}}
\newcommand{\mclL}{{\mathcal{L}}}

\newcommand{\tv}{{\tilde v}}
\newcommand{\cN}{\mathcal{N}}

\newcommand{\kp}{\kappa}

\newcommand{\bhmp}{\bhm^\prime}

\newcommand{\bigchiPifin}{\bigchiPi^{{\rm{fin}}}}

\newcommand{\sbhm}{s^\bhm}

\newcommand{\newtr}{\tau} 
\newcommand{\ppbigchiPifin}{\bigchiPi^{\prime\prime,{\rm{fin}}}}

\newcommand{\Abhm}{{{\mathcal{G}}(\bhm)}}

\newcommand{\pbigchiPifin}{\bigchiPi^{\prime,{\rm{fin}}}}

\newcommand{\newtrAbhm}{\newtr^\Abhm}

\newcommand{\Cbhm}{{{\mathcal{C}}_\bhm}}

\newcommand{\Rbhm}{{{\mathcal{R}}_\bhm}}

\newcommand{\bhms}{{1^\bhm s}}
\newcommand{\bhmps}{{1^{\bhmp} s}}

\newcommand{\bhmpp}{\bhm^{\prime\prime}}
\newcommand{\lgst}{{1^\bhm w_0}}

\newcommand{\bartrincl}{\bar{\trincl}}
\newcommand{\bzeta}{{\hat \zeta}}

\newcommand{\hT}{{\hat T}}

\newcommand{\mbbS}{{\mathbb{S}}}

\newcommand{\trincl}{\xi}

\newcommand{\hatbeta}{{\hat{\beta}}}

\newcommand{\hatw}{{\hat w}}
\newcommand{\hats}{{\hat s}}

\newcommand{\matMNR}{{{{\mathrm{M}}_\hatN}(\bR)}}
\newcommand{\matGLNR}{{{{\mathrm{GL}}_\hatN}(\bR)}}

\newcommand{\vecRN}{{\bR^\hatN}}
\newcommand{\vece}{\stbsis} 

\newcommand{\rsystem}{{\hat R}} 
\newcommand{\rsysbase}{{\hat \Pi}} 
\newcommand{\orderrsysbase}{{(\rsysbase)}} 

\newcommand{\rsysp}{{\rsystem^+(\rsysbase)}} 

\newcommand{\hatiota}{{\hat{\trincl}}}

\newcommand{\hatWPi}{{{\hat W}(\rsysbase)}} 

\newcommand{\hatSPi}{{{\hat S}(\rsysbase)}} 

\newcommand{\hatlng}{{\hat \ell}} 

\newcommand{\whLng}{{\widehat{\mathfrak{L}}}}

\newcommand{\hatwo}{{\hatw_0}}

\newcommand{\hatb}{{\hat b}}

\newcommand{\dotXNCartan}{{{\dot {\mathcal{X}}}^{\mathrm{Cartan}}_\rkN}}
\newcommand{\dotXNSuper}{{{\dot {\mathcal{X}}}^{\mathrm{Super}}_\rkN}}
\newcommand{\dotXNExtra}{{{\dot {\mathcal{X}}}^{\mathrm{Extra}}_\rkN}}

\newcommand{\tC}{{\tilde{C}}}
\newcommand{\tc}{{\tilde{c}}}

\newcommand{\mcA}{{\mathcal{A}}}
\newcommand{\tvsigma}{{\tilde{\tau}}}

\newcommand{\mcC}{{\mathcal{C}}}

\newcommand{\tR}{{\tilde{R}}}
\newcommand{\tRp}{\tR^+}

\newcommand{\tils}{{\tilde{s}}}

\newcommand{\mcQ}{\mathcal{Q}}
\newcommand{\mcK}{\mathcal{K}}
\newcommand{\toolmap}{\nabla }
\newcommand{\tltoolmap}{{\tilde{\toolmap}}}

\newcommand{\figHWM}{
\setlength{\unitlength}{1mm}
\begin{picture}(70,25)(-7,0)

\put(39, 22){{\scriptsize{$\Lambda$}}}\put(47, 10){{\scriptsize{$\newtr_1\Lambda$}}}
\put(-5, 22){{\scriptsize{$\newtr_2\Lambda$}}} 
\put(-17, 10){{\scriptsize{$\newtr_1 \newtr_2 \Lambda$}}}

\put(39, -4){{\scriptsize{$\newtr_2 \newtr_1 \Lambda$}}}
\put(-10, -4){{\scriptsize{$\newtr_1\newtr_2\newtr_1\Lambda$}}}
\put(-10, -7){{\scriptsize{$=\newtr_2\newtr_1\newtr_2\Lambda$}}}

\put(0,20){\circle*{1}}\multiput(0,20)(10,0){5}{\circle*{1}}
\put(-5,10){\circle*{1}}\multiput(-5,10)(10,0){6}{\circle*{1}}
\multiput(0,0)(10,0){5}{\circle*{1}}
\put(-5,10){\line(1,2){5}}\put(45,10){\line(-1,2){5}}
\put(0,0){\line(-1,2){5}}\put(40,0){\line(1,2){5}}
\multiput(0,20)(10,0){4}{\line(1,0){10}}
\multiput(0,0)(10,0){4}{\line(1,0){10}}
\end{picture}}

\newcommand{\figAzero}{
\setlength{\unitlength}{1mm}
\begin{picture}(40,10)(-7,0)

\put(1, 2){{\scriptsize{$\al_1$}}}
\put(12, 2){{\scriptsize{$\al_2$}}}
\put(23, 2){{\scriptsize{$\al_3$}}}
\put(34, 2){{\scriptsize{$\al_4$}}}

\put(1, 9){{\scriptsize{$\bq^{-2}$}}}
\put(2,6){\circle{3}}

\put(7, 8){{\scriptsize{$\bq^2$}}}
\put(3.5,6){\line(1,0){8}}

\put(12, 9){{\scriptsize{$-1$}}}
\put(13,6){\circle{3}}

\put(17, 8){{\scriptsize{$\bq^{-2}$}}}
\put(14.5,6){\line(1,0){8}}

\put(23, 9){{\scriptsize{$\bq^2$}}}
\put(24,6){\circle{3}}

\put(28, 8){{\scriptsize{$\bq^{-2}$}}}
\put(25.5,6){\line(1,0){8}}

\put(34, 9){{\scriptsize{$\bq^2$}}}
\put(35,6){\circle{3}}

\end{picture}}

\newcommand{\figAone}{
\setlength{\unitlength}{1mm}
\begin{picture}(40,10)(-7,0)

\put(1, 2){{\scriptsize{$\al_1$}}}
\put(12, 2){{\scriptsize{$\al_2$}}}
\put(23, 2){{\scriptsize{$\al_3$}}}
\put(34, 2){{\scriptsize{$\al_4$}}}

\put(1, 9){{\scriptsize{$-1$}}}
\put(2,6){\circle{3}}

\put(6, 8){{\scriptsize{$\bq^{-2}$}}}
\put(3.5,6){\line(1,0){8}}

\put(12, 9){{\scriptsize{$-1$}}}
\put(13,6){\circle{3}}

\put(18, 8){{\scriptsize{$\bq^2$}}}
\put(14.5,6){\line(1,0){8}}

\put(23, 9){{\scriptsize{$-1$}}}
\put(24,6){\circle{3}}

\put(28, 8){{\scriptsize{$\bq^{-2}$}}}
\put(25.5,6){\line(1,0){8}}

\put(34, 9){{\scriptsize{$\bq^2$}}}
\put(35,6){\circle{3}}

\end{picture}}

\newcommand{\figAtwo}{
\setlength{\unitlength}{1mm}
\begin{picture}(40,10)(-7,0)

\put(1, 2){{\scriptsize{$\al_1$}}}
\put(12, 2){{\scriptsize{$\al_2$}}}
\put(23, 2){{\scriptsize{$\al_3$}}}
\put(34, 2){{\scriptsize{$\al_4$}}}

\put(1, 9){{\scriptsize{$-1$}}}
\put(2,6){\circle{3}}

\put(6, 8){{\scriptsize{$\bq^{-2}$}}}
\put(3.5,6){\line(1,0){8}}

\put(12, 9){{\scriptsize{$\bq^2$}}}
\put(13,6){\circle{3}}

\put(17, 8){{\scriptsize{$\bq^{-2}$}}}
\put(14.5,6){\line(1,0){8}}

\put(23, 9){{\scriptsize{$-1$}}}
\put(24,6){\circle{3}}

\put(29, 8){{\scriptsize{$\bq^2$}}}
\put(25.5,6){\line(1,0){8}}

\put(34, 9){{\scriptsize{$-1$}}}
\put(35,6){\circle{3}}

\end{picture}}

\newcommand{\figAthree}{
\setlength{\unitlength}{1mm}
\begin{picture}(40,10)(-7,0)

\put(1, 2){{\scriptsize{$\al_1$}}}
\put(12, 2){{\scriptsize{$\al_2$}}}
\put(23, 2){{\scriptsize{$\al_3$}}}
\put(34, 2){{\scriptsize{$\al_4$}}}

\put(1, 9){{\scriptsize{$-1$}}}
\put(2,6){\circle{3}}

\put(6, 8){{\scriptsize{$\bq^{-2}$}}}
\put(3.5,6){\line(1,0){8}}

\put(12, 9){{\scriptsize{$\bq^2$}}}
\put(13,6){\circle{3}}

\put(17, 8){{\scriptsize{$\bq^{-2}$}}}
\put(14.5,6){\line(1,0){8}}

\put(23, 9){{\scriptsize{$\bq^2$}}}
\put(24,6){\circle{3}}

\put(28, 8){{\scriptsize{$\bq^{-2}$}}}
\put(25.5,6){\line(1,0){8}}

\put(34, 9){{\scriptsize{$-1$}}}
\put(35,6){\circle{3}}

\end{picture}}

\newcommand{\figAfour}{
\setlength{\unitlength}{1mm}
\begin{picture}(40,10)(-7,0)

\put(1, 2){{\scriptsize{$\al_1$}}}
\put(12, 2){{\scriptsize{$\al_2$}}}
\put(23, 2){{\scriptsize{$\al_3$}}}
\put(34, 2){{\scriptsize{$\al_4$}}}

\put(1, 9){{\scriptsize{$-1$}}}
\put(2,6){\circle{3}}

\put(7, 8){{\scriptsize{$\bq^2$}}}
\put(3.5,6){\line(1,0){8}}

\put(12, 9){{\scriptsize{$-1$}}}
\put(13,6){\circle{3}}

\put(17, 8){{\scriptsize{$\bq^{-2}$}}}
\put(14.5,6){\line(1,0){8}}

\put(23, 9){{\scriptsize{$\bq^2$}}}
\put(24,6){\circle{3}}

\put(28, 8){{\scriptsize{$\bq^{-2}$}}}
\put(25.5,6){\line(1,0){8}}

\put(34, 9){{\scriptsize{$-1$}}}
\put(35,6){\circle{3}}

\end{picture}}

\newcommand{\figAfive}{
\setlength{\unitlength}{1mm}
\begin{picture}(40,10)(-7,0)

\put(1, 2){{\scriptsize{$\al_1$}}}
\put(12, 2){{\scriptsize{$\al_2$}}}
\put(23, 2){{\scriptsize{$\al_3$}}}
\put(34, 2){{\scriptsize{$\al_4$}}}

\put(1, 9){{\scriptsize{$\bq^2$}}}
\put(2,6){\circle{3}}

\put(6, 8){{\scriptsize{$\bq^{-2}$}}}
\put(3.5,6){\line(1,0){8}}

\put(12, 9){{\scriptsize{$-1$}}}
\put(13,6){\circle{3}}

\put(18, 8){{\scriptsize{$\bq^2$}}}
\put(14.5,6){\line(1,0){8}}

\put(23, 9){{\scriptsize{$-1$}}}
\put(24,6){\circle{3}}

\put(28, 8){{\scriptsize{$\bq^{-2}$}}}
\put(25.5,6){\line(1,0){8}}

\put(34, 9){{\scriptsize{$-1$}}}
\put(35,6){\circle{3}}

\end{picture}}

\newcommand{\figAsix}{
\setlength{\unitlength}{1mm}
\begin{picture}(40,10)(-7,0)

\put(1, 2){{\scriptsize{$\al_1$}}}
\put(12, 2){{\scriptsize{$\al_2$}}}
\put(23, 2){{\scriptsize{$\al_3$}}}
\put(34, 2){{\scriptsize{$\al_4$}}}

\put(1, 9){{\scriptsize{$\bq^2$}}}
\put(2,6){\circle{3}}

\put(6, 8){{\scriptsize{$\bq^{-2}$}}}
\put(3.5,6){\line(1,0){8}}

\put(12, 9){{\scriptsize{$\bq^2$}}}
\put(13,6){\circle{3}}

\put(17, 8){{\scriptsize{$\bq^{-2}$}}}
\put(14.5,6){\line(1,0){8}}

\put(23, 9){{\scriptsize{$-1$}}}
\put(24,6){\circle{3}}

\put(29, 8){{\scriptsize{$\bq^2$}}}
\put(25.5,6){\line(1,0){8}}

\put(34, 9){{\scriptsize{$\bq^{-2}$}}}
\put(35,6){\circle{3}}

\end{picture}}

\newcommand{\figureA}{
\setlength{\unitlength}{1mm}
\begin{picture}(180,55)(5,0)

\put(10,55){$\bhm_{f,t_2},\,(t_2\in\fkJ_{0,4})$}
\put(60,55){$\bhm_{f,5}$}\put(110,55){$\bhm_{f,6}$}

\put(0,40){\figAzero} \put(49,46){\line(1,0){8}}\put(52,47){{\scriptsize{$2$}}}
\put(50,40){\figAone} \put(99,46){\line(1,0){8}}\put(102,47){{\scriptsize{$3$}}}
\put(100,40){\figAtwo}

\put(145,33){\line(0,1){6}}\put(146,35){{\scriptsize{$4$}}}

\put(10,35){$\bhm_{f,9}$}
\put(60,35){$\bhm_{f,8}$}\put(110,35){$\bhm_{f,7}$}

\put(0,20){\figAfive} \put(49,26){\line(1,0){8}}\put(52,27){{\scriptsize{$2$}}}
\put(50,20){\figAfour} \put(99,26){\line(1,0){8}}\put(102,27){{\scriptsize{$1$}}}
\put(100,20){\figAthree}

\put(34,13){\line(0,1){6}}\put(35,15){{\scriptsize{$3$}}}

\put(10,15){$\bhm_{f,10}$}

\put(0,0){\figAsix}
\end{picture}}

\newcommand{\figBzero}{
\setlength{\unitlength}{1mm}
\begin{picture}(40,10)(-7,0)

\put(1, 2){{\scriptsize{$\al_1$}}}
\put(12, 2){{\scriptsize{$\al_2$}}}
\put(23, 2){{\scriptsize{$\al_3$}}}
\put(34, 2){{\scriptsize{$\al_4$}}}

\put(1, 9){{\scriptsize{$\bq^{-2}$}}}
\put(2,6){\circle{3}}

\put(7, 8){{\scriptsize{$\bq^2$}}}
\put(3.5,6){\line(1,0){8}}

\put(12, 9){{\scriptsize{$-1$}}}
\put(13,6){\circle{3}}

\put(17, 8){{\scriptsize{$\bq^{-2}$}}}
\put(14.5,6){\line(1,0){8}}

\put(23, 9){{\scriptsize{$\bq^2$}}}
\put(24,6){\circle{3}}

\put(28, 8){{\scriptsize{$\bq^{-2}$}}}
\put(25.5,6){\line(1,0){8}}

\put(34, 9){{\scriptsize{$\bq$}}}
\put(35,6){\circle{3}}

\end{picture}}

\newcommand{\figBone}{
\setlength{\unitlength}{1mm}
\begin{picture}(40,10)(-7,0)

\put(1, 2){{\scriptsize{$\al_1$}}}
\put(12, 2){{\scriptsize{$\al_2$}}}
\put(23, 2){{\scriptsize{$\al_3$}}}
\put(34, 2){{\scriptsize{$\al_4$}}}

\put(1, 9){{\scriptsize{$-1$}}}
\put(2,6){\circle{3}}

\put(6, 8){{\scriptsize{$\bq^{-2}$}}}
\put(3.5,6){\line(1,0){8}}

\put(12, 9){{\scriptsize{$-1$}}}
\put(13,6){\circle{3}}

\put(18, 8){{\scriptsize{$\bq^2$}}}
\put(14.5,6){\line(1,0){8}}

\put(23, 9){{\scriptsize{$-1$}}}
\put(24,6){\circle{3}}

\put(28, 8){{\scriptsize{$\bq^{-2}$}}}
\put(25.5,6){\line(1,0){8}}

\put(34, 9){{\scriptsize{$\bq$}}}
\put(35,6){\circle{3}}

\end{picture}}

\newcommand{\figBtwo}{
\setlength{\unitlength}{1mm}
\begin{picture}(40,10)(-7,0)

\put(1, 2){{\scriptsize{$\al_1$}}}
\put(12, 2){{\scriptsize{$\al_2$}}}
\put(23, 2){{\scriptsize{$\al_3$}}}
\put(34, 2){{\scriptsize{$\al_4$}}}

\put(1, 9){{\scriptsize{$-1$}}}
\put(2,6){\circle{3}}

\put(6, 8){{\scriptsize{$\bq^{-2}$}}}
\put(3.5,6){\line(1,0){8}}

\put(12, 9){{\scriptsize{$\bq^2$}}}
\put(13,6){\circle{3}}

\put(17, 8){{\scriptsize{$\bq^{-2}$}}}
\put(14.5,6){\line(1,0){8}}

\put(23, 9){{\scriptsize{$-1$}}}
\put(24,6){\circle{3}}

\put(29, 8){{\scriptsize{$\bq^2$}}}
\put(25.5,6){\line(1,0){8}}

\put(34, 9){{\scriptsize{$-\bq^{-1}$}}}
\put(35,6){\circle{3}}

\end{picture}}

\newcommand{\figureB}{
\setlength{\unitlength}{1mm}
\begin{picture}(180,15)(5,0)

\put(10,15){\scriptsize{$\bhm_{f,u_1}\,(u_1\in\fkJ_{0,5}\cup\{10,11,16\})$}}
\put(60,15){\scriptsize{$\bhm_{f,u_2}\,(u_2\in\{6,9,12,15\})$}}
\put(110,15){\scriptsize{$\bhm_{f,u_3}\,(u_3\in\{7,8,13,14\})$}}

\put(0,00){\figBzero} \put(49,06){\line(1,0){8}}\put(52,07){{\scriptsize{$2$}}}
\put(50,00){\figBone} \put(99,06){\line(1,0){8}}\put(102,07){{\scriptsize{$3$}}}
\put(100,00){\figBtwo}

\end{picture}}

\newcommand{\figCzero}{
\setlength{\unitlength}{1mm}
\begin{picture}(40,10)(-7,0)

\put(1, 2){{\scriptsize{$\al_1$}}}
\put(12, 2){{\scriptsize{$\al_2$}}}
\put(23, 2){{\scriptsize{$\al_3$}}}
\put(34, 2){{\scriptsize{$\al_4$}}}

\put(1, 9){{\scriptsize{$-1$}}}
\put(2,6){\circle{3}}

\put(6, 8){{\scriptsize{$\bq^{-2}$}}}
\put(3.5,6){\line(1,0){8}}

\put(12, 9){{\scriptsize{$\bq^2$}}}
\put(13,6){\circle{3}}

\put(17, 8){{\scriptsize{$\bq^{-2}$}}}
\put(14.5,6){\line(1,0){8}}

\put(23, 9){{\scriptsize{$\bq^2$}}}
\put(24,6){\circle{3}}

\put(28, 8){{\scriptsize{$\bq^{-4}$}}}
\put(25.5,6){\line(1,0){8}}

\put(34, 9){{\scriptsize{$\bq^4$}}}
\put(35,6){\circle{3}}

\end{picture}}

\newcommand{\figCone}{
\setlength{\unitlength}{1mm}
\begin{picture}(40,10)(-7,0)

\put(1, 2){{\scriptsize{$\al_1$}}}
\put(12, 2){{\scriptsize{$\al_2$}}}
\put(23, 2){{\scriptsize{$\al_3$}}}
\put(34, 2){{\scriptsize{$\al_4$}}}

\put(1, 9){{\scriptsize{$-1$}}}
\put(2,6){\circle{3}}

\put(7, 8){{\scriptsize{$\bq^2$}}}
\put(3.5,6){\line(1,0){8}}

\put(12, 9){{\scriptsize{$-1$}}}
\put(13,6){\circle{3}}

\put(17, 8){{\scriptsize{$\bq^{-2}$}}}
\put(14.5,6){\line(1,0){8}}

\put(23, 9){{\scriptsize{$\bq^2$}}}
\put(24,6){\circle{3}}

\put(28, 8){{\scriptsize{$\bq^{-4}$}}}
\put(25.5,6){\line(1,0){8}}

\put(34, 9){{\scriptsize{$\bq^4$}}}
\put(35,6){\circle{3}}

\end{picture}}

\newcommand{\figCtwo}{
\setlength{\unitlength}{1mm}
\begin{picture}(40,10)(-7,0)

\put(1, 2){{\scriptsize{$\al_1$}}}
\put(12, 2){{\scriptsize{$\al_2$}}}
\put(23, 2){{\scriptsize{$\al_3$}}}
\put(34, 2){{\scriptsize{$\al_4$}}}

\put(1, 9){{\scriptsize{$\bq^2$}}}
\put(2,6){\circle{3}}

\put(6, 8){{\scriptsize{$\bq^{-2}$}}}
\put(3.5,6){\line(1,0){8}}

\put(12, 9){{\scriptsize{$-1$}}}
\put(13,6){\circle{3}}

\put(18, 8){{\scriptsize{$\bq^2$}}}
\put(14.5,6){\line(1,0){8}}

\put(23, 9){{\scriptsize{$-1$}}}
\put(24,6){\circle{3}}

\put(28, 8){{\scriptsize{$\bq^{-4}$}}}
\put(25.5,6){\line(1,0){8}}

\put(34, 9){{\scriptsize{$\bq^4$}}}
\put(35,6){\circle{3}}

\end{picture}}

\newcommand{\figCthree}{
\setlength{\unitlength}{1mm}
\begin{picture}(27,10)(-7,0)

\put(1, 2){{\scriptsize{$\al_1$}}}
\put(12, 2){{\scriptsize{$\al_2$}}}
\put(27, 11.5){{\scriptsize{$\al_3$}}}
\put(27, -0.5){{\scriptsize{$\al_4$}}}

\put(1, 9){{\scriptsize{$\bq^2$}}}
\put(2,6){\circle{3}}

\put(6, 8){{\scriptsize{$\bq^{-2}$}}}
\put(3.5,6){\line(1,0){8}}

\put(12, 9){{\scriptsize{$\bq^2$}}}
\put(13,6){\circle{3}}

\put(16, -1.5){{\scriptsize{$\bq^{-2}$}}}
\put(16, 10){{\scriptsize{$\bq^{-2}$}}}
\put(14,5){\line(2,-1){8.5}}
\put(14,7){\line(2,1){8.5}}

\put(23, -4.5){{\scriptsize{$-1$}}}
\put(23, 15){{\scriptsize{$-1$}}}
\put(24,0){\circle{3}}

\put(25.5, 5.5){{\scriptsize{$\bq^4$}}}
\put(24,1.5){\line(0,1){9}}

\put(24,12){\circle{3}}

\end{picture}}

\newcommand{\figCfour}{
\setlength{\unitlength}{1mm}
\begin{picture}(40,10)(-7,0)

\put(1, 2){{\scriptsize{$\al_1$}}}
\put(12, 2){{\scriptsize{$\al_2$}}}
\put(23, 2){{\scriptsize{$\al_4$}}}
\put(34, 2){{\scriptsize{$\al_3$}}}

\put(1, 9){{\scriptsize{$\bq^2$}}}
\put(2,6){\circle{3}}

\put(6, 8){{\scriptsize{$\bq^{-2}$}}}
\put(3.5,6){\line(1,0){8}}

\put(12, 9){{\scriptsize{$-1$}}}
\put(13,6){\circle{3}}

\put(18, 8){{\scriptsize{$\bq^2$}}}
\put(14.5,6){\line(1,0){8}}

\put(23, 9){{\scriptsize{$-1$}}}
\put(24,6){\circle{3}}

\put(28, 8){{\scriptsize{$\bq^{-4}$}}}
\put(25.5,6){\line(1,0){8}}

\put(34, 9){{\scriptsize{$\bq^4$}}}
\put(35,6){\circle{3}}

\end{picture}}

\newcommand{\figCfive}{
\setlength{\unitlength}{1mm}
\begin{picture}(40,10)(-7,0)

\put(1, 2){{\scriptsize{$\al_1$}}}
\put(12, 2){{\scriptsize{$\al_2$}}}
\put(23, 2){{\scriptsize{$\al_4$}}}
\put(34, 2){{\scriptsize{$\al_3$}}}

\put(1, 9){{\scriptsize{$-1$}}}
\put(2,6){\circle{3}}

\put(7, 8){{\scriptsize{$\bq^2$}}}
\put(3.5,6){\line(1,0){8}}

\put(12, 9){{\scriptsize{$-1$}}}
\put(13,6){\circle{3}}

\put(17, 8){{\scriptsize{$\bq^{-2}$}}}
\put(14.5,6){\line(1,0){8}}

\put(23, 9){{\scriptsize{$\bq^2$}}}
\put(24,6){\circle{3}}

\put(28, 8){{\scriptsize{$\bq^{-4}$}}}
\put(25.5,6){\line(1,0){8}}

\put(34, 9){{\scriptsize{$\bq^4$}}}
\put(35,6){\circle{3}}

\end{picture}}

\newcommand{\figCsix}{
\setlength{\unitlength}{1mm}
\begin{picture}(40,10)(-7,0)

\put(1, 2){{\scriptsize{$\al_1$}}}
\put(12, 2){{\scriptsize{$\al_2$}}}
\put(23, 2){{\scriptsize{$\al_4$}}}
\put(34, 2){{\scriptsize{$\al_3$}}}

\put(1, 9){{\scriptsize{$-1$}}}
\put(2,6){\circle{3}}

\put(6, 8){{\scriptsize{$\bq^{-2}$}}}
\put(3.5,6){\line(1,0){8}}

\put(12, 9){{\scriptsize{$\bq^2$}}}
\put(13,6){\circle{3}}

\put(17, 8){{\scriptsize{$\bq^{-2}$}}}
\put(14.5,6){\line(1,0){8}}

\put(23, 9){{\scriptsize{$\bq^2$}}}
\put(24,6){\circle{3}}

\put(28, 8){{\scriptsize{$\bq^{-4}$}}}
\put(25.5,6){\line(1,0){8}}

\put(34, 9){{\scriptsize{$\bq^4$}}}
\put(35,6){\circle{3}}

\end{picture}}

\newcommand{\figureC}{
\setlength{\unitlength}{1mm}
\begin{picture}(180,70)(5,0)

\put(10,65){$\bhm_{f,t_2}\,(t_2\in\fkJ_{0,9})$}
\put(60,65){$\bhm_{f,10}$}\put(110,65){$\bhm_{f,11}$}

\put(0,50){\figCzero} \put(49,56){\line(1,0){8}}\put(52,57){{\scriptsize{$1$}}}
\put(50,50){\figCone} \put(99,56){\line(1,0){8}}\put(102,57){{\scriptsize{$2$}}}
\put(100,50){\figCtwo}

\put(134,43){\line(0,1){7}}\put(135,46){{\scriptsize{$3$}}}
\put(134,12){\line(0,1){7}}\put(135,15){{\scriptsize{$4$}}}

\put(60,15){$\bhm_{f,14}$}
\put(110,40){$\bhm_{f,12}$}
\put(110,15){$\bhm_{f,13}$}

\put(49,06){\line(1,0){8}}\put(52,07){{\scriptsize{$1$}}}
\put(99,06){\line(1,0){8}}\put(102,07){{\scriptsize{$2$}}}
\put(100,25){\figCthree}

\put(10,15){$\bhm_{f,15}$}

\put(0,0){\figCsix}\put(50,0){\figCfive}\put(100,0){\figCfour}
\end{picture}}

\newcommand{\figDzero}{
\setlength{\unitlength}{1mm}
\begin{picture}(40,10)(-7,0)

\put(1, 2){{\scriptsize{$\al_1$}}}
\put(12, 2){{\scriptsize{$\al_2$}}}
\put(23, 2){{\scriptsize{$\al_3$}}}
\put(38, 11.5){{\scriptsize{$\al_4$}}}
\put(38, -0.5){{\scriptsize{$\al_5$}}}

\put(1, 9){{\scriptsize{$\bq^{-2}$}}}
\put(2,6){\circle{3}}

\put(7, 8){{\scriptsize{$\bq^2$}}}
\put(3.5,6){\line(1,0){8}}

\put(12, 9){{\scriptsize{$-1$}}}
\put(13,6){\circle{3}}

\put(17, 8){{\scriptsize{$\bq^{-2}$}}}
\put(14.5,6){\line(1,0){8}}

\put(23, 9){{\scriptsize{$\bq^2$}}}
\put(24,6){\circle{3}}

\put(27, -1.5){{\scriptsize{$\bq^{-2}$}}}
\put(27, 10){{\scriptsize{$\bq^{-2}$}}}
\put(25,5){\line(2,-1){8.5}}
\put(25,7){\line(2,1){8.5}}

\put(34, -4.5){{\scriptsize{$\bq^2$}}}
\put(34, 15){{\scriptsize{$\bq^2$}}}
\put(35,0){\circle{3}}

\put(35,12){\circle{3}}

\end{picture}}

\newcommand{\figDone}{
\setlength{\unitlength}{1mm}
\begin{picture}(40,10)(-7,0)

\put(1, 2){{\scriptsize{$\al_1$}}}
\put(12, 2){{\scriptsize{$\al_2$}}}
\put(23, 2){{\scriptsize{$\al_3$}}}
\put(38, 11.5){{\scriptsize{$\al_4$}}}
\put(38, -0.5){{\scriptsize{$\al_5$}}}

\put(1, 9){{\scriptsize{$-1$}}}
\put(2,6){\circle{3}}

\put(6, 8){{\scriptsize{$\bq^{-2}$}}}
\put(3.5,6){\line(1,0){8}}

\put(12, 9){{\scriptsize{$-1$}}}
\put(13,6){\circle{3}}

\put(18, 8){{\scriptsize{$\bq^2$}}}
\put(14.5,6){\line(1,0){8}}

\put(23, 9){{\scriptsize{$-1$}}}
\put(24,6){\circle{3}}

\put(27, -1.5){{\scriptsize{$\bq^{-2}$}}}
\put(27, 10){{\scriptsize{$\bq^{-2}$}}}
\put(25,5){\line(2,-1){8.5}}
\put(25,7){\line(2,1){8.5}}

\put(34, -4.5){{\scriptsize{$\bq^2$}}}
\put(34, 15){{\scriptsize{$\bq^2$}}}
\put(35,0){\circle{3}}

\put(35,12){\circle{3}}

\end{picture}}

\newcommand{\figDtwo}{
\setlength{\unitlength}{1mm}
\begin{picture}(40,10)(-7,0)

\put(1, 2){{\scriptsize{$\al_1$}}}
\put(12, 2){{\scriptsize{$\al_2$}}}
\put(23, 2){{\scriptsize{$\al_3$}}}
\put(38, 11.5){{\scriptsize{$\al_4$}}}
\put(38, -0.5){{\scriptsize{$\al_5$}}}

\put(1, 9){{\scriptsize{$-1$}}}
\put(2,6){\circle{3}}

\put(6, 8){{\scriptsize{$\bq^{-2}$}}}
\put(3.5,6){\line(1,0){8}}

\put(12, 9){{\scriptsize{$\bq^2$}}}
\put(13,6){\circle{3}}

\put(17, 8){{\scriptsize{$\bq^{-2}$}}}
\put(14.5,6){\line(1,0){8}}

\put(23, 9){{\scriptsize{$-1$}}}
\put(24,6){\circle{3}}

\put(27.5, -1){{\scriptsize{$\bq^2$}}}
\put(27.5, 11){{\scriptsize{$\bq^2$}}}
\put(25,5){\line(2,-1){8.5}}
\put(25,7){\line(2,1){8.5}}

\put(34, -4.5){{\scriptsize{$-1$}}}
\put(34, 15){{\scriptsize{$-1$}}}
\put(35,0){\circle{3}}

\put(36.5, 5.5){{\scriptsize{$\bq^{-4}$}}}
\put(35,1.5){\line(0,1){9}}

\put(35,12){\circle{3}}

\end{picture}}

\newcommand{\figDthree}{
\setlength{\unitlength}{1mm}
\begin{picture}(40,10)(-7,0)

\put(-10, 2){{\scriptsize{$\al_1$}}}
\put(1, 2){{\scriptsize{$\al_2$}}}
\put(12, 2){{\scriptsize{$\al_3$}}}
\put(23, 2){{\scriptsize{$\al_4$}}}
\put(34, 2){{\scriptsize{$\al_5$}}}

\put(-10, 9){{\scriptsize{$-1$}}}
\put(-9,6){\circle{3}}

\put(-5, 8){{\scriptsize{$\bq^{-2}$}}}
\put(-7.5,6){\line(1,0){8}}

\put(1, 9){{\scriptsize{$\bq^2$}}}
\put(2,6){\circle{3}}

\put(6, 8){{\scriptsize{$\bq^{-2}$}}}
\put(3.5,6){\line(1,0){8}}

\put(12, 9){{\scriptsize{$\bq^2$}}}
\put(13,6){\circle{3}}

\put(17, 8){{\scriptsize{$\bq^{-2}$}}}
\put(14.5,6){\line(1,0){8}}

\put(23, 9){{\scriptsize{$-1$}}}
\put(24,6){\circle{3}}

\put(29, 8){{\scriptsize{$\bq^4$}}}
\put(25.5,6){\line(1,0){8}}

\put(34, 9){{\scriptsize{$\bq^{-4}$}}}
\put(35,6){\circle{3}}

\end{picture}}

\newcommand{\figureD}{
\setlength{\unitlength}{1mm}
\begin{picture}(135,50)(-3,0)

\put(10,49){$\bhm_{f,u_1}$}
\put(10,45){\scriptsize{$(u_1\in\fkJ_{0,7}\cup\{14,15,22\})$}}

\put(80,49){$\bhm_{f,u_2}$}
\put(80,45){\scriptsize{$(u_2\in\{8,13,16,21\})$}}

\put(80,19){$\bhm_{f,u_3}$}
\put(80,15){\scriptsize{$(u_3\in\{9,12,17,20\})$}}

\put(10,19){$\bhm_{f,u_4}$}
\put(10,15){\scriptsize{$(u_4\in\{10,11,18,19\})$}}

\put(60,36){\line(1,0){8}}\put(63,37){{\scriptsize{$2$}}}
 \put(65,6){\line(1,0){8}}\put(68,7){{\scriptsize{$4$}}}

\put(104,19){\line(0,1){7}}\put(105,22){{\scriptsize{$3$}}}

\put(0,30){\figDzero}\put(70,30){\figDone}
\put(11,0){\figDthree}\put(70,0){\figDtwo}

\end{picture}}

\newcommand{\figFzero}{
\setlength{\unitlength}{1mm}
\begin{picture}(40,10)(-7,0)

\put(1, 2){{\scriptsize{$\al_1$}}}
\put(12, 2){{\scriptsize{$\al_2$}}}
\put(23, 2){{\scriptsize{$\al_3$}}}
\put(34, 2){{\scriptsize{$\al_4$}}}

\put(1, 9){{\scriptsize{$-1$}}}
\put(2,6){\circle{3}}

\put(6, 8){{\scriptsize{$\bq^{-2}$}}}
\put(3.5,6){\line(1,0){8}}

\put(12, 9){{\scriptsize{$\bq^2$}}}
\put(13,6){\circle{3}}

\put(17, 8){{\scriptsize{$\bq^{-4}$}}}
\put(14.5,6){\line(1,0){8}}

\put(23, 9){{\scriptsize{$\bq^4$}}}
\put(24,6){\circle{3}}

\put(28, 8){{\scriptsize{$\bq^{-4}$}}}
\put(25.5,6){\line(1,0){8}}

\put(34, 9){{\scriptsize{$\bq^4$}}}
\put(35,6){\circle{3}}

\end{picture}}

\newcommand{\figFone}{
\setlength{\unitlength}{1mm}
\begin{picture}(40,10)(-7,0)

\put(1, 2){{\scriptsize{$\al_1$}}}
\put(12, 2){{\scriptsize{$\al_2$}}}
\put(23, 2){{\scriptsize{$\al_3$}}}
\put(34, 2){{\scriptsize{$\al_4$}}}

\put(1, 9){{\scriptsize{$-1$}}}
\put(2,6){\circle{3}}

\put(7, 8){{\scriptsize{$\bq^2$}}}
\put(3.5,6){\line(1,0){8}}

\put(12, 9){{\scriptsize{$-1$}}}
\put(13,6){\circle{3}}

\put(17, 8){{\scriptsize{$\bq^{-4}$}}}
\put(14.5,6){\line(1,0){8}}

\put(23, 9){{\scriptsize{$\bq^4$}}}
\put(24,6){\circle{3}}

\put(28, 8){{\scriptsize{$\bq^{-4}$}}}
\put(25.5,6){\line(1,0){8}}

\put(34, 9){{\scriptsize{$\bq^4$}}}
\put(35,6){\circle{3}}

\end{picture}}

\newcommand{\figFtwo}{
\setlength{\unitlength}{1mm}
\begin{picture}(40,10)(-7,0)

\put(1, 2){{\scriptsize{$\al_4$}}}
\put(12, 2){{\scriptsize{$\al_3$}}}
\put(27, 11.5){{\scriptsize{$\al_1$}}}
\put(27, -0.5){{\scriptsize{$\al_2$}}}

\put(1, 9){{\scriptsize{$\bq^4$}}}
\put(2,6){\circle{3}}

\put(6, 8){{\scriptsize{$\bq^{-4}$}}}
\put(3.5,6){\line(1,0){8}}

\put(12, 9){{\scriptsize{$-1$}}}
\put(13,6){\circle{3}}

\put(17, -1.5){{\scriptsize{$\bq^4$}}}
\put(17, 11){{\scriptsize{$\bq^{-2}$}}}
\put(14,5){\line(2,-1){8.5}}
\put(14,7){\line(2,1){8.5}}

\put(23, -4.5){{\scriptsize{$-1$}}}
\put(23, 15){{\scriptsize{$\bq^2$}}}
\put(24,0){\circle{3}}

\put(25.5, 5.5){{\scriptsize{$\bq^{-2}$}}}
\put(24,1.5){\line(0,1){9}}

\put(24,12){\circle{3}}

\end{picture}}

\newcommand{\figFthree}{
\setlength{\unitlength}{1mm}
\begin{picture}(27,10)(-7,0)

\put(1, 2){{\scriptsize{$\al_2$}}}
\put(12, 2){{\scriptsize{$\al_3$}}}
\put(27, 11.5){{\scriptsize{$\al_1$}}}
\put(27, -0.5){{\scriptsize{$\al_4$}}}

\put(1, 9){{\scriptsize{$\bq^4$}}}
\put(2,6){\circle{3}}

\put(6, 8){{\scriptsize{$\bq^{-4}$}}}
\put(3.5,6){\line(1,0){8}}

\put(12, 9){{\scriptsize{$-1$}}}
\put(13,6){\circle{3}}

\put(17, -1.5){{\scriptsize{$\bq^4$}}}
\put(17, 11){{\scriptsize{$\bq^2$}}}
\put(14,5){\line(2,-1){8.5}}
\put(14,7){\line(2,1){8.5}}

\put(23, -4.5){{\scriptsize{$-1$}}}
\put(23, 15){{\scriptsize{$-1$}}}
\put(24,0){\circle{3}}

\put(25.5, 5.5){{\scriptsize{$\bq^{-6}$}}}
\put(24,1.5){\line(0,1){9}}

\put(24,12){\circle{3}}

\end{picture}}

\newcommand{\figFfour}{
\setlength{\unitlength}{1mm}
\begin{picture}(40,10)(-7,0)

\put(1, 2){{\scriptsize{$\al_2$}}}
\put(12, 2){{\scriptsize{$\al_3$}}}
\put(23, 2){{\scriptsize{$\al_4$}}}
\put(34, 2){{\scriptsize{$\al_1$}}}

\put(1, 9){{\scriptsize{$\bq^4$}}}
\put(2,6){\circle{3}}

\put(6, 8){{\scriptsize{$\bq^{-4}$}}}
\put(3.5,6){\line(1,0){8}}

\put(12, 9){{\scriptsize{$\bq^4$}}}
\put(13,6){\circle{3}}

\put(17, 8){{\scriptsize{$\bq^{-4}$}}}
\put(14.5,6){\line(1,0){8}}

\put(23, 9){{\scriptsize{$-1$}}}
\put(24,6){\circle{3}}

\put(29, 8){{\scriptsize{$\bq^6$}}}
\put(25.5,6){\line(1,0){8}}

\put(34, 9){{\scriptsize{$\bq^{-6}$}}}
\put(35,6){\circle{3}}

\end{picture}}

\newcommand{\figureF}{
\setlength{\unitlength}{1mm}
\begin{picture}(160,45)(5,0)

\put(10,44){$\bhm_{f,u_1}$}
\put(10,40){\scriptsize{$(u_1\in\fkJ_{0,9}\cup\{18\})$}}

\put(60,44){$\bhm_{f,u_2}$}
\put(60,40){\scriptsize{$(u_2\in\{10,17\})$}}

\put(110,44){$\bhm_{f,u_3}$}
\put(110,40){\scriptsize{$(u_3\in\{11,16\})$}}

\put(110,19){$\bhm_{f,u_4}$}
\put(110,15){\scriptsize{$(u_4\in\{12,15\})$}}

\put(60,19){$\bhm_{f,u_5}$}
\put(60,15){\scriptsize{$(u_5\in\fkJ_{13,14})$}}

\put(49,31){\line(1,0){8}}\put(52,32){{\scriptsize{$1$}}}
\put(99,31){\line(1,0){8}}\put(102,32){{\scriptsize{$2$}}}
\put(123,18){\line(0,1){7}}\put(124,21){{\scriptsize{$3$}}}
\put(99,6){\line(1,0){8}}\put(102,7){{\scriptsize{$4$}}}

\put(0,25){\figFzero}
\put(50,25){\figFone}
\put(100,25){\figFtwo}
\put(100,0){\figFthree}
\put(50,0){\figFfour}

\end{picture}}

\newcommand{\figGzero}{
\setlength{\unitlength}{1mm}
\begin{picture}(25,10)(-7,0)

\put(1, 2){{\scriptsize{$\al_1$}}}
\put(12, 2){{\scriptsize{$\al_2$}}}
\put(23, 2){{\scriptsize{$\al_3$}}}

\put(1, 9){{\scriptsize{$-1$}}}
\put(2,6){\circle{3}}

\put(6, 8){{\scriptsize{$\bq^{-2}$}}}
\put(3.5,6){\line(1,0){8}}

\put(12, 9){{\scriptsize{$\bq^2$}}}
\put(13,6){\circle{3}}

\put(17, 8){{\scriptsize{$\bq^{-6}$}}}
\put(14.5,6){\line(1,0){8}}

\put(23, 9){{\scriptsize{$\bq^6$}}}
\put(24,6){\circle{3}}

\end{picture}}

\newcommand{\figGone}{
\setlength{\unitlength}{1mm}
\begin{picture}(25,10)(-7,0)

\put(1, 2){{\scriptsize{$\al_1$}}}
\put(12, 2){{\scriptsize{$\al_2$}}}
\put(23, 2){{\scriptsize{$\al_3$}}}

\put(1, 9){{\scriptsize{$-1$}}}
\put(2,6){\circle{3}}

\put(7, 8){{\scriptsize{$\bq^2$}}}
\put(3.5,6){\line(1,0){8}}

\put(12, 9){{\scriptsize{$-1$}}}
\put(13,6){\circle{3}}

\put(17, 8){{\scriptsize{$\bq^{-6}$}}}
\put(14.5,6){\line(1,0){8}}

\put(23, 9){{\scriptsize{$\bq^6$}}}
\put(24,6){\circle{3}}

\end{picture}}

\newcommand{\figGtwo}{
\setlength{\unitlength}{1mm}
\begin{picture}(40,10)(-7,0)

\put(1, 2){{\scriptsize{$\al_2$}}}
\put(16, 11.5){{\scriptsize{$\al_1$}}}
\put(16, -0.5){{\scriptsize{$\al_3$}}}

\put(1, 9){{\scriptsize{$-1$}}}
\put(2,6){\circle{3}}

\put(6, -1.5){{\scriptsize{$\bq^6$}}}
\put(6, 11){{\scriptsize{$\bq^{-2}$}}}
\put(3,5){\line(2,-1){8.5}}
\put(3,7){\line(2,1){8.5}}

\put(12, -4.5){{\scriptsize{$-1$}}}
\put(12, 15){{\scriptsize{$\bq^2$}}}
\put(13,0){\circle{3}}

\put(14.5, 5.5){{\scriptsize{$\bq^{-4}$}}}
\put(13,1.5){\line(0,1){9}}

\put(13,12){\circle{3}}

\end{picture}}

\newcommand{\figGthree}{
\setlength{\unitlength}{1mm}
\begin{picture}(25,10)(-7,0)

\put(1, 2){{\scriptsize{$\al_2$}}}
\put(12, 2){{\scriptsize{$\al_3$}}}
\put(23, 2){{\scriptsize{$\al_1$}}}

\put(1, 9){{\scriptsize{$\bq^6$}}}
\put(2,6){\circle{3}}

\put(6, 8){{\scriptsize{$\bq^{-6}$}}}
\put(3.5,6){\line(1,0){8}}

\put(12, 9){{\scriptsize{$-1$}}}
\put(13,6){\circle{3}}

\put(18, 8){{\scriptsize{$\bq^4$}}}
\put(14.5,6){\line(1,0){8}}

\put(23, 9){{\scriptsize{$-\bq^{-2}$}}}
\put(24,6){\circle{3}}

\end{picture}}

\newcommand{\figureG}{
\setlength{\unitlength}{1mm}
\begin{picture}(180,25)(5,0)

\put(10,19){$\bhm_{f,u_1}$}
\put(10,15){\scriptsize{$(u_1\in\fkJ_{0,6}\cup\{13\})$}}

\put(49,19){$\bhm_{f,u_2}$}
\put(49,15){\scriptsize{$(u_2\in\{7,12\})$}}

\put(87,24){$\bhm_{f,u_3}$}
\put(87,20){\scriptsize{$(u_3\in\{8,11\})$}}

\put(121,19){$\bhm_{f,u_4}$}
\put(121,15){\scriptsize{$(u_4\fkJ_{9,10})$}}

\put(38,6){\line(1,0){8}}\put(41,7){{\scriptsize{$1$}}}
\put(77,6){\line(1,0){8}}\put(80,7){{\scriptsize{$2$}}}
\put(110,6){\line(1,0){8}}\put(113,7){{\scriptsize{$4$}}}

\put(0,0){\figGzero}
\put(39,0){\figGone}
\put(77,0){\figGtwo}
\put(111,0){\figGthree}

\end{picture}}

\newcommand{\figDXzero}{
\setlength{\unitlength}{1mm}
\begin{picture}(25,10)(-7,0)

\put(1, 2){{\scriptsize{$\al_1$}}}
\put(12, 2){{\scriptsize{$\al_2$}}}
\put(23, 2){{\scriptsize{$\al_3$}}}

\put(1, 9){{\scriptsize{$\bq$}}}
\put(2,6){\circle{3}}

\put(6, 8){{\scriptsize{$\bq^{-1}$}}}
\put(3.5,6){\line(1,0){8}}

\put(12, 9){{\scriptsize{$-1$}}}
\put(13,6){\circle{3}}

\put(17, 8){{\scriptsize{$\br^{-1}$}}}
\put(14.5,6){\line(1,0){8}}

\put(23, 9){{\scriptsize{$\br$}}}
\put(24,6){\circle{3}}

\end{picture}}

\newcommand{\figDXone}{
\setlength{\unitlength}{1mm}
\begin{picture}(25,10)(-7,0)

\put(1, 2){{\scriptsize{$\al_2$}}}
\put(16, 11.5){{\scriptsize{$\al_1$}}}
\put(16, -0.5){{\scriptsize{$\al_3$}}}

\put(1, 9){{\scriptsize{$-1$}}}
\put(2,6){\circle{3}}

\put(6.5, -0.5){{\scriptsize{$\br$}}}
\put(6.5, 11){{\scriptsize{$\bq$}}}
\put(3,5){\line(2,-1){8.5}}
\put(3,7){\line(2,1){8.5}}

\put(12, -4.5){{\scriptsize{$-1$}}}
\put(12, 15){{\scriptsize{$-1$}}}
\put(13,0){\circle{3}}

\put(14.5, 5.5){{\scriptsize{$(\bq\br)^{-1}$}}}
\put(13,1.5){\line(0,1){9}}

\put(13,12){\circle{3}}

\end{picture}}

\newcommand{\figDXtwo}{
\setlength{\unitlength}{1mm}
\begin{picture}(40,10)(-7,0)

\put(1, 2){{\scriptsize{$\al_2$}}}
\put(12, 2){{\scriptsize{$\al_1$}}}
\put(23, 2){{\scriptsize{$\al_3$}}}

\put(1, 9){{\scriptsize{$\bq$}}}
\put(2,6){\circle{3}}

\put(6, 8){{\scriptsize{$\bq^{-1}$}}}
\put(3.5,6){\line(1,0){8}}

\put(12, 9){{\scriptsize{$-1$}}}
\put(13,6){\circle{3}}

\put(17, 8){{\scriptsize{$\bq\br$}}}
\put(14.5,6){\line(1,0){8}}

\put(21, 9){{\scriptsize{$(\bq\br)^{-1}$}}}
\put(24,6){\circle{3}}

\end{picture}}

\newcommand{\figureDX}{
\setlength{\unitlength}{1mm}
\begin{picture}(110,25)(5,0)

\put(10,19){$\bhm_{f,u_1}$}
\put(10,15){\scriptsize{$(u_1\in\fkJ_{0,2}\cup\{7\})$}}

\put(49,24){$\bhm_{f,u_2}$}
\put(49,20){\scriptsize{$(u_2\in\{3,6\})$}}

\put(84,19){$\bhm_{f,u_3}$}
\put(84,15){\scriptsize{$(u_3\in\fkJ_{4,5})$}}

\put(38,6){\line(1,0){8}}\put(41,7){{\scriptsize{$2$}}}
\put(74,6){\line(1,0){8}}\put(77,7){{\scriptsize{$1$}}}

\put(0,0){\figDXzero}
\put(39,0){\figDXone}
\put(74,0){\figDXtwo}

\end{picture}}

\newcommand{\figZTzero}{
\setlength{\unitlength}{1mm}
\begin{picture}(25,10)(-7,0)

\put(1, 2){{\scriptsize{$\al_1$}}}
\put(12, 2){{\scriptsize{$\al_2$}}}

\put(1, 9){{\scriptsize{$\bzeta$}}}
\put(2,6){\circle{3}}

\put(6, 8){{\scriptsize{$\bq^{-1}$}}}
\put(3.5,6){\line(1,0){8}}

\put(12, 9){{\scriptsize{$\bq$}}}
\put(13,6){\circle{3}}

\end{picture}}

\newcommand{\figZTone}{
\setlength{\unitlength}{1mm}
\begin{picture}(25,10)(-7,0)

\put(1, 2){{\scriptsize{$\al_1$}}}
\put(12, 2){{\scriptsize{$\al_2$}}}

\put(1, 9){{\scriptsize{$\bzeta$}}}
\put(2,6){\circle{3}}

\put(5, 8){{\scriptsize{$\bq\bzeta^{-1}$}}}
\put(3.5,6){\line(1,0){8}}

\put(12, 9){{\scriptsize{$\bzeta\bq^{-1}$}}}
\put(13,6){\circle{3}}

\end{picture}}

\newcommand{\figureZT}{
\setlength{\unitlength}{1mm}
\begin{picture}(65,25)(5,0)

\put(10,19){$\bhm_{f,u_1}$}
\put(10,15){\scriptsize{$(u_1\in\fkJ_{0,1}\cup\{4\})$}}

\put(49,19){$\bhm_{f,u_2}$}
\put(49,15){\scriptsize{$(u_2\in\fkJ_{2,3})$}}

\put(34,6){\line(1,0){8}}\put(37,7){{\scriptsize{$1$}}}

\put(0,0){\figZTzero}
\put(39,0){\figZTone}

\end{picture}}

\newcommand{\figRFzero}{
\setlength{\unitlength}{1mm}
\begin{picture}(40,10)(-7,0)

\put(1, 2){{\scriptsize{$\al_1$}}}
\put(12, 2){{\scriptsize{$\al_2$}}}
\put(23, 2){{\scriptsize{$\al_3$}}}
\put(34, 2){{\scriptsize{$\al_4$}}}

\put(1, 9){{\scriptsize{$\bq$}}}
\put(2,6){\circle{3}}

\put(6, 8){{\scriptsize{$\bq^{-1}$}}}
\put(3.5,6){\line(1,0){8}}

\put(12, 9){{\scriptsize{$\bq$}}}
\put(13,6){\circle{3}}

\put(17, 8){{\scriptsize{$\bq^{-1}$}}}
\put(14.5,6){\line(1,0){8}}

\put(23, 9){{\scriptsize{$-1$}}}
\put(24,6){\circle{3}}

\put(28, 8){{\scriptsize{$-\bq$}}}
\put(25.5,6){\line(1,0){8}}

\put(34, 9){{\scriptsize{$-\bq^{-1}$}}}
\put(35,6){\circle{3}}

\end{picture}}

\newcommand{\figRFone}{
\setlength{\unitlength}{1mm}
\begin{picture}(40,10)(-7,0)

\put(1, 2){{\scriptsize{$\al_1$}}}
\put(12, 2){{\scriptsize{$\al_2$}}}
\put(27, 11.5){{\scriptsize{$\al_3$}}}
\put(27, -0.5){{\scriptsize{$\al_4$}}}

\put(1, 9){{\scriptsize{$\bq$}}}
\put(2,6){\circle{3}}

\put(6, 8){{\scriptsize{$\bq^{-1}$}}}
\put(3.5,6){\line(1,0){8}}

\put(12, 9){{\scriptsize{$-1$}}}
\put(13,6){\circle{3}}

\put(16, -0.5){{\scriptsize{$-1$}}}
\put(18, 11){{\scriptsize{$\bq$}}}
\put(14,5){\line(2,-1){8.5}}
\put(14,7){\line(2,1){8.5}}

\put(23, -4.5){{\scriptsize{$-1$}}}
\put(23, 15){{\scriptsize{$-1$}}}
\put(24,0){\circle{3}}

\put(25.5, 5.5){{\scriptsize{$-\bq^{-1}$}}}
\put(24,1.5){\line(0,1){9}}

\put(24,12){\circle{3}}

\end{picture}}

\newcommand{\figRFtwo}{
\setlength{\unitlength}{1mm}
\begin{picture}(40,10)(-7,0)

\put(1, 2){{\scriptsize{$\al_1$}}}
\put(12, 2){{\scriptsize{$\al_2$}}}
\put(23, 2){{\scriptsize{$\al_4$}}}
\put(34, 2){{\scriptsize{$\al_3$}}}

\put(1, 9){{\scriptsize{$\bq$}}}
\put(2,6){\circle{3}}

\put(6, 8){{\scriptsize{$\bq^{-1}$}}}
\put(3.5,6){\line(1,0){8}}

\put(12, 9){{\scriptsize{$-1$}}}
\put(13,6){\circle{3}}

\put(17, 8){{\scriptsize{$-1$}}}
\put(14.5,6){\line(1,0){8}}

\put(23, 9){{\scriptsize{$-1$}}}
\put(24,6){\circle{3}}

\put(28, 8){{\scriptsize{$-\bq$}}}
\put(25.5,6){\line(1,0){8}}

\put(34, 9){{\scriptsize{$-\bq^{-1}$}}}
\put(35,6){\circle{3}}

\end{picture}}

\newcommand{\figRFthree}{
\setlength{\unitlength}{1mm}
\begin{picture}(27,10)(-7,0)

\put(1, 2){{\scriptsize{$\al_3$}}}
\put(12, 2){{\scriptsize{$\al_4$}}}
\put(27, 11.5){{\scriptsize{$\al_1$}}}
\put(27, -0.5){{\scriptsize{$\al_2$}}}

\put(1, 9){{\scriptsize{$-\bq^{-1}$}}}
\put(2,6){\circle{3}}

\put(6, 8){{\scriptsize{$-\bq$}}}
\put(3.5,6){\line(1,0){8}}

\put(12, 9){{\scriptsize{$-1$}}}
\put(13,6){\circle{3}}

\put(16, -0.5){{\scriptsize{$-1$}}}
\put(17, 11){{\scriptsize{$\bq^{-1}$}}}
\put(14,5){\line(2,-1){8.5}}
\put(14,7){\line(2,1){8.5}}

\put(23, -4.5){{\scriptsize{$-1$}}}
\put(23, 15){{\scriptsize{$-1$}}}
\put(24,0){\circle{3}}

\put(25.5, 5.5){{\scriptsize{$\bq$}}}
\put(24,1.5){\line(0,1){9}}

\put(24,12){\circle{3}}

\end{picture}}

\newcommand{\figRFfour}{
\setlength{\unitlength}{1mm}
\begin{picture}(40,10)(-7,0)

\put(1, 2){{\scriptsize{$\al_2$}}}
\put(12, 2){{\scriptsize{$\al_1$}}}
\put(23, 2){{\scriptsize{$\al_4$}}}
\put(34, 2){{\scriptsize{$\al_3$}}}

\put(1, 9){{\scriptsize{$\bq$}}}
\put(2,6){\circle{3}}

\put(6, 8){{\scriptsize{$\bq^{-1}$}}}
\put(3.5,6){\line(1,0){8}}

\put(12, 9){{\scriptsize{$-1$}}}
\put(13,6){\circle{3}}

\put(17, 8){{\scriptsize{$-\bq$}}}
\put(14.5,6){\line(1,0){8}}

\put(23, 9){{\scriptsize{$-\bq^{-1}$}}}
\put(24,6){\circle{3}}

\put(28, 8){{\scriptsize{$-\bq$}}}
\put(25.5,6){\line(1,0){8}}

\put(34, 9){{\scriptsize{$-\bq^{-1}$}}}
\put(35,6){\circle{3}}

\end{picture}}

\newcommand{\figRFfive}{
\setlength{\unitlength}{1mm}
\begin{picture}(40,10)(-7,0)

\put(1, 2){{\scriptsize{$\al_3$}}}
\put(12, 2){{\scriptsize{$\al_2$}}}
\put(27, 11.5){{\scriptsize{$\al_1$}}}
\put(27, -0.5){{\scriptsize{$\al_4$}}}

\put(1, 9){{\scriptsize{$\bq$}}}
\put(2,6){\circle{3}}

\put(6, 8){{\scriptsize{$\bq^{-1}$}}}
\put(3.5,6){\line(1,0){8}}

\put(12, 9){{\scriptsize{$-1$}}}
\put(13,6){\circle{3}}

\put(16, -0.5){{\scriptsize{$-1$}}}
\put(18, 11){{\scriptsize{$\bq$}}}
\put(14,5){\line(2,-1){8.5}}
\put(14,7){\line(2,1){8.5}}

\put(23, -4.5){{\scriptsize{$-1$}}}
\put(23, 15){{\scriptsize{$-1$}}}
\put(24,0){\circle{3}}

\put(25.5, 5.5){{\scriptsize{$-\bq^{-1}$}}}
\put(24,1.5){\line(0,1){9}}

\put(24,12){\circle{3}}

\end{picture}}

\newcommand{\figRFsix}{
\setlength{\unitlength}{1mm}
\begin{picture}(40,10)(-7,0)

\put(1, 2){{\scriptsize{$\al_3$}}}
\put(12, 2){{\scriptsize{$\al_2$}}}
\put(23, 2){{\scriptsize{$\al_1$}}}
\put(34, 2){{\scriptsize{$\al_4$}}}

\put(1, 9){{\scriptsize{$\bq$}}}
\put(2,6){\circle{3}}

\put(6, 8){{\scriptsize{$\bq^{-1}$}}}
\put(3.5,6){\line(1,0){8}}

\put(12, 9){{\scriptsize{$\bq$}}}
\put(13,6){\circle{3}}

\put(17, 8){{\scriptsize{$\bq^{-1}$}}}
\put(14.5,6){\line(1,0){8}}

\put(23, 9){{\scriptsize{$-1$}}}
\put(24,6){\circle{3}}

\put(28, 8){{\scriptsize{$-\bq$}}}
\put(25.5,6){\line(1,0){8}}

\put(34, 9){{\scriptsize{$-\bq^{-1}$}}}
\put(35,6){\circle{3}}

\end{picture}}

\newcommand{\figureRF}{
\setlength{\unitlength}{1mm}
\begin{picture}(180,70)(5,0)

\put(10,69){$\bhm_{f,u_1}$}
\put(10,65){\scriptsize{$(u_1\in\fkJ_{0,4})$}}

\put(60,69){$\bhm_{f,u_2}$}
\put(60,65){\scriptsize{$(u_2\in\{5,13\})$}}

\put(110,69){$\bhm_{f,u_3}$}
\put(110,65){\scriptsize{$(u_3\in\{6,12\})$}}

\put(110,44){$\bhm_{f,u_4}$}
\put(110,40){\scriptsize{$(u_4\in\{7,11\})$}}

\put(110,19){$\bhm_{f,u_5}$}
\put(110,15){\scriptsize{$(u_5\in\fkJ_{8,10})$}}

\put(60,40){$\bhm_{f,14}$}

\put(60,15){$\bhm_{f,15}$}

\put(0,50){\figRFzero} \put(49,56){\line(1,0){8}}\put(52,57){{\scriptsize{$3$}}}
\put(50,50){\figRFone} \put(96,56){\line(1,0){8}}\put(99,57){{\scriptsize{$4$}}}
\put(100,50){\figRFtwo}

\put(73,43){\line(0,1){7}}\put(74,46){{\scriptsize{$2$}}}
\put(84,12){\line(0,1){7}}\put(85,15){{\scriptsize{$1$}}}
\put(123,43){\line(0,1){7}}\put(124,46){{\scriptsize{$2$}}}
\put(134,12){\line(0,1){7}}\put(135,15){{\scriptsize{$1$}}}

\put(100,25){\figRFthree}
\put(100,0){\figRFfour}

\put(50,25){\figRFfive}
\put(50,0){\figRFsix}
\end{picture}}

\newcommand{\figzero}{
\setlength{\unitlength}{1mm}
\begin{picture}(55,10)(-7,0)

\put(0, 1){{\footnotesize{$\al_1$}}}
\put(1, 9){{\scriptsize{$\bq$}}}
\put(2,6){\circle{3}}

\put(8, 8){{\scriptsize{$\bq^{-1}$}}}
\put(3.5,6){\line(1,0){12}}

\put(15, 1){{\footnotesize{$\al_2$}}}
\put(16, 9){{\scriptsize{$\bq$}}}
\put(17,6){\circle{3}}

\put(23, 8){{\scriptsize{$\bq^{-1}$}}}
\put(18.5,6){\line(1,0){12}}

\put(30, 1){{\footnotesize{$\al_3$}}}
\put(30, 9){{\scriptsize{$-1$}}}
\put(32,6){\circle{3}}

\put(38, 8){{\scriptsize{$-\bq$}}}
\put(33.5,6){\line(1,0){12}}

\put(45, 1){{\footnotesize{$\al_4$}}}
\put(45, 9){{\scriptsize{$-\bq^{-1}$}}}
\put(47,6){\circle{3}}

\end{picture}}

\newcommand{\figone}{
\setlength{\unitlength}{1mm}
\begin{picture}(55,17.5)(-7,0)

\put(0, 1){{\footnotesize{$\al_1$}}}
\put(1, 9){{\scriptsize{$\bq$}}}
\put(2,6){\circle{3}}

\put(8, 8){{\scriptsize{$\bq^{-1}$}}}
\put(3.5,6){\line(1,0){12}}

\put(15, 1){{\footnotesize{$\al_2$}}}
\put(15.5, 9){{\scriptsize{$-1$}}}
\put(17,6){\circle{3}}

\put(23, 3){{\scriptsize{$-1$}}}
\put(18.5, 6){\line(1,0){12}}

\put(30, 1){{\footnotesize{$\al_4$}}}
\put(28, 10){{\scriptsize{$-\bq^{-1}$}}}
\put(32, 6){\circle{3}}\put(34.5, 6){{\scriptsize{$-1$}}}

\put(18,7){\line(1,1){5.5}}\put(19.5, 11.5){{\scriptsize{$\bq$}}}
\put(31,7){\line(-1,1){5.5}}
\put(24.5, 13.5){\circle{3}}\put(27, 13.5){{\scriptsize{$-1$}}}
\put(23, 17){{\footnotesize{$\al_3$}}}

\end{picture}}

\newcommand{\figtwo}{
\setlength{\unitlength}{1mm}
\begin{picture}(55,17.5)(-7,0)

\put(0, 1){{\footnotesize{$\al_3$}}}
\put(1, 9){{\scriptsize{$\bq$}}}
\put(2,6){\circle{3}}

\put(8, 8){{\scriptsize{$\bq^{-1}$}}}
\put(3.5,6){\line(1,0){12}}

\put(15, 1){{\footnotesize{$\al_2$}}}
\put(15.5, 9){{\scriptsize{$-1$}}}
\put(17,6){\circle{3}}

\put(23, 3){{\scriptsize{$-1$}}}
\put(18.5, 6){\line(1,0){12}}

\put(30, 1){{\footnotesize{$\al_4$}}}
\put(28, 10){{\scriptsize{$-\bq^{-1}$}}}
\put(32, 6){\circle{3}}\put(34.5, 6){{\scriptsize{$-1$}}}

\put(18,7){\line(1,1){5.5}}\put(19.5, 11.5){{\scriptsize{$\bq$}}}
\put(31,7){\line(-1,1){5.5}}
\put(24.5, 13.5){\circle{3}}\put(27, 13.5){{\scriptsize{$-1$}}}
\put(23, 17){{\footnotesize{$\al_1$}}}

\end{picture}}

\newcommand{\figthree}{
\setlength{\unitlength}{1mm}
\begin{picture}(55,10)(-7,0)

\put(0, 1){{\footnotesize{$\al_1$}}}
\put(1, 9){{\scriptsize{$\bq$}}}
\put(2,6){\circle{3}}

\put(8, 8){{\scriptsize{$\bq^{-1}$}}}
\put(3.5,6){\line(1,0){12}}

\put(15, 1){{\footnotesize{$\al_2$}}}
\put(15, 9){{\scriptsize{$-1$}}}
\put(17,6){\circle{3}}

\put(22.5, 8){{\scriptsize{$-1$}}}
\put(18.5,6){\line(1,0){12}}

\put(30, 1){{\footnotesize{$\al_4$}}}
\put(30, 9){{\scriptsize{$-1$}}}
\put(32,6){\circle{3}}

\put(38, 8){{\scriptsize{$-\bq$}}}
\put(33.5,6){\line(1,0){12}}

\put(45, 1){{\footnotesize{$\al_3$}}}
\put(45, 9){{\scriptsize{$-\bq^{-1}$}}}
\put(47,6){\circle{3}}

\end{picture}}

\newcommand{\figfour}{
\setlength{\unitlength}{1mm}
\begin{picture}(55,10)(-7,0)

\put(0, 1){{\footnotesize{$\al_3$}}}
\put(1, 9){{\scriptsize{$\bq$}}}
\put(2,6){\circle{3}}

\put(8, 8){{\scriptsize{$\bq^{-1}$}}}
\put(3.5,6){\line(1,0){12}}

\put(15, 1){{\footnotesize{$\al_2$}}}
\put(16, 9){{\scriptsize{$\bq$}}}
\put(17,6){\circle{3}}

\put(23, 8){{\scriptsize{$\bq^{-1}$}}}
\put(18.5,6){\line(1,0){12}}

\put(30, 1){{\footnotesize{$\al_1$}}}
\put(30, 9){{\scriptsize{$-1$}}}
\put(32,6){\circle{3}}

\put(38, 8){{\scriptsize{$-\bq$}}}
\put(33.5,6){\line(1,0){12}}

\put(45, 1){{\footnotesize{$\al_4$}}}
\put(45, 9){{\scriptsize{$-\bq^{-1}$}}}
\put(47,6){\circle{3}}

\end{picture}}

\newcommand{\figfive}{
\setlength{\unitlength}{1mm}
\begin{picture}(55,10)(-7,0)

\put(0, 1){{\footnotesize{$\al_3$}}}
\put(1, 9){{\scriptsize{$\bq$}}}
\put(2,6){\circle{3}}

\put(8, 8){{\scriptsize{$\bq^{-1}$}}}
\put(3.5,6){\line(1,0){12}}

\put(15, 1){{\footnotesize{$\al_2$}}}
\put(15, 9){{\scriptsize{$-1$}}}
\put(17,6){\circle{3}}

\put(22.5, 8){{\scriptsize{$-1$}}}
\put(18.5,6){\line(1,0){12}}

\put(30, 1){{\footnotesize{$\al_4$}}}
\put(30, 9){{\scriptsize{$-1$}}}
\put(32,6){\circle{3}}

\put(38, 8){{\scriptsize{$-\bq$}}}
\put(33.5,6){\line(1,0){12}}

\put(45, 1){{\footnotesize{$\al_1$}}}
\put(45, 9){{\scriptsize{$-\bq^{-1}$}}}
\put(47,6){\circle{3}}

\end{picture}}

\newcommand{\figsix}{
\setlength{\unitlength}{1mm}
\begin{picture}(55,17.5)(-7,0)

\put(0, 1){{\footnotesize{$\al_2$}}}
\put(0.5, 9){{\scriptsize{$-1$}}}
\put(2,6){\circle{3}}

\put(8, 3){{\scriptsize{$-1$}}}
\put(3.5, 6){\line(1,0){12}}

\put(15, 1){{\footnotesize{$\al_4$}}}
\put(13, 10){{\scriptsize{$-\bq^{-1}$}}}
\put(17, 6){\circle{3}}\put(18, 8){{\scriptsize{$-1$}}}

\put(3,7){\line(1,1){5.5}}\put(4.5, 11.5){{\scriptsize{$\bq$}}}
\put(16,7){\line(-1,1){5.5}}
\put(9.5, 13.5){\circle{3}}\put(12, 13.5){{\scriptsize{$-1$}}}
\put(8, 17){{\footnotesize{$\al_3$}}}

\put(30, 1){{\footnotesize{$\al_1$}}}
\put(23, 8){{\scriptsize{$-\bq$}}}
\put(18.5,6){\line(1,0){12}}
\put(30, 9){{\scriptsize{$-\bq^{-1}$}}}
\put(32,6){\circle{3}}

\end{picture}}

\newcommand{\figseven}{
\setlength{\unitlength}{1mm}
\begin{picture}(55,17.5)(-7,0)

\put(0, 1){{\footnotesize{$\al_2$}}}
\put(0.5, 9){{\scriptsize{$-1$}}}
\put(2,6){\circle{3}}

\put(8, 3){{\scriptsize{$-1$}}}
\put(3.5, 6){\line(1,0){12}}

\put(15, 1){{\footnotesize{$\al_4$}}}
\put(13, 10){{\scriptsize{$-\bq^{-1}$}}}
\put(17, 6){\circle{3}}\put(18, 8){{\scriptsize{$-1$}}}

\put(3,7){\line(1,1){5.5}}\put(4.5, 11.5){{\scriptsize{$\bq$}}}
\put(16,7){\line(-1,1){5.5}}
\put(9.5, 13.5){\circle{3}}\put(12, 13.5){{\scriptsize{$-1$}}}
\put(8, 17){{\footnotesize{$\al_1$}}}

\put(30, 1){{\footnotesize{$\al_3$}}}
\put(23, 8){{\scriptsize{$-\bq$}}}
\put(18.5,6){\line(1,0){12}}
\put(30, 9){{\scriptsize{$-\bq^{-1}$}}}
\put(32,6){\circle{3}}

\end{picture}}

\newcommand{\figeight}{
\setlength{\unitlength}{1mm}
\begin{picture}(55,10)(-7,0)

\put(0, 1){{\footnotesize{$\al_2$}}}
\put(1, 9){{\scriptsize{$\bq$}}}
\put(2,6){\circle{3}}

\put(8, 8){{\scriptsize{$\bq^{-1}$}}}
\put(3.5,6){\line(1,0){12}}

\put(15, 1){{\footnotesize{$\al_1$}}}
\put(15, 9){{\scriptsize{$-1$}}}
\put(17,6){\circle{3}}

\put(23, 8){{\scriptsize{$-\bq$}}}
\put(18.5,6){\line(1,0){12}}

\put(30, 1){{\footnotesize{$\al_4$}}}
\put(30, 9){{\scriptsize{$-\bq^{-1}$}}}
\put(32,6){\circle{3}}

\put(38, 8){{\scriptsize{$-\bq$}}}
\put(33.5,6){\line(1,0){12}}

\put(45, 1){{\footnotesize{$\al_3$}}}
\put(45, 9){{\scriptsize{$-\bq^{-1}$}}}
\put(47,6){\circle{3}}

\end{picture}}

\newcommand{\fignine}{
\setlength{\unitlength}{1mm}
\begin{picture}(55,10)(-7,0)

\put(0, 1){{\footnotesize{$\al_2$}}}
\put(1, 9){{\scriptsize{$\bq$}}}
\put(2,6){\circle{3}}

\put(8, 8){{\scriptsize{$\bq^{-1}$}}}
\put(3.5,6){\line(1,0){12}}

\put(15, 1){{\footnotesize{$\al_3$}}}
\put(15, 9){{\scriptsize{$-1$}}}
\put(17,6){\circle{3}}

\put(23, 8){{\scriptsize{$-\bq$}}}
\put(18.5,6){\line(1,0){12}}

\put(30, 1){{\footnotesize{$\al_4$}}}
\put(30, 9){{\scriptsize{$-\bq^{-1}$}}}
\put(32,6){\circle{3}}

\put(38, 8){{\scriptsize{$-\bq$}}}
\put(33.5,6){\line(1,0){12}}

\put(45, 1){{\footnotesize{$\al_1$}}}
\put(45, 9){{\scriptsize{$-\bq^{-1}$}}}
\put(47,6){\circle{3}}

\end{picture}}

\newcommand{\figDallzero}{
\setlength{\unitlength}{1mm}
\begin{picture}(27,10)(-7,0)

\put(1, 2){{\scriptsize{$\al_1$}}}
\put(12, 2){{\scriptsize{$\al_2$}}}
\put(27, 11.5){{\scriptsize{$\al_3$}}}
\put(27, -0.5){{\scriptsize{$\al_4$}}}

\put(1, 9){{\scriptsize{$\bq^{-2}$}}}
\put(2,6){\circle{3}}

\put(7, 8){{\scriptsize{$\bq^2$}}}
\put(3.5,6){\line(1,0){8}}

\put(11, 9){{\scriptsize{$-1$}}}
\put(13,6){\circle{3}}

\put(16, -1.5){{\scriptsize{$\bq^{-2}$}}}
\put(16, 10){{\scriptsize{$\bq^{-2}$}}}
\put(14,5){\line(2,-1){8.5}}
\put(14,7){\line(2,1){8.5}}

\put(23, -5){{\scriptsize{$\bq^2$}}}
\put(23, 15){{\scriptsize{$\bq^2$}}}
\put(24,0){\circle{3}}

\put(24,12){\circle{3}}

\end{picture}}

\newcommand{\figDallone}{
\setlength{\unitlength}{1mm}
\begin{picture}(27,10)(-7,0)

\put(1, 2){{\scriptsize{$\al_1$}}}
\put(12, 2){{\scriptsize{$\al_2$}}}
\put(27, 11.5){{\scriptsize{$\al_3$}}}
\put(27, -0.5){{\scriptsize{$\al_4$}}}

\put(0, 9){{\scriptsize{$-1$}}}
\put(2,6){\circle{3}}

\put(6, 8){{\scriptsize{$\bq^2$}}}
\put(3.5,6){\line(1,0){8}}

\put(11, 9){{\scriptsize{$\bq^{-2}$}}}
\put(13,6){\circle{3}}

\put(17, -1.5){{\scriptsize{$\bq^2$}}}
\put(17, 11){{\scriptsize{$\bq^2$}}}
\put(14,5){\line(2,-1){8.5}}
\put(14,7){\line(2,1){8.5}}

\put(22, -5){{\scriptsize{$-1$}}}
\put(22, 15){{\scriptsize{$-1$}}}
\put(24,0){\circle{3}}

\put(25.5, 5.5){{\scriptsize{$\bq^{-4}$}}}
\put(24,1.5){\line(0,1){9}}

\put(24,12){\circle{3}}

\end{picture}}

\newcommand{\figDalltwo}{
\setlength{\unitlength}{1mm}
\begin{picture}(27,10)(-7,0)

\put(1, 2){{\scriptsize{$\al_1$}}}
\put(12, 2){{\scriptsize{$\al_2$}}}
\put(27, 11.5){{\scriptsize{$\al_3$}}}
\put(27, -0.5){{\scriptsize{$\al_4$}}}

\put(0, 9){{\scriptsize{$-1$}}}
\put(2,6){\circle{3}}

\put(5, 8){{\scriptsize{$\bq^{-2}$}}}
\put(3.5,6){\line(1,0){8}}

\put(11, 9){{\scriptsize{$-1$}}}
\put(13,6){\circle{3}}

\put(17, -1.5){{\scriptsize{$\bq^2$}}}
\put(17, 11){{\scriptsize{$\bq^2$}}}
\put(14,5){\line(2,-1){8.5}}
\put(14,7){\line(2,1){8.5}}

\put(22, -5){{\scriptsize{$-1$}}}
\put(22, 15){{\scriptsize{$-1$}}}
\put(24,0){\circle{3}}

\put(25.5, 5.5){{\scriptsize{$\bq^{-4}$}}}
\put(24,1.5){\line(0,1){9}}

\put(24,12){\circle{3}}

\end{picture}}

\newcommand{\figDallthree}{
\setlength{\unitlength}{1mm}
\begin{picture}(40,10)(-7,0)

\put(1, 2){{\scriptsize{$\al_1$}}}
\put(12, 2){{\scriptsize{$\al_2$}}}
\put(23, 2){{\scriptsize{$\al_3$}}}
\put(34, 2){{\scriptsize{$\al_4$}}}

\put(0, 9){{\scriptsize{$-1$}}}
\put(2,6){\circle{3}}

\put(7, 8){{\scriptsize{$\bq^2$}}}
\put(3.5,6){\line(1,0){8}}

\put(11, 9){{\scriptsize{$-1$}}}
\put(13,6){\circle{3}}

\put(16, 8){{\scriptsize{$\bq^{-2}$}}}
\put(14.5,6){\line(1,0){8}}

\put(22, 9){{\scriptsize{$-1$}}}
\put(24,6){\circle{3}}

\put(29, 8){{\scriptsize{$\bq^4$}}}
\put(25.5,6){\line(1,0){8}}

\put(34, 9){{\scriptsize{$\bq^{-4}$}}}
\put(35,6){\circle{3}}

\end{picture}}

\newcommand{\figDallfour}{
\setlength{\unitlength}{1mm}
\begin{picture}(40,10)(-7,0)

\put(1, 2){{\scriptsize{$\al_1$}}}
\put(12, 2){{\scriptsize{$\al_2$}}}
\put(23, 2){{\scriptsize{$\al_4$}}}
\put(34, 2){{\scriptsize{$\al_3$}}}

\put(1, 9){{\scriptsize{$-1$}}}
\put(2,6){\circle{3}}

\put(7, 8){{\scriptsize{$\bq^2$}}}
\put(3.5,6){\line(1,0){8}}

\put(11, 9){{\scriptsize{$-1$}}}
\put(13,6){\circle{3}}

\put(16, 8){{\scriptsize{$\bq^{-2}$}}}
\put(14.5,6){\line(1,0){8}}

\put(22, 9){{\scriptsize{$-1$}}}
\put(24,6){\circle{3}}

\put(29, 8){{\scriptsize{$\bq^4$}}}
\put(25.5,6){\line(1,0){8}}

\put(34, 9){{\scriptsize{$\bq^{-4}$}}}
\put(35,6){\circle{3}}

\end{picture}}

\newcommand{\figDallfive}{
\setlength{\unitlength}{1mm}
\begin{picture}(40,10)(-7,0)

\put(1, 2){{\scriptsize{$\al_1$}}}
\put(12, 2){{\scriptsize{$\al_2$}}}
\put(23, 2){{\scriptsize{$\al_3$}}}
\put(34, 2){{\scriptsize{$\al_4$}}}

\put(1, 9){{\scriptsize{$-1$}}}
\put(2,6){\circle{3}}

\put(6, 8){{\scriptsize{$\bq^{-2}$}}}
\put(3.5,6){\line(1,0){8}}

\put(12, 9){{\scriptsize{$\bq^2$}}}
\put(13,6){\circle{3}}

\put(16, 8){{\scriptsize{$\bq^{-2}$}}}
\put(14.5,6){\line(1,0){8}}

\put(22, 9){{\scriptsize{$-1$}}}
\put(24,6){\circle{3}}

\put(29, 8){{\scriptsize{$\bq^4$}}}
\put(25.5,6){\line(1,0){8}}

\put(34, 9){{\scriptsize{$\bq^{-4}$}}}
\put(35,6){\circle{3}}

\end{picture}}

\newcommand{\figDallsix}{
\setlength{\unitlength}{1mm}
\begin{picture}(40,10)(-7,0)

\put(1, 2){{\scriptsize{$\al_1$}}}
\put(12, 2){{\scriptsize{$\al_2$}}}
\put(23, 2){{\scriptsize{$\al_3$}}}
\put(34, 2){{\scriptsize{$\al_4$}}}

\put(1, 9){{\scriptsize{$\bq^2$}}}
\put(2,6){\circle{3}}

\put(5, 8){{\scriptsize{$\bq^{-2}$}}}
\put(3.5,6){\line(1,0){8}}

\put(11, 9){{\scriptsize{$-1$}}}
\put(13,6){\circle{3}}

\put(18, 8){{\scriptsize{$\bq^2$}}}
\put(14.5,6){\line(1,0){8}}

\put(23, 9){{\scriptsize{$\bq^{-2}$}}}
\put(24,6){\circle{3}}

\put(29, 8){{\scriptsize{$\bq^4$}}}
\put(25.5,6){\line(1,0){8}}

\put(34, 9){{\scriptsize{$\bq^{-4}$}}}
\put(35,6){\circle{3}}

\end{picture}}

\newcommand{\figDallseven}{
\setlength{\unitlength}{1mm}
\begin{picture}(40,10)(-7,0)

\put(1, 2){{\scriptsize{$\al_1$}}}
\put(12, 2){{\scriptsize{$\al_2$}}}
\put(23, 2){{\scriptsize{$\al_4$}}}
\put(34, 2){{\scriptsize{$\al_3$}}}

\put(1, 9){{\scriptsize{$\bq^2$}}}
\put(2,6){\circle{3}}

\put(5, 8){{\scriptsize{$\bq^{-2}$}}}
\put(3.5,6){\line(1,0){8}}

\put(11, 9){{\scriptsize{$-1$}}}
\put(13,6){\circle{3}}

\put(18, 8){{\scriptsize{$\bq^2$}}}
\put(14.5,6){\line(1,0){8}}

\put(23, 9){{\scriptsize{$\bq^{-2}$}}}
\put(24,6){\circle{3}}

\put(29, 8){{\scriptsize{$\bq^4$}}}
\put(25.5,6){\line(1,0){8}}

\put(34, 9){{\scriptsize{$\bq^{-4}$}}}
\put(35,6){\circle{3}}

\end{picture}}

\newcommand{\figDalleight}{
\setlength{\unitlength}{1mm}
\begin{picture}(40,10)(-7,0)

\put(1, 2){{\scriptsize{$\al_1$}}}
\put(12, 2){{\scriptsize{$\al_2$}}}
\put(23, 2){{\scriptsize{$\al_4$}}}
\put(34, 2){{\scriptsize{$\al_3$}}}

\put(1, 9){{\scriptsize{$-1$}}}
\put(2,6){\circle{3}}

\put(6, 8){{\scriptsize{$\bq^{-2}$}}}
\put(3.5,6){\line(1,0){8}}

\put(12, 9){{\scriptsize{$\bq^2$}}}
\put(13,6){\circle{3}}

\put(16, 8){{\scriptsize{$\bq^{-2}$}}}
\put(14.5,6){\line(1,0){8}}

\put(22, 9){{\scriptsize{$-1$}}}
\put(24,6){\circle{3}}

\put(29, 8){{\scriptsize{$\bq^4$}}}
\put(25.5,6){\line(1,0){8}}

\put(34, 9){{\scriptsize{$\bq^{-4}$}}}
\put(35,6){\circle{3}}

\end{picture}}

\newcommand{\figureDall}{
\setlength{\unitlength}{1mm}
\begin{picture}(180,100)(3,0)

\put(22,71){\line(0,1){7}}\put(23,74){{\scriptsize{$2$}}}
\put(72,71){\line(0,1){7}}\put(73,74){{\scriptsize{$3$}}}
\put(22,41){\line(0,1){7}}\put(23,44){{\scriptsize{$3$}}}
\put(72,41){\line(0,1){7}}\put(73,44){{\scriptsize{$2$}}}
\put(122,41){\line(0,1){7}}\put(123,44){{\scriptsize{$2$}}}

\put(49,68){\line(1,1){8}}\put(52,74){{\scriptsize{$1$}}}
\put(49,38){\line(1,1){8}}\put(52,44){{\scriptsize{$1$}}}

\put(99,76){\line(1,-1){8}}\put(103,74){{\scriptsize{$4$}}}

\put(10,51){\line(-1,-2){5}}\put(5,41){\line(0,-1){25}}\put(5,16){\line(1,-2){5}}
\put(2,28){{\scriptsize{$4$}}}

\put(110,51){\line(-1,-2){5}}\put(105,41){\line(0,-1){25}}
\put(105,16){\line(-1,-1){15}}\put(90,01){\line(-1,0){40}}
\put(97,12){{\scriptsize{$1$}}}

\put(10,97){$(1,1,0,0,0)$}
\put(60,97){$(0,1,1,0,0)$}

\put(0,80){\figDallzero} 
\put(50,80){\figDallone}

\put(10,67){$(1,0,1,0,0)$}
\put(60,65){$(0,1,0,1,0)$}\put(110,65){$(0,1,0,1,1)$}

\put(0,50){\figDalltwo} 
\put(50,50){\figDallthree} 
\put(100,50){\figDallfour}

\put(10,35){$(1,0,0,1,0)$}\put(60,35){$(0,0,1,1,0)$}\put(110,35){$(0,1,0,1,1)$}

\put(0,20){\figDallfive}\put(50,20){\figDallsix}\put(100,20){\figDallseven}

\put(10,10){$(1,0,0,1,1)$}

\put(0,-5){\figDalleight}
\end{picture}}

\begin{document}

\pagestyle{myheadings}
\title{
Classification of Finite Dimensional Irreducible Representations of
Generalized Quantum Groups via Weyl Groupoids
}

\author{Saeid Azam$^*$, Hiroyuki Yamane$^\dagger$, Malihe Yousofzadeh$^*$}
\markboth{\centerline{\em S. AZAM, H. YAMANE, M.
YOUSOFZADEH}}{\centerline{\em IRREDUCIBLE REPRESENTATIONS}}

 \maketitle
\let\thefootnote\relax\footnote{{\it Date:}}
\let\thefootnote\relax
\let\thefootnote\relax
\let\thefootnote\relax\footnote{$^*$This research was in part supported by a grant from IPM
(No. 90170217) and (No. 90170031). The authors would like to
thank the Center of Excellence for Mathematics, University of
Isfahan.}
\let\thefootnote\relax\footnote{$^\dagger$This research was also in part supported by
Japan's Grand-in-Aid for Scientific Research (C), 22540020.}

\begin{abstract}
Let $U(\bhm)$ be a generalized
quantum group such that $\dim U^+(\bhm)=\infty$,
$|R^+(\bhm)|<\infty$, and $R^+(\bhm)$ is irreducible, where
$U^+(\bhm)$ is the positive part of $U(\bhm)$, and $R^+(\bhm)$
is the Kharchenko's positive root system of $U^+(\bhm)$. In
this paper, we give a list of finite-dimensional irreducible
highest weight $U(\bhm)$-modules, relying on a {\it{special}}
reduced expression of the longest element of the Weyl groupoid
of $R(\bhm):=R^+(\bhm)\cup -R^+(\bhm)$.
\end{abstract}

\section*{{Introduction}}
In this paper, we give a list of finite-dimensional irreducible highest weight modules
of a generalized quantum group $U(\bhm)$ whose positive part
$U^+(\bhm)$ is infinite-dimensional and
has a Kharchenko's PBW-basis with a finite irreducible positive root system. We call such $U(\bhm)$
a generalized quantum group
of {\it{finite-and-infinite-dimensional-type}} (FID-type, for short).

We begin with recalling some facts of Lie superalgebras.
The class of
{\it{contragredient Lie
superalgebras}} \cite[Subsection~2.5.1]{Kac77} is defined in a way similar to that for Kac-Moody Lie algebras.
Kac classified the
finite-dimensional
simple Lie superalgebras
\cite[Theorem~5]{Kac77}, where
the finite-dimensional irreducible
contragredient Lie superalgebras played crucial roles; those are
\newline\par
(1) Simple Lie algebras of type
$X_\rkN$, where $X=A,\ldots,G$,

(2) $sl(m+1|n+1)$ ($m+n\geq 2$),

(3) $B(m,n)$ ($m\geq 0$, $n\geq 1$), $C(n)$
($n \geq 3$), $D(m,n)$ ($m\geq 2$, $n\geq 1$), $D(2,1;x)$
($x\ne 0,\,-1$), $F(4)$, $G(3)$.
\newline\newline
The ones in (1) and (3) are simple.
The simple Lie superalgebras $A(m,n)$ are defined by
$sl(m+1|n+1)$ if $m\ne n$, and otherwise
$A(n,n):=sl(n+1|n+1)/\mathfrak{i}$,
where $\mathfrak{i}$ is its unique one-dimensional
ideal.

Bases of the root systems of the Lie superalgebras of (2)-(3) are not conjugate under the action
of their Weyl groups. However each two of them are transformed to each other under the action of
their {\it{Weyl groupoids}} $W$,
whose axiomatic treatment was introduced by Heckenberger and the second author
\cite{HY08}.
Kac \cite[Theorem~8~(c)]{Kac77} gave a list
of finite-dimensional irreducible highest weight modules of the Lie superalgebras in (2)-(3) above.
In the same way as in this paper,
we can
have a new proof of recovering the list;
our idea is to use a specially good one among the reduced
expressions of the longest element (with a `standard' end domain)
of the Weyl groupoid $W$, 
see also Remark~\ref{remark:lasttwo}.

Let $\mathfrak{g}:=sl(m+1|n+1)$ or $C(n)$ for example.
Let $\mathfrak{h}$ be a Cartan subalgebra of $\mathfrak{g}$
such that the Dynkin diagram of $(\mathfrak{g},\mathfrak{h})$ is a standard one.
Let $\Pi=\{\al_i|1\leq i\leq \dim\mathfrak{h}\}$ be the set of simple roots $\al_i$
corresponding to $\mathfrak{h}$.
Let $w_0$ be
the longest element of $W$ of $\mathfrak{g}$
whose end domain is corresponding to $\mathfrak{h}$.
Then the length $\ell(w_0)$ of $w_0$ is equal to the number of positive roots
of $\mathfrak{g}$.
Let $k$ be the the number of even positive roots
of $\mathfrak{g}$.
The key fact used in this paper
is that there exists a reduced
expression $s_{i_1}\cdots s_{i_{\ell(w_0)}}$ of $w_0$ such that
$s_{i_1}\cdots s_{i_{x-1}}(\al_{i_x})$, $1\leq x\leq k$, are even positive roots,
and $s_{i_1}\cdots s_{i_{y-1}}(\al_{i_y})$, $k+1\leq y\leq \ell(w_0)$, are odd positive roots.
We claim that this is essential to the fact that
an irreducible highest weight $\mathfrak{g}$-module of  highest weight
$\Lambda$ is finite-dimensional if and only if
${\frac {2\langle \Lambda , \al_i\rangle} {\langle \al_i, \al_i \rangle}}\in\bZgeqo$
for all even simple roots $\al_i$, where $\langle\,,\,\rangle$ is the bilinear form
coming from the Killing form of $\mathfrak{g}$.

Motivated by Andruskiewitsch and Schneider's theory \cite{AS98}, \cite{AS10} toward the classification
of pointed Hopf algebras,
Heckenberger \cite{Hec09} classified
the Nichols algebras of diagonal-type.
Let $\bK$ be a characteristic zero field.
Let $U(\bhm)$ be the $\bK$-algebra 
defined in the same manner as in the Lusztig's book \cite[3.1.1 (a)-(e)]{b-Lusztig93}
for any bi-homomorphism $\bhm:\bZPi\times\bZPi\to\bKt$,
where $\Pi=\{\al_i|i\in\fkI\}$ is the set of
simple roots of the Kharchencko's positive root system $R^+(\bhm)$
associated with $\bhm$.
We call $U(\bhm)$
{\it{the generalized quantum group}}.
We say that $\bhm$ (or $U(\bhm)$)
is of {\it{finite-type}} if $R^+(\bhm)$ is finite and irreducible.
We say that $\bhm$ (or $U(\bhm)$)
is of {\it{finite-and-infinite-dimensional-type}} (FID-type, for short)
if $\bhm$ is of finite-type and $\dim U^+(\bhm)=\infty$.
A Nichols algebra of diagonal-type
is isomorphic to the positive part $U^+(\bhm)$ of $U(\bhm)$
for some $\bhm$ of finite-type.
If $U(\bhm)$ is of FID-type, then it is
a multi-parameter quantum algebra of
a simple Lie algebra in (1),
a multi-parameter quantum superalgebra of
a simple Lie superalgebra in (2) or (3),
or one of two algebras
\cite[Table~1-Row~5, Table~3-Row~14]{Hec09}.
Our main result, Theorems in Section~\label{section:MainTh}, gives a list of
finite-dimensional irreducible
highest weight modules of such $U(\bhm)$ in the way mentioned above.

Studying representation theory of $U(\bhm)$ must be interesting and fruitful
since
the factorization formula of Shapovalov determinants of
any $U(\bhm)$ of finite-type has been obtained by
Heckenberger and the second author in \cite{HY10}.
We believe that it would help us to give a new way to study
Lusztig's conjecture~\cite{Lusztig90}.

This paper is organized as follows.

In Section~1, we correct general facts of Weyl groupoids.
In Section~2, we give exmaples of reduced expressions of the longest elements
of the Weyl groups, which will be used in Section~3.
In Section~3, we give reduced expressions of the Weyl groupoids
associated to Lie superalgebras of $ABCD$-types.
In Section~4, we give a definition of
generalized quantum groups $U(\bhm)$ associated with any bi-homomorphism $\bhm$,
explain the Kharchenko's PBW-theorem of $U(\bhm)$,
and introduce the Weyl groupoids
associated with $\bhm$.
In Section~5, we explain properties of Weyl groupoids associated
with $U(\bhm)$ of finite-type, and Heckenberger's classification
of
$U(\bhm)$'s of FID-type.
In Section~6,
we give a key criterion when an irreducible highest $U(\bhm)$-module
is finite-dimensional, see Lemma~\ref{lemma:lgelmd}.
In Section~7, we give a list of finite-dimensional
irreducible highest weight $U(\bhm)$-modules
for $U(\bhm)$ having a standard Dynkin diagram, see Theorems there.

In \cite{Y06}, the second author has given a result similar to Theorem~\ref{theorem:MainSec}~(4).

\section{Weyl groupoids}\label{section:WGpoid}
\subsection{Basic terminology}\label{subsection:BTerm}

For a set $\mathfrak{s}$, let
$|\mathfrak{s}|$ denote the cardinality of $\mathfrak{s}$.

Let $\bN$ denote the set of positive integers.
Let $\bZ$ denote the ring of integers.
Let $\bQ$ denote the field of rational numbers.
For $x$, $y\in\bQ$, let
$\fkJ_{x,y}:=\{n\in\bZ|x\leq n\leq y\}$.
Note that $\fkJ_{x,y}$ is empty if $x>y$,
or $n<x\leq y<n+1$ for some $n\in\bZ$.
Moreover, if $\fkJ_{x,y}$ is empty, then as usual,
we let $\sum_{n \in \fkJ_{x,y}}z_n:=0$
and $\prod_{n \in \fkJ_{x,y}}z_n:=1$ with any symbol $z_n$,
where $0$ and $1$ are the zero element
and the unit element of $\bZ$ respectively.
For $x\in\bQ$, let $\fkJ_{x,\infty}:=\{m\in\bZ|m\geq x\}$,
and let $\fkJ_{-\infty,x}:=\{m\in\bZ|m\leq x\}$.
Then $\bN =\fkJ_{1,\infty}$. Let
$\bZgeqo :=\fkJ_{0,\infty}$. Let $\bR$ denote the field of real numbers.

Throughout this paper, we use the fixed notation below.
\begin{equation}\label{eqn:fixdfI}
\begin{array}{l}
\mbox{Let $\rkN\in\bN$ be a fixed positive integer. Let $\fkI:=\fkJ_{1,\rkN}$.} \\
\mbox{Let $\bV$ be a fixed $\rkN$-dimensional $\bR$-linear space.} \\
\mbox{Let $\orderedPi=(\al_i|i\in\fkI)$ be a fixed ordered $\bR$-basis of $\bV$.} \\
\mbox{Let $\Pi:=\{\al_i|i\in\fkI\}$, so $\Pi$ is a (set) $\bR$-basis of $\bV$.} \\
\mbox{Let $\bK$ be a fixed filed whose characteristic is zero.
Let $\bKt:=\bK\setminus\{0\}$.}
\end{array}
\end{equation}

Let $\bZPi:=\oplus_{i\in\fkI}\bZ\al_i(\subsetneq\bV)$, i.e.,
$\bZPi$ is the free $\bZ$-module
with the basis $\Pi$.
Then $\rmrk_\bZ\bZPi=\rkN$.
Let $\bZgeqoPi:=\oplus_{i\in\fkI}\bZgeqo\al_i(\subsetneq\bZPi)$.

For $n\in\bN\cup\{\infty\}$, let $\funcI_n$ be the set of maps from
$\fkJ_{1,n}$ to $\fkI$. Let $\funcI_0$ be the set composed of
a unique element $\phi$, i.e.,
$|\funcI_0|=1$ and $\phi\in\funcI_0$.

For a unital $\bK$-algebra $\mathcal{G}$, let $\rmCh(\mathcal{G})$
denote the set of $\bK$-algebra homomorphisms from $\mathcal{G}$ to $\bK$.

Let $t\in\bK$. For
$m\in\bZgeqo$, let $(m)_t:=\sum_{j\in \fkJ_{0,m-1}}t^j$ and
$(m)_t!:=\prod_{j\in \fkJ_{1,m}}(j)_t$.
For $m\in\bZgeqo$ and $n\in\fkJ_{1,m-1}$,
let ${m \choose 0}_t:={m \choose m}_t:=1$,
and ${m \choose n}_t:={m-1 \choose n-1}_t +t^n{m-1 \choose n}_t$.
Then ${m \choose n}_t(n)_t!(m-n)_t!=(m)_t!$, and
${m \choose n}_t:=t^{m-n}{m-1 \choose n-1}_t +t^n{m-1 \choose n}_t$.

For $m\in\bN$,
let $\bKt_m:=\{r\in\bKt|r^m=1,\,r^t\ne 1 (t\in \fkJ_{1,m-1})\}$.
Let $\bKtinf:=\bKt\setminus\cup_{m\in\bN}\bKt_m$.

For an associative $\bK$-algebra $\mathfrak{a}$
and $X$, $Y\in\mathfrak{a}$,
let $[X,Y]:=
XY-YX$.

Let $\uplus$ mean disjoint union of sets.

For $\bZ$-modules $\mathfrak{b}$ and $\mathfrak{c}$,
let $\rmHom_\bZ(\mathfrak{b},\mathfrak{c})$ be
the $\bZ$-module formed by the $\bZ$-module
homomorphisms from $\mathfrak{b}$ to $\mathfrak{c}$.

The symbols $\delta_{ij}$, $\delta_{i,j}$,
and $\delta(i,j)$, denote Kronecker's
delta, that is, $\delta_{ij}=1$ if $i=j$, and $\delta_{ij}=0$ otherwise.

\subsection{Modification of axioms of generalized root systems} 
\label{subsection:ModifiGRS}

Keep the notation as in \eqref{eqn:fixdfI}.

We call an $\rkN\times\rkN$ matrix 
$\tC=[\tc_{ij}]_{i,j\in\fkI}$ over $\bZ$
{\it{a generalized Cartan matrix}}
if $({\rm{M}}1)$ and $({\rm{M}}2)$ below hold.\newline\par
$({\rm{M}}1)$ $\tc_{ii}=2$ ($i\in\fkI$).  \par 
$({\rm{M}}2)$ $\tc_{jk}\leq 0$, $\delta(\tc_{jk},0)=\delta(\tc_{kj},0)$ 
($j$, $k\in\fkI$, $j\ne k$).
\newline\par
Let $\mcA$ be a non-empty set. 
Let $\tvsigma_i:\mcA\to\mcA$ be maps 
($i\in\fkI$).
Let $\tC^a=[\tc^a_{ij}]_{i,j\in\fkI}$ be generalized Cartan matrices
($a\in\mcA$). We call the data
\begin{equation*}
\mcC=\mcC(\fkI,\mcA,(\tvsigma_i)_{i\in\fkI},(\tC^a)_{a\in\mcA})
\end{equation*} {\it{a {\rm{(}}rank-$\rkN$) Cartan scheme{\rm{}}}} if
$({\rm{C}}1)$ and $({\rm{C}}2)$ below hold.
\newline\par
$({\rm{C}}1)$  $\tvsigma_i^2=\rmid_\mcA$ ($i\in\fkI$). \par
$({\rm{C}}2)$ $\tc^{\tvsigma_i(a)}_{ij}=\tc^a_{ij}$ ($i\in\fkI$).
\newline\par
Let $\mcC=\mcC(\fkI,\mcA,(\tvsigma_i)_{i\in\fkI},(\tC^a)_{a\in\mcA})$
be a Cartan scheme.
Define $\tils^a_i\in\rmGL(\bV)$ ($a\in\mcA$, $i\in\fkI$)
by 
\begin{equation}\label{eqn:deftilsa}
\tils^a_i(\al_j)=\al_j-\tc^a_{ij}\al_i\quad(j\in\fkI).
\end{equation} 
Then
\begin{equation}\label{eqn:tilsa-sq}
(\tils^a_i)^2=\tils^{\tvsigma_i(a)}_i\tils^a_i=\rmid_\bV
\quad(a\in\mcA,\,i\in\fkI).
\end{equation} 

\begin{notation}\label{notation:Ab-chift-1} Let $\mcC=\mcC(\fkI,\mcA,(\tvsigma_i)_{i\in\fkI},(\tC^a)_{a\in\mcA})$
be a Cartan scheme.

{\rm{(1)}}
For $a\in\mcA$ and $f\in\funcI_n$ for some $n\in\bZgeqo\cup\{\infty\}$,
let 
\begin{equation*}
\begin{array}{ll}
a_{f,0}:=a,\,\,1^a\tils_{f,0}:=\rmid_\bV,\\
a_{f,t}:=\tvsigma_{f(t)}(a_{f,t-1}),\,\
1^a\tils_{f,t}:=1^a\tils_{f,t-1}\tils^{a_t}_{f(t)}
\quad (t\in\fkJ_{1,n}).
\end{array}
\end{equation*}
{\rm{(2)}}
For $a$, $a^\prime\in\mcA$, let
\begin{equation}\label{eqn:AbrHmWbhm}
\mcH(a,a^\prime) 
:=\{\,1^a\tils_{f,t}\,|\,f\in\funcI_\infty,\,
t\in\bZgeqo,\,a_{f,t}=a^\prime\,\},
\end{equation} as a subset of $\rmGL(\bV)$.
\end{notation}

We say that a Cartan scheme $\mcC=\mcC(\fkI,\mcA,(\tvsigma_i)_{i\in\fkI},(\tC^a)_{a\in\mcA})$
is {\it{connected}} if
$|\mcH(a,a^\prime)|\geq 1$ for all $a$, $a^\prime\in\mcA$. 

\begin{definition}\label{definition:GnRtSys}
Let $\mcC=\mcC(\fkI,\mcA,(\tvsigma_i)_{i\in\fkI},(\tC^a)_{a\in\mcA})$
be a Cartan scheme.
For each $a\in\mcA$, let
$\tR(a)$ be a subset of $\bV=\oplus_{i\in\fkI}\bR\al_i$,
and $\tRp(a):=\tR(a)\cap(\oplus_{i\in\fkI}\bZgeqo\al_i)$.
We call the data
\begin{equation*}
\mcR=\mcR(\mcC,(\tR(a))_{a\in\mcA})
\end{equation*} {\it{a generalized root system of type $\mcC$}} if
$({\rm{R}}1)$-$({\rm{R}}4)$ below hold.
\newline\par
$({\rm{R}}1)$ \quad $\tR(a)=\tRp(a)\cup-\tRp(a)$ \quad {\rm{(}}$a\in\mcA${\rm{)}}.\par
$({\rm{R}}2)$ \quad $\tR(a)\cap\bZ\al_i=\{\,\al_i,\,-\al_i\,\}$ \quad 
{\rm{(}}$a\in\mcA$, $i\in\fkI${\rm{)}}.\par
$({\rm{R}}3)$ \quad $\tils^a_i(\tR(a))=\tR(\tvsigma_i(a))$ \quad 
{\rm{(}}$a\in\mcA$, $i\in\fkI${\rm{)}}. \par
$({\rm{R}}4)$ \quad For $a$, $a^\prime\in\mcA$,
if $\rmid_\bV\in\mcH(a,a^\prime)$, then $a=a^\prime$.
\end{definition}

Let $\mcR=\mcR(\mcC,(\tR(a))_{a\in\mcA})$ be
a generalized root system of type $\mcC$. By $({\rm{R}}1)$, 
$({\rm{R}}2)$, $({\rm{R}}3)$ and the definition of 
$\tils^a_i$, we have
\begin{equation}\label{eqn:tsaitRp}
\tils^a_i(\tRp(a)\setminus\{\al_i\})=\tRp(\tvsigma_i(a))\setminus\{\al_i\},
\end{equation} and
\begin{equation}\label{eqn:tsaitRp-2}
-\tc^a_{ij}=\max\{\,k\in\bZgeqo\,|\,\al_j+k\al_i\in\tRp(a)\,\}
\quad(i,j\in\fkI,\,i\ne j).
\end{equation}
If $\mcC$ is connected, we say that $\mcR$ is {\it{connected}}.

\begin{lemma}\label{lemma:preISOM}
Let $\mcR=\mcR(\mcC,(\tR(a))_{a\in\mcA})$ and
$\mcR^\prime=\mcR(\mcC^\prime,(\tR^\prime(a^\prime))_{a^\prime\in\mcA})$ be
generalized root systems of types 
$\mcC$ and $\mcC^\prime$ respectively.
Let ${\check{a}}\in\mcA$ and ${\check{a}}^\prime\in\mcA^\prime$.
Assume that $\tR({\check{a}})=\tR^\prime({\check{a}}^\prime)$.
Then we have
\begin{equation*}
\tR({\check{a}}_{f,n})=\tR^\prime({\check{a}}^\prime_{f,n}),\,\,
\tils^{{\check{a}}_{f,n}}_i=\tils^{{\check{a}}^\prime_{f,n}}_i
\quad(n\in\bN,\,f\in\funcI_n,\,i\in\fkI).
\end{equation*}
\end{lemma}
{\it{Proof.}} This lemma follows easily from
\eqref{eqn:tilsa-sq}, \eqref{eqn:tsaitRp-2}
and $({\rm{R}}3)$.
\hfill $\Box$

\begin{definition}\label{definition:ISOM}
{\rm{Let $\mcR=\mcR(\mcC,(\tR(a))_{a\in\mcA})$ and
$\mcR^\prime=\mcR(\mcC^\prime,(\tR^\prime(a^\prime))_{a^\prime\in\mcA^\prime})$ be
connected generalized root systems of types 
$\mcC$ and $\mcC^\prime$ respectively.
Let ${\check{a}}\in\mcA$ and ${\check{a}}^\prime\in\mcA^\prime$.

{\rm{(1)}} We say that the pair $(\mcR,{\check{a}})$
is {\it{quasi-isomorphic}} to the pair $(\mcR^\prime,{\check{a}}^\prime)$
if $\tR({\check{a}})=\tR^\prime({\check{a}}^\prime)$.

{\rm{(2)}} We say that the pair $(\mcR,{\check{a}})$
is {\it{isomorphic}} to the pair $(\mcR^\prime,{\check{a}}^\prime)$
if $\tR({\check{a}})=\tR^\prime({\check{a}}^\prime)$
and for all $n\in\bN$ and all $f\in\funcI_n$, it follows that 
${\check{a}}_{f,n}={\check{a}}$
if and only if ${\check{a}}^\prime_{f,n}={\check{a}}^\prime$.
}}
\end{definition}

\begin{lemma} \label{lemma:rp-of-ax}
Let $\mcC=\mcC(\fkI,\mcA,(\tvsigma_i)_{i\in\fkI},(\tC^a)_{a\in\mcA})$
be a Cartan scheme. Let $\mcR=\mcR(\mcC,(\tR(a))_{a\in\mcA})$ be
a generalized root system of type $\mcC$.
Let $a\in\mcA$ and $i$, $j\in\fkI$ with $i\ne j$. 
Let $m:=|\tRp(a)\cap(\bR\al_i\oplus\bR\al_j)|\in\fkJ_{2,\infty}\cup\{\,\infty\,\}$.
Assume $m<\infty$.
Define $f\in\funcI_{2m}$ by $f(2x-1):=i$ and $f(2x):=j$
$(x\in\fkJ_{1,m})$. Then  
\newline\par
$({\rm{R}}4)^\prime$ \quad $a_{f,2m}=a$ and $1^a\tils_{f,2m}=\rmid_\bV$.

\end{lemma}
{\it{Proof.}}
For $x\in\fkJ_{1,m}$, let $\beta_x:=1^a\tils_{f,x-1}(\al_{f(x)})$.
For $x\in\fkJ_{0,m}$, let 
$Z_x:=\tRp(a_{f,x})\cap(\bR\al_i\oplus\bR\al_j)$ and
$Y_x:=Z_0\cap-1^a\tils_{f,x}(Z_x)$.
By \eqref{eqn:deftilsa} and \eqref{eqn:tsaitRp}, 
\begin{equation}\label{eqn:preYpri}
|Z_x|=m\quad (x\in\fkJ_{0,m}).
\end{equation}
We show that for $x\in\fkJ_{1,m}$,
\newline\par
$(*)_x$ \quad $|Y_x|=x$ and $Y_x=\{\,\beta_y\,|\,y\in\fkJ_{1,x}\,\}$.
\newline\newline
Then $(*)_1$ follows from \eqref{eqn:tsaitRp}. Assume that $x\in\fkJ_{2,m}$ and 
$(*)_{x-1}$ holds. Then
\begin{equation}\label{eqn:Ypri}
\begin{array}{lcl}
Y_x & = & Z_0\cap(-1^a\tils_{f,x}(Z_x\setminus\{\,\al_{f(x)}\,\})
\cup\{\,-1^a\tils_{f,x}(\al_{f(x)})\,\}) \\
& = & Z_0\cap(-1^a\tils_{f,x-1}(Z_{x-1}\setminus\{\,\al_{f(x)}\,\})
\cup\{\,\beta_x\,\}) \\
& & \quad\mbox{(by \eqref{eqn:deftilsa} and \eqref{eqn:tsaitRp})} \\
& = & (Y_{x-1}\setminus\{\,-\beta_x\,\})\cup(Z_0\cap\{\,\beta_x\,\}).
\end{array}
\end{equation} 
Since $-1^a\tils_{f,x}(Z_x\setminus
\{\,\al_{f(x)}\,\})
\cup\{\,-1^a\tils_{f,x}(\al_{f(x)})\,\}=\emptyset$,
we have $(Y_{x-1}\setminus\{\,-\beta_x\,\})\cap(Z_0\cap\{\,\beta_x\,\})=\emptyset$.
Hence, by \eqref{eqn:Ypri},
\begin{equation}\label{eqn:Yprid}
Y_x=(Y_{x-1}\setminus\{\,-\beta_x\,\})\uplus(Z_0\cap\{\,\beta_x\,\}).
\end{equation}
Assume that $\beta_x\notin Z_0$. 
Then $\beta_x\in -Z_0$, so
$1^a\tils_{f,x-1}(\al_{f(x)})=\beta_x\in -Z_0$.   Since $\beta_{x-1}\in Z_0$, 
$1^a\tils_{f,x-1}(\al_{f(x-1)})=-1^a\tils_{f,x-2}(\al_{f(x-1)})= -\beta_{x-1}\in-Z_0$.
Since $\{\,f(x-1),\,f(x)\,\}=\{\,i,\,j\,\}$, we have 
$Y_{x-1}=Z_0$. Hence $x-1=m$, contradiction. So
$\beta_x\in Z_0$.
From \eqref{eqn:Yprid}, we obtain $(*)_x$.

By \eqref{eqn:preYpri} and $(*)_m$, $1^a\tils_{f,m}(Z_m)=-Z_0$.
Hence $1^a\tils_{f,m}(\{\,\al_i,\,\al_j\,\})=\{\,-\al_i,\,-\al_j\,\}$.
By the same argument, letting $f^\prime\in\funcI_{m}$
by $f^\prime(y):=f(m+y)$,
we have 
$1^{a_{f,m}}\tils_{f^\prime,m}(\{\,\al_i,\,\al_j\,\})=\{\,-\al_i,\,-\al_j\,\}$.
Hence $1^a\tils_{f,2m}(\{\,\al_i,\,\al_j\,\})=\{\,\al_i,\,\al_j\,\}$.
By \eqref{eqn:deftilsa}, 
the determinant of the $2\times2$-matrix $(\tils^{a_{f,x}}_{f(x)})_{|\bR\al_i\oplus\bR\al_j}$ is $-1$
for every $x\in\fkJ_{1,2m}$. So $1^a\tils_{f,2m}(\al_k)=\al_k$ for 
$k\in\{\,i,\,j\,\}$. By \eqref{eqn:deftilsa}, 
for 
$k\in\fkI\setminus\{\,i,\,j\,\}$, $1^a\tils_{f,2m}(\al_k)\in\al_k
+(\bZgeqo\al_i\oplus\bZgeqo\al_j)$. From $({\rm{R}}1)$,
we obtain the second claim of $({\rm{R}}4)^\prime$.
From (R4),
we obtain the first claim of $({\rm{R}}4)^\prime$.
This completes the proof.
\hfill $\Box$

\begin{remark} \label{remark:eqAXs}
Original definition of the generalized root systems was 
given in terms of $({\rm{R}}1)$-$({\rm{R}}3)$, $({\rm{R}}4)^\prime$,
see \cite{HY08}, \cite{CH09}.
From \cite[Lemma~8~(iii)]{HY08} and Lemma~\ref{lemma:rp-of-ax},
it follows that the definition due to $({\rm{R}}1)$-$({\rm{R}}5)$
is equivalent to the one due to
$({\rm{R}}1)$-$({\rm{R}}3)$, $({\rm{R}}4)^\prime$, $({\rm{R}}5)$.
\end{remark} 

\begin{definition} \label{definition:defWGpoid}
{\rm{Let $\mcC=\mcC(\fkI,\mcA,(\tvsigma_i)_{i\in\fkI},(\tC^a)_{a\in\mcA})$
be a Cartan scheme. 
Let $\mcWmcC$ be the category defined by (${\rm{cat}}1$)-(${\rm{cat}}3$) below. 
\newline\par
(${\rm{cat}}1$) The collection $\rmOb(\mcWmcC)$
of its objects is the set $\mcA$. \par
(${\rm{cat}}2$) For $a$,  
$a^\prime\in\mcA$, 
the set $\rmHom_\mcWmcC(a,a^\prime)$ of its morphisms from $a^\prime$ to $a$ is defined by 
\begin{equation*}
\rmHom_\mcWmcC(a,a^\prime)
:=\{\,(a,w,a^\prime)\,|\,w\in\mcH(a,a^\prime)\,\},
\end{equation*}
as
a subset of
$\mcA\times\rmGL(\bV)\times\mcA$.
\par
(${\rm{cat}}3$) For $a$,  
$a^\prime$, $a^{\prime\prime}$, 
the composition
\begin{equation*}
\rmHom_\mcWmcC(a,a^\prime)\times
\rmHom_\mcWmcC(a^\prime,a^{\prime\prime})
\to\rmHom_\mcWmcC(a,a^{\prime\prime})
\end{equation*}
of morphisms is defined by 
\begin{equation*}
(a,w,a^\prime)\circ(a^\prime,w^\prime,a^{\prime\prime})
:=(a,ww^\prime,a^{\prime\prime}),
\end{equation*}
where $ww^\prime$ means the product of elements $w$
and $w^\prime$ of the group $\rmGL(\bV)$.
We call $\mcWmcC$ {\it{a Weyl groupoid}} of $\mcC$.
If $\mcR$ is a generalized root system of type $\mcC$,
we let $\mcWmcR:=\mcWmcC$ and call it {\it{a Weyl groupoid}} of $\mcR$.
}}
\end{definition}

\subsection{Length function of a Weyl groupoid}\label{subsection:Lfunc}
In Subsection~\ref{subsection:Lfunc},
let $\mcC=\mcC(\fkI,\mcA,(\tvsigma_i)_{i\in\fkI},(\tC^a)_{a\in\mcA})$
and $\mcR=\mcR(\mcC,(\tR(a))_{a\in\mcA})$
be a Cartan scheme and
a generalized root system of type $\mcC$
respectively, and let $a\in\mcA$.

Let $\mcH(a,\_):=\cup_{a^\prime\in\mcA}\mcH(a,a^\prime)$.

For $w\in\mcH(a,\_)$, let
\begin{equation}\label{eqn:prellw}
\preatilell(w):=\{\,\beta\in \tRp(a)\,|
\,w^{-1}(\beta)\in -\bZgeqoPi\,\}.
\end{equation} 

Define the map $\atilell:\mcH(a,\_)\to\bZgeqo$ by 
\begin{equation}\label{eqn:defbhmell}
\atilell(w):=|\preatilell(w)|.
\end{equation}

\begin{lemma}\label{lemma:non-initial} 
{\rm{(1)}} Let $w\in\mcH(a,\_)$. Then the following three conditions
{\rm{(1-i)}}, {\rm{(1-ii)}} and {\rm{(1-iii)}} are equivalent.
\newline\par
{\rm{(1-i)}}\quad $\preatilell(w)=\tRp(a)$.
\par
{\rm{(1-ii)}}\quad $w^{-1}(\Pi)\subseteq -\bZgeqoPi\setminus\{0\}$, i.e.,
$\Pi\subseteq\preatilell(w)$.
\par
{\rm{(1-iii)}}\quad $w(\Pi)=-\Pi$.
\newline\par
{\rm{(2)}} For $a^\prime\in\mcA$ and $w\in\mcH(a,a^\prime)$,
we have 
\begin{equation*}
\preatilell(w)=\{\,\beta\in \tRp(a)\,|
\,w^{-1}(\beta)\in -\tRp(a^\prime)\,\}.
\end{equation*} 

{\rm{(3)}} Let $a^\prime$, $a^{\prime\prime}\in\mcA$, 
$w\in \mcH(a,a^\prime)$ and $w^\prime\in \mcH(a,a^{\prime\prime})$.
If $w=w^\prime$, then
$a^\prime=a^{\prime\prime}$. 

{\rm{(4)}} For $w\in\mcH(a,\_)$, we have
\begin{equation}\label{eqn:dfell-rnl}
\atilell(w)={\rm{Min}}\{l\in\bZgeqo\,|\,
\exists f\in\funcI_l,\,
1^a\tils_{f,l}=w\,\}.
\end{equation}

{\rm{(5)}} Let $a^\prime\in\mcA$. For $w\in\mcH(a,a^\prime)$
and $i\in\fkI$, we have
\begin{equation}\label{eqn:ellwc}
\atilell(w \tils^{\tvsigma_ia}_i)=\left\{\begin{array}{ll}
\atilell(w)+1 & \quad\mbox{if $w(\al_i)\in \tRp(a)$}, \\
\atilell(w)-1 & \quad\mbox{if $w(\al_i)\in -\tRp(a)$}.
\end{array}\right.
\end{equation}

{\rm{(6)}} Let $w\in\mcH(a,\_)$ and $l:=\atilell(w)$.
Let $f \in\funcI_l$ be such that $w=1^a\tils_{f,l}$. Then we have
\begin{equation}\label{eqn:ellwf}
\preatilell(w)=\{\,1^a\tils_{f,r-1}(\al_{f(r)})\,|\,r\in\fkJ_{1,l}\,\}.
\end{equation}
\end{lemma}
{\it{Proof.}} The claim (1) is clear from 
Definition~\ref{definition:GnRtSys}~$({\rm{R}}1)$, 
$({\rm{R}}2)$ and 
\eqref{eqn:prellw}. 
The claim (2) follows from Definition~\ref{definition:GnRtSys}~$({\rm{R}}1)$.
The claim (3) follows from Definition~\ref{definition:GnRtSys}~$({\rm{R}}4)$.
The claim (4) (resp. the claim (5), resp. the claim (6))
follows from Definition~\ref{definition:GnRtSys},
Lemma~\ref{lemma:rp-of-ax}
and \cite[Lemma~8~(iii)]{HY08}
(resp. 
\cite[Corollary~3]{HY08}, resp.
\cite[Corollary~2]{HY08}).
\hfill $\Box$

\subsection{Longest elements of a finite Weyl groupoid}\label{subsection:Lgst}
In Subsection~\ref{subsection:Lgst},
let $\mcC=\mcC(\fkI,\mcA,(\tvsigma_i)_{i\in\fkI},(\tC^a)_{a\in\mcA})$
and $\mcR=\mcR(\mcC,(\tR(a))_{a\in\mcA})$
be a Cartan scheme and
a generalized root system of type $\mcC$
respectively, and let $a\in\mcA$.

\begin{lemma} {\rm{(}}{\rm{\cite[{\it{Corollary}}~5]{HY08}.}}{\rm{)}} \label{lemma:lgslm}
Assume $|\tRp(a)|<\infty$.
Let $n:=|\tRp(a)|$.

{\rm{(1)}}
There exists a unique $1^aw_0\in \mcH(a,\_)$
such that $\atilell(1^aw_0)=n$. 

{\rm{(2)}} For $w\in \mcH(a,\_)$, it follows that 
\begin{equation}\label{eqn:lgsta}
w=1^aw_0\,\,\mbox{if and only if}\,\, w(\Pi)=-\Pi.
\end{equation}

{\rm{(3)}} 
For $a^\prime\in\mcA$ and $w\in\mcH(a,a^\prime)$, we have 
$n=\atilell(w)
+{\tilde {\ell}}_{a^\prime}(w^{-1}1^aw_0)$.
\end{lemma}
{\it {Proof.}}  
Let $a\in\mcA$. Let $w\in\mcH(a,a^\prime)$. Assume $\atilell(w)<n$.
By \eqref{eqn:dfell-rnl},
we have ${\tilde {\ell}}_{a^\prime}(w^{-1})=\atilell(w)$.
By Definition~\ref{definition:GnRtSys}~$({\rm{R}}3)$ and Lemme~\ref{lemma:non-initial}~(1),
there exists $i\in\fkI$ such that $w(\al_i)\in\tRp(a)$.
By \eqref{eqn:ellwc}, $\atilell(w \tils^{\tvsigma_ia^\prime}_i)=\atilell(w)+1$.
Thus we see the existence of $1^aw_0$.
Let $w^\prime$, $w^{\prime\prime}\in\mcH(a,\_)$
be such that $\atilell(w^\prime)=\atilell(w^{\prime\prime})=n$.
By Lemme~\ref{lemma:non-initial}~(1), we see that 
$w^\prime(\Pi)=w^{\prime\prime}(\Pi)=-\Pi$.
Since $(w^{\prime\prime})^{-1}w^\prime(\Pi)=\Pi$, by 
\eqref{eqn:dfell-rnl},
we have $(w^{\prime\prime})^{-1}w^\prime=\rmid_\bV$.
Hence $w^\prime=w^{\prime\prime}$. Thus we obtain the uniqueness of $1^aw_0$.
Thus we obtain the claims~(1) and (2).
Let $w_1\in\mcH(a,a^\prime)$. By \eqref{eqn:prellw}, $\atilell(w_1)\leq n$.
Assume $w_1\ne1^aw_0$. By the claim~(1), $\atilell(w_1)< n$.
By an argument as above, there exists $w_2\in\mcH(a^\prime,\_)$
such that $\atilell(w_1w_2)=n$
and ${\tilde{\ell}}_{a^\prime}(w_2)\leq n-\atilell(w_1)$.
By \eqref{eqn:dfell-rnl}, ${\tilde{\ell}}_{a^\prime}(w_2)= n-\atilell(w_1)$
since $\atilell(w_1w_2)=n$.
By the claim~(1), $w_1w_2=1^aw_0$.
Thus we obtain the claim~(3).
\hfill $\Box$
\newline\par
If $|\tRp(a)|<\infty$,
we call the only element
$1^aw_0$ (or more precisely, the pair $(a, 1^aw_0)$),
as in Lemma~\ref{lemma:lgslm}~(1), {\it{the longest element ending up with $a$}}.

\begin{lemma} \label{lemma:prlgelmmm}
{\rm{(}}See {\rm{\cite[{\it{Proposition}}~2.12]{CH09}.}}{\rm{)}}
Assume $|\tRp(a)|<\infty$. Let $n:=|\tRp(a)|$. Let $f\in\funcI_n$ be such that 
$1^a\tils_{f,n}=1^aw_0$.
Then
\begin{equation}\label{eqn:ellRp}
\tRp(a)=\{\,1^a\tils_{f,r-1}(\al_{f(r)})\,|\,r\in\fkJ_{1,n}\,\}.
\end{equation} In particular,
\begin{equation}\label{eqn:ellRp-a}
\tR(a)=\bigcup_{k=0}^\infty\bigcup_{f^\prime\in\funcI_k}1^a\tils_{f^\prime,k}(\Pi)\,.
\end{equation} 
\end{lemma}

{\it{Proof.}}
The equation \eqref{eqn:ellRp} 
follows from \eqref{eqn:defbhmell},
Lemma~\ref{lemma:lgslm}~(1), and \eqref{eqn:ellwf}.
The equation \eqref{eqn:ellRp-a} is clear from 
\eqref{eqn:ellRp}.
\hfill $\Box$

\begin{lemma} \label{lemma:lgelmmm}
Let $n\in\bN$ , $f\in\funcI_n$ and
$X:=\{\,1^a\tils_{f,r-1}(\al_{f(r)})\,
|\,r\in\fkJ_{1,n}\,\}(\subset\tR(a))$.
Assume
\begin{equation}\label{eqn:crtlg}
\Pi\subseteq X\subseteq\bZgeqoPi.
\end{equation}
Then $n=|\tRp(a)|$, $1^a\tils_{f,n}=1^aw_0$
and $\tRp(a)=X$.
\end{lemma}
{\it{Proof.}}
Let $w:=1^a\tils_{f,n}$.
It follows from \eqref{eqn:defbhmell},
Definition~\ref{definition:GnRtSys}~$({\rm{R}}3)$ 
and \eqref{eqn:ellwc} that
$\atilell(w)=n$.
By \eqref{eqn:ellwf}, $X=\preatilell(w)$.
Hence $\Pi\subseteq\preatilell(w)$.
By Lemma~\ref{lemma:non-initial}~(1), $X=\tRp(a)$.
Since $X=\preatilell(w)$,
we have $|X|=\atilell(w)=n$ by \eqref{eqn:defbhmell}.
By Lemma~\ref{lemma:lgslm}~(1), we have
$w=1^aw_0$. 
This completes the proof. 
\hfill $\Box$

\subsection{Technical fact}\label{subsection:techlem}
By Lemmas~\ref{lemma:preISOM} and \ref{lemma:lgelmmm}, we have
\begin{lemma}\label{lemma:techlm}
Keep the notation as in Definition~{\rm{\ref{definition:ISOM}}}.
Assume $|\tRp({\check{a}})|<\infty$. Then $(\mcR,{\check{a}})$
is quasi-isomorphic to $(\mcR^\prime,{\check{a}}^\prime)$ if and only if
\begin{equation*}
\tils^{{\check{a}}_{f,n}}_i=\tils^{{\check{a}^\prime}_{f,n}}_i\quad(n\in\bZgeqo,\,f\in\funcI_n,\,i\in\fkI).
\end{equation*}
In particular, $1^{{\check{a}}}w_0=1^{{\check{a}}^\prime}w_0$
and  $\tRp({\check{a}})=(\tR^\prime)^+({\check{a}}^\prime)$.
\end{lemma}

\section{Longest elements of finite Weyl groups}
\label{section:EcfWg-I}

\subsection{Root systems of type $\mathrm{A}-\mathrm{G}$}\label{subsection:RsAG}

In this section, we mention some longest elements of the finite 
Weyl group, or the crystallographic finite Coxeter group,
which will be used to study $\bhm$ treated in Theorem~\ref{theorem:cl}~(1)-(6)
below.

Let $\hatN\in\bN$.
Let $\bR^\hatN$ denote the $\hatN$-dimensional
$\bR$-linear space of $\hatN$-tuple column vectors, that is
$\bR^\hatN=\{\,{^t[}\underbrace{x_1,\ldots,x_\hatN}_\hatN]\,|\,x_i\in\bR\,(i\in\fkJ_{1,\hatN})\,\}$.
For $i\in\fkJ_{1,\hatN}$, 
let $\stbsis_i:={^t[}\underbrace{0,\ldots,0}_{i-1},1,\underbrace{0,\ldots,0}_{\hatN-i}]\in\bR^\hatN$,
that is, $\{e_i|i\in\fkJ_{1,\hatN}\}$ is the standard $\bR$-basis of $\bR^\hatN$.
For $m\in\fkJ_{1,\hatN}$, we regard $\bR^m$ as the $\bR$-linear subspace
$\oplus_{r=1}^m\bR\stbsis_r$ of $\bR^\hatN$.
For a subset $X$ of $\fkJ_{1,\hatN}$,  define the
$\bR$-linear map 
$\stProj_X:\bR^\hatN\to\bR^\hatN$
by $\stProj_X(\stbsis_i):=\stbsis_i$ ($i\in X$) and $\stProj_X(\stbsis_j):=0$ 
($j\in \fkJ_{1,\hatN}\setminus X$).
Let $\matMNR$ be the $\bR$-algebra of $\hatN\times\hatN$-matrices.
Let $\matGLNR$ be the group of invertible $\hatN\times\hatN$-matrices.
Let $
\hateta:\vecRN\times\vecRN\to\bR$ be the 
$\bR$-bilinear map
defined by $\hateta(\vece_k,\vece_r):=\delta_{kr}$.
For $v\in\vecRN\setminus\{0\}$, define $\hats_v\in\matGLNR$
by $\hats_v(u):=u-{\frac {2\hateta(u,v)} {\hateta(v,v)}}v$ $(u\in\vecRN)$,
that is $\hats_v$ is the {\it{reflection}} with respect to $v$. Note that
\begin{equation}\label{eqn:svsqur}
\hats_v^2=\rmid_{\vecRN}\quad(v\in\vecRN\setminus\{0\}),
\end{equation} and
\begin{equation}\label{eqn:svsqur-b}
\hateta(\hats_v(u),\hats_v(u^\prime))=\hateta(u,u^\prime)\quad(v\in\vecRN\setminus\{0\},
\,u,u^\prime\in\vecRN).
\end{equation} Using \eqref{eqn:svsqur} and \eqref{eqn:svsqur-b},
we have
\begin{equation}\label{eqn:conjsv}
\hats_v\hats_{v^\prime}\hats_v=\hats_{\hats_v(v^\prime)}
\quad(v,v^\prime\in\vecRN\setminus\{0\}).
\end{equation}

We say that a finite subset $\rsystem$ of $\vecRN\setminus\{0\}$ is
{\it{a crystallographic root system}} (in $\vecRN$) if $|\rsystem|<\infty$, 
$\hats_v(\rsystem)=\rsystem$ and
$\bR v\cap \rsystem=\{v,-v\}$
for all $v\in\rsystem$, 
and 
${\frac {2\hateta(v^\prime,v^{\prime\prime})}
{\hateta(v^\prime,v^\prime)}}\in\bZ$ 
for all $v^\prime$,
$v^{\prime\prime}\in\rsystem$, see \cite[1.2, 2.9]{Hum}.

Let $\rsystem$ be a crystallographic root system in $\vecRN$. 
We call $\rsystem$ {\it{irreducible}}
if for all $\hatbeta$, $\hatbeta^\prime\in\rsystem$,
there exist $r\in\bN$, and 
$\hatbeta_t\in\rsystem$ ($t\in\fkJ_{1,r}$) 
such that $\hateta(\hatbeta,\hatbeta_1)\ne 0$,
$\hateta(\hatbeta_t,\hatbeta_{t+1})\ne 0$ ($t\in\fkJ_{1,r-1}$)
and $\hateta(\hatbeta_r,\hatbeta^\prime)\ne 0$,
see \cite[2.2]{Hum} and \eqref{eqn:conjWeyl}.
We say that a subset $\rsysbase$ of $\rsystem$ is
{\it{a root basis}} of $\rsystem$ if $\rsysbase$ is a (set) $\bR$-basis of 
$\rmSpan_\bR(\rsysbase)$ and 
$\rsystem\subset\rmSpan_{\bZ_{\geq 0}}(\rsysbase)\cup-\rmSpan_{\bZ_{\geq 0}}(\rsysbase)$
(this is called {\it{a simple system}} in \cite[1.3, 2.9]{Hum}).

Let $\rsysbase$ be a root basis of $\rsystem$.
We call $\dim_\bR\rmSpan_\bR(\rsysbase)=|\rsysbase|$
the {\it{rank}} of $\rsystem$. 
By \cite[Corollary~1.5]{Hum}, we have
\begin{equation}\label{eqn:conjWeyl}
\rsystem=\hatWPi\cdot\rsysbase.
\end{equation} 
Let $\hatWPi$ be the subgroup of $\matGLNR$
generated by all $\hats_v$ with $v\in\rsysbase$.
We call $\hatWPi$ {\it{the Coxeter group associated with
$(\rsystem,\rsysbase)$}}.
Let $\hatSPi:=\{\,\hats_v\in\hatWPi\,|\,v\in\rsysbase\,\}$.
We call $(\hatWPi,\hatSPi)$
{\it{the Coxeter system associated with 
$(\rsystem,\rsysbase)$}}, see
\cite[1.9 and Theorem~1.5]{Hum}.
Let $\rsysbase$ be a root basis of $\rsystem$.
Let $\rsysp:=\rsystem\cap\rmSpan_{\bZ_{\geq 0}}(\rsysbase)$.
We call $\rsysp$ {\it{a positive root system of $\rsystem$
associated with $\rsysbase$}}
(this is called {\it{a positive system}} in \cite[1.3]{Hum}).

\begin{definition}\label{definition:strtbsis} {\rm{(See \cite[2.10]{Hum}.) 
Let $\rkN$ 
and $\fkI=\fkJ_{1,\rkN}$ be the one of 
\eqref{eqn:fixdfI}.
Let $\hatN\in\fkJ_{\rkN,\infty}$.
Let $\rsystem$ be a rank-$\rkN$ crystallographic root system in $\vecRN$.
Let $\rsysbase=\{\,\hatal_i\,|\,i\in\fkI\,\}$ be a root basis of $\rsystem$. 
Let $\orderrsysbase:=(\,\hatal_1,\hatal_2,\ldots,\hatal_\rkN)\in
\underbrace{\vecRN\times\cdots\times\vecRN}_\rkN$, so
$\orderrsysbase$ is an ordered $\bR$-basis of $\rmSpan_\bR(\rsystem)$.

{\rm{(1)}} Assume that $\rkN\geq 1$ and $\hatN=\rkN+1$. 
We call $\rsystem$ {\it{the ${\mathrm{A}}_\rkN$-type
standard root system}} if
\begin{equation}\label{eqn:str-A}
\rsystem=\{\,\vece_x-\vece_y\,|\,x,y\in\fkJ_{1,\rkN+1},\,x\ne y\,\}.
\end{equation}
We call $\orderrsysbase$ {\it{the ${\mathrm{A}}_\rkN$-data}} if 
$\hatal_i=\vece_i-\vece_{i+1}$ ($i\in\fkI$).

{\rm{(2)}} Assume $\rkN=\hatN\geq 2$. We call $\rsystem$ the 
{\it{${\mathrm{B}}_\rkN$-type
standard root system}} if
\begin{equation}\label{eqn:str-B}
\begin{array}{lcl}
\rsystem &= &\{\,c\vece_x+c^\prime\vece_y\,|\,x,y\in\fkJ_{1,\rkN},\,x< y,\,
c,\,c^\prime\in\{1,-1\}\,\} \\
& &\,\, \cup\{\,c^{\prime\prime}\vece_z\,|\,z\in\fkJ_{1,\rkN},\,
c^{\prime\prime}\in\{1,-1\}\,\}.
\end{array}
\end{equation}
We call $\orderrsysbase$ {\it{the ${\mathrm{B}}_\rkN$-data}} if 
$\hatal_i=\vece_i-\vece_{i+1}$ ($i\in\fkJ_{1,\rkN-1}$)
and $\hatal_\rkN=\vece_\rkN$.

{\rm{(3)}} Assume $\rkN=\hatN\geq 3$. We call $\rsystem$ 
{\it{the ${\mathrm{C}}_\rkN$-type
standard root system}} if
\begin{equation}\label{eqn:str-C}
\begin{array}{lcl}
\rsystem &=&\{\,c\vece_x+c^\prime\vece_y\,|\,x,y\in\fkJ_{1,\rkN},\,x< y,\,
c,\,c^\prime\in\{1,-1\}\,\} \\
& &\,\, \cup\{\,2 c^{\prime\prime}\vece_z\,|
\,z\in\fkJ_{1,\rkN},\,
c^{\prime\prime}\in\{1,-1\}\,\}.
\end{array}
\end{equation}
We call $\orderrsysbase$ {\it{the ${\mathrm{C}}_\rkN$-data}} if 
$\hatal_i=\vece_i-\vece_{i+1}$ ($i\in\fkJ_{1,\rkN-1}$)
and $\hatal_\rkN=2\vece_\rkN$.

{\rm{(4)}} Assume $\rkN=\hatN\geq 4$. We call $\rsystem$ {\it{the ${\mathrm{D}}_\rkN$-type
standard root system}} if
\begin{equation}\label{eqn:str-D}
\rsystem=\{\,c\vece_x+c^\prime\vece_y\,|\,x,y\in\fkJ_{1,\rkN},\,x< y,\,
c,\,c^\prime\in\{1,-1\}\,\}.
\end{equation}
We call $\orderrsysbase$ {\it{the ${\mathrm{D}}_\rkN$-data}} if 
$\hatal_i=\vece_i-\vece_{i+1}$ ($i\in\fkJ_{1,\rkN-1}$)
and $\hatal_\rkN=\vece_{\rkN-1}+\vece_\rkN$.

{\rm{(5)}} Assume that $\rkN=6$ and $\hatN=8$. We call $\rsystem$ {\it{the ${\mathrm{E}}_6$-type
standard root system}} if
\begin{equation*}
\begin{array}{lcl}
\rsystem &= &\{\,c\vece_x+c^\prime\vece_y\,|\,x,y\in\fkJ_{1,5},\,x< y,\,
c,\,c^\prime\in\{1,-1\}\,\} \\
& & \,\, 
\cup\{\,{\frac 1 2}((\sum_{r=1}^5 c_r\vece_r)+(\prod_{k=1}^5c_k)(\vece_6-\vece_7-\vece_8))\,|\,
c_r\in\{1,-1\}\,(r\in\fkJ_{1,5})\,\}.
\end{array}
\end{equation*}
We call $\orderrsysbase$ {\it{the ${\mathrm{E}}_6$-data}} if 
$\hatal_1={\frac 1 2}(\vece_1+\vece_8-\sum_{r=2}^7 \vece_r))$,
$\hatal_2=\vece_1+\vece_2$ and
$\hatal_i=
\vece_{i-1}-\vece_{i-2}$ ($i\in\fkJ_{3,6}$).

{\rm{(6)}} Assume that $\rkN=7$ and $\hatN=8$. We call $\rsystem$ 
{\it{the ${\mathrm{E}}_7$-data}} if 
\begin{equation*}
\begin{array}{lcl}
\rsystem &= &\{\,c\vece_x+c^\prime\vece_y\,|\,x,y\in\fkJ_{1,6},\,x< y,\,
c,\,c^\prime\in\{1,-1\}\,\} \\
& & \,\, \cup\{\,c^{\prime\prime}(\vece_7-\vece_8)
\,|\,c^{\prime\prime}\in\{1,-1\}\,\} \\
& & \,\, 
\cup\{\,{\frac 1 2}((\sum_{r=1}^6 c_r\vece_r)-(\prod_{k=1}^6c_k)(\vece_7-\vece_8))\,|\,
c_r\in\{1,-1\}\,(r\in\fkJ_{1,6})\,\}.
\end{array}
\end{equation*}

We call $\orderrsysbase$ {\it{the ${\mathrm{E}}_7$-data}} if 
$\hatal_1={\frac 1 2}(\vece_1+\vece_8-\sum_{r=2}^7 \vece_r)$,
$\hatal_2=\vece_1+\vece_2$ and
$\hatal_i=
\vece_{i-1}-\vece_{i-2}$ ($i\in\fkJ_{3,7}$).

{\rm{(7)}} Assume $\rkN=\hatN=8$. We call $\rsystem$ {\it{the ${\mathrm{E}}_8$-type
standard root system}} if
\begin{equation*}
\begin{array}{lcl}
\rsystem &= &\{\,c\vece_x+c^\prime\vece_y\,|\,x,y\in\fkJ_{1,8},\,x< y,\,
c,\,c^\prime\in\{1,-1\}\,\} 
 \\
& & \,\, 
\cup\{\,{\frac 1 2}((\sum_{r=1}^7 c_r\vece_r)+(\prod_{k=1}^7c_k)\vece_8)\,|\,
c_r\in\{1,-1\}\,(r\in\fkJ_{1,7})\,\}.
\end{array}
\end{equation*}

We call $\orderrsysbase$ {\it{the ${\mathrm{E}}_8$-data}} if 
$\hatal_1={\frac 1 2}(\vece_1+\vece_8-\sum_{r=2}^7 \vece_r))$,
$\hatal_2=\vece_1+\vece_2$ and
$\hatal_i=
\vece_{i-1}-\vece_{i-2}$ ($i\in\fkJ_{3,8}$).

{\rm{(8)}} Assume $\rkN=\hatN=4$. We call $\rsystem$ {\it{the ${\mathrm{F}}_4$-type
standard root system}} if
\begin{equation*}
\begin{array}{lcl}
\rsystem &= &\{\,c\vece_x+c^\prime\vece_y\,|\,x,y\in\fkJ_{1,4},\,x< y,\,
c,\,c^\prime\in\{1,-1\}\,\} \\
& & \,\, 
\cup\{\,c^{\prime\prime}\vece_z\,|
\,z\in\fkJ_{1,4},\,
c^{\prime\prime}\in\{1,-1\}\,\} \\
& & \,\, 
\cup\{\,{\frac 1 2}\sum_{r=1}^4 c_r\vece_r\,|\,
c_r\in\{1,-1\}\,(r\in\fkJ_{1,4})\,\}.
\end{array}
\end{equation*}
We call $\orderrsysbase$ {\it{the ${\mathrm{F}}_4$-data}} if
$\hatal_1=\vece_2-\vece_3$, $\hatal_2=\vece_3-\vece_4$, 
$\hatal_3=\vece_4$ and
$\hatal_4={\frac 1 2}(\vece_1-\vece_2-\vece_3-\vece_4)$.

{\rm{(9)}} Assume  that $\rkN=2$ and $\hatN=3$. 
We call $\rsystem$ {\it{the ${\mathrm{G}}_2$-type
standard root system}} if
\begin{equation*}
\begin{array}{lcl}
\rsystem &= &\{\,c(\vece_x-\vece_y)\,|\,x,y\in\fkJ_{1,3},\,x< y,\,
c\in\{1,-1\}\,\} \\
& & \,\, 
\cup\{\,c^\prime(2\vece_{z_1}-\vece_{z_2}-\vece_{z_3})\,|
\,\,\{\,z_1,\,z_2,\,z_3\,\}=\fkJ_{1,3},\,
c^\prime\in\{1,-1\}\,\}. 
\end{array}
\end{equation*}
We call $\orderrsysbase$ {\it{the ${\mathrm{G}}_2$-data}} if 
$\hatal_1=\vece_1-\vece_2$ and 
$\hatal_2=-2\vece_1+\vece_2+\vece_3$. 

{\rm{(10)}} Let $\rsystem$ and $\rsysbase$ be the ones treated 
in the above (1)-(9).
We call $\rsystem$ 
{\it{a rank-$\rkN$ standard irreducible root system}}.
We call $\orderrsysbase$ 
{\it{a rank-$\rkN$ Cartan data.}}

}}
\end{definition}

It is well-known that rank-$\rkN$ irreducible crystallographic root systems
are isomorphic to exactly the rank-$\rkN$ standard irreducible root systems
(cf. \cite[2.10]{Hum}).

\begin{definition}\label{definition:mcRordbs}
Let $\orderrsysbase$ be a rank-$\rkN$ Cartan data.
Define the $\bR$-linear isomorphism 
\begin{equation*}
\hatiota:\rmSpan_\bR(\rsysbase)\to\bV
\end{equation*}
by $\hatiota(\hatal_i)=\al_i$
$(i\in\fkI)$.
Let $\mcR_\orderrsysbase=\mcR(\mcC_\orderrsysbase,(\tR(a))_{a\in\mcA_\orderrsysbase})$
be the generalized root systems such that 
$|\mcA_\orderrsysbase|=1$ and
\begin{equation} \label{eqn:crrCox-1}
\tils^a_i=\hatiota\circ\hats_{\hatal_i}\circ\hatiota^{-1}
\quad(a\in\mcA_\orderrsysbase,\,i\in\fkI).
\end{equation} By \eqref{eqn:ellRp-a} and
\eqref{eqn:conjWeyl}, we have
\begin{equation} \label{eqn:crrCox-2}
\hatiota(\rsystem)=\tR(a)\quad(a\in\mcA_\orderrsysbase),
\end{equation} where $\rsystem$ is the rank-$\rkN$ root system
corresponding to $\orderrsysbase$.
\end{definition}

\begin{proposition} \label{proposition:crrCox}
The correspondence $\orderrsysbase\mapsto\mcR_\orderrsysbase$
from
the set of all rank-$\rkN$ Cartan datas
to
the family of all connected generalized root systems
$\mcR=\mcR(\mcC,(\tR(a))_{a\in\mcA})$
with $|\mcA|=1$ {\rm{(}}and $|\fkI|=\rkN${\rm{)}}
is injective.
Moreover it is surjective up to isomorphisms
in the sense of Definition~{\rm{\ref{definition:ISOM}~(2)}}.
\end{proposition}
{\it{Proof.}} See \cite[Proposition~4.9]{Kac90}.
\hfill $\Box$

\subsection{Longest elements of irreducible Weyl groups}
\label{subsection:LeoiWG}
In this subsection,
let $\orderrsysbase=(\hatal_i|i\in\fkI)$ 
be a rank-$\rkN$ Cartan data, and
we also treat $\mcR_\orderrsysbase=\mcR(\mcC_\orderrsysbase,(\tR(a))_{a\in\mcA_\orderrsysbase})$,
see \eqref{eqn:crrCox-1}. Let $a\in\mcA_\orderrsysbase$.
Define the map $\hatlng:\hatWPi\to\bZgeqo$
in the following way, see \cite[1.6]{Hum}.
Let $\hatlng(1):=0$, where $1$ is a unit of $\hatWPi$.
Note that an arbitrary $\hatw\in\hatWPi$
can be written as a product
of finite $\hats_\hatbeta$'s with some $\hatbeta\in\rsysbase$,
say $\hatw=\underbrace{\hats_{\hatbeta_1}\cdots\hats_{\hatbeta_r}}_r$
for some $r\in\bN$ and some $\hatbeta_x\in\rsysbase$
($x\in\fkJ_{1,r}$). If $\hatw\ne 1$, let
$\hatlng(\hatw)$ be the smallest $r$ for which such an expression exists,
and call the expression {\it{reduced}}.
By \eqref{eqn:dfell-rnl} and \eqref{eqn:crrCox-1}, we have
$\hatlng(\hatw)=\atilell(\hatiota\circ\hatw\circ\hatiota^{-1})$. 
We call $\hatlng(\hatw)$ {\it{the length of $\hatw$}}.
Let
\begin{equation*}
\whLng(\hatw):=\{\,\hatbeta\in\rsysp\,|\,\hatw(\hatbeta)\in-\rsysp\}\quad(\hatw\in\hatWPi),
\end{equation*}
so $\whLng(\hatw)=\preatilell(\hatiota\circ\hatw\circ\hatiota^{-1})$
by
\eqref{eqn:prellw} and \eqref{eqn:crrCox-1}.
By \eqref{eqn:defbhmell} and \eqref{eqn:crrCox-1},
\begin{equation}\label{eqn:proflgth-a}
\hatlng(\hatw)=|\whLng(\hatw)|,
\end{equation} see also \cite[Corollary~1.7]{Hum}.
By \eqref{eqn:hecpronedash} and \eqref{eqn:crrCox-1},
\begin{equation}\label{eqn:proflgth-b}
\hats_\hatal(\rsysp\setminus\{\hatal\})=\rsysp\setminus\{\hatal\}
\quad(\hatal\in\rsysbase),
\end{equation} see also \cite[Propsoition~1.4]{Hum}. 
By \eqref{eqn:ellwc} and \eqref{eqn:crrCox-1},
\begin{equation}\label{eqn:proflgth-c}
\hatlng(\hatw\hats_\hatal)=
\left\{\begin{array}{ll}
\hatlng(\hatw)+1 & \mbox{if $\hatw(\hatal)\in\rsysp$}, \\
\hatlng(\hatw)-1 & \mbox{if $\hatw(\hatal)\in-\rsysp$},
\end{array}\right.
\end{equation} for $\hatal\in\rsysbase$, 
see also \cite[Lemma~1.6 and Corollary~1.7]{Hum}.

Assume that $|\rsystem|<\infty$.
By the above properties, we can see that there exists a 
unique $\hatwo\in\hatWPi$ such that 
$\hatwo(\rsysbase)=-\rsysbase$, see \cite[1.8]{Hum}.
It is well-known that $\hatlng(\hatwo)=|\rsysp|$, that
$\hatwo$ is the only element of $\hatWPi$ that $\hatlng(\hatw)\leq\hatlng(\hatwo)$
for all $\hatw\in\hatWPi$, and that
\begin{equation}\label{eqn:progrouplg} \hatlng(\hatw)=\hatlng(\hatwo)-\hatlng(\hatwo\hatw^{-1})
\quad\mbox{for all $\hatw\in\hatWPi$}.
\end{equation} 
We call $\hatwo$ 
{\it{the longest element of the Coxeter system}} of
$(\hatWPi,\hatSPi)$.
Note that 
\begin{equation}\label{eqn:proflgth-d}
\hatwo=\hatiota^{-1}\circ 1^a w_0\circ\hatiota\quad\mbox{and}\quad
\hatlng(\hatwo)=\atilell(1^a w_0).
\end{equation}
It is well-known that
\begin{equation}\label{eqn:ctprpoflg}
\hatlng(\hatwo)=|\rsysp|,
\end{equation} see also Lemma~\ref{lemma:lgslm}.
Let $n:=\hatlng(\hatwo)$, and let 
$\underbrace{\hats_{\hatbeta_1}\hats_{\hatbeta_2}\cdots\hats_{\hatbeta_n}}_n$
be the reduced expression of $\hatwo$,
where $\hatbeta_k$'s are some elements of $\rsysbase$.
Then we know well that 
\begin{equation}\label{eqn:ctprpoflg-b}
\rsysp:=\{\,\underbrace{\hats_{\hatbeta_1}\hats_{\hatbeta_2}\cdots\hats_{\hatbeta_{k-1}}}_{k-1}(\hatbeta_k)\,|\,k\in\fkJ_{1,n}\,\},
\end{equation} see also \eqref{eqn:ellRp}.

\begin{proposition}\label{proposition:concretewo}
Let $\orderrsysbase=(\,\hatal_1,\hatal_2,\ldots,\hatal_\rkN)$ be 
a rank-$\rkN$ Cartan data. Let $\rsysbase=\{\,\hatal_i\,|\,i\in\fkI\,\}$.
Let $\hats_i:=\hats_{\hatal_i}$ $(i\in\fkI)$.
Let $\hatwo$ be the longest element of the Coxeter system
$(\hatWPi,\hatSPi)$. 
Let $h:={\frac {2|\rsysp|} {\rkN}}$.
Let $V:=\rmSpan_\bR(\rsysbase)$.

{\rm{(1)}} $h\in\bN$.

{\rm{(2)}} 
Assume that $\orderrsysbase$ is neither the ${\mathrm{A}}_\rkN$-data,
the ${\mathrm{D}}_\rkN$-data, nor the ${\mathrm{E}}_6$-data.
Then there is no bijection $u:\fkI\to\fkI$ 
such that $u\ne \rmid_\fkI$ and $\hateta(\hatal_i,\hatal_j)=
\hateta(\hatal_{u(i)},\hatal_{u(j)})$
$(i,\,j\in\fkI)$.
Moreover $(\hatwo)_{|V}=-\rmid_V$.
Furthermore $h\in 2\bN$, and
$(\underbrace{\hats_1\hats_2\cdots\hats_\rkN}_\rkN)^{\frac h 2}$
is a reduced expression of $\hatwo$.

{\rm{(3)}}
Assume that $\orderrsysbase$ is the ${\mathrm{A}}_\rkN$-data. Then
$\hatwo(\vece_x)=\vece_{\hatN-x+1}$ 
$(x\in\fkJ_{1,\hatN})$, so $\hatwo(\hatal_i)=-\hatal_{\rkN-i+1}$ 
$(i\in\fkI)$. Moreover
\begin{equation}\label{eqn:concretewo-AN}
\hatwo=(\underbrace{\hats_1\hats_2\cdots\hats_\rkN}_\rkN)(\underbrace{\hats_1\hats_2\cdots\hats_{\rkN-1}}_{\rkN-1})
\cdots(\underbrace{\hats_1\hats_2}_2)\underbrace{\hats_1}_1, 
\end{equation} and the RHS of \eqref{eqn:concretewo-AN} is a reduced expression of $\hatwo$.

{\rm{(4)}}
Assume that $\orderrsysbase$ is the ${\mathrm{D}}_\rkN$-data. 
If $\rkN\in 2\bN$, $\hatwo=-\rmid_{\vecRN}$. 
If $\rkN\in 2\bN-1$, then $\hatwo(\hatal_i)=-\hatal_i$
$(i\in\fkJ_{1,\rkN-2})$,
$\hatwo(\hatal_{\rkN-1})=-\hatal_\rkN$ and $\hatwo(\hatal_\rkN)=-\hatal_{\rkN-1}$.
Moreover $(\underbrace{\hats_1\hats_2\cdots\hats_\rkN}_\rkN)^{\rkN-1}$ is 
a reduced expression of $\hatwo$. Furthermore, for 
$r\in\fkJ_{1,\rkN-1}$, we have
\begin{equation} \label{eqn:hatsrsN}
(\underbrace{\hats_r\hats_{r+1}\cdots\hats_\rkN}_{\rkN-r+1})^{\rkN-r}
=\stProj_{\fkJ_{1,r-1}}-\stProj_{\fkJ_{r,\rkN-1}}+(-1)^{\rkN-r}\stProj_{\fkJ_{\rkN,\rkN}}.
\end{equation}

{\rm{(5)}} 
If $\rkN=6$, and $\orderrsysbase$ is 
the ${\mathrm{E}}_6$-data.
then $h=12$ and $(\hats_1\hats_3\hats_5\hats_2\hats_4\hats_6)^6$ is 
a reduced expression of $\hatwo$.

\end{proposition}
{\it{Proof.}}
Let $\hatb:=\underbrace{\hats_1\hats_2\cdots\hats_\rkN}_\rkN$.
Let $h^\prime$ be the order of $\hatb$.
Then $\hatb$ and $h^\prime$ are called {\it{a Coxeter element}}
and {\it{the Coxeter number}} respectively, see \cite[Exercise~3.19]{Hum}. 
By \cite[Proposition~3.18]{Hum}, we have
\begin{equation} \label{eqn:fctexp}
h^\prime=h.
\end{equation} 
Fix $\zeta\in \bC^\times_{h^\prime}$.
It is clear from \eqref{eqn:svsqur-b} that 
$\hatb$ acts on the $\rkN$-dimensional $\bC$-linear space $V\otimes_\bR\bC$
as a diagonalizable linear map whose eigenvalues are
$\zeta^m$ for some $m\in\fkJ_{0,h^\prime-1}$;
these integers $m$'s are called {\it{the exponents}},
see \cite[3.16]{Hum}.

(1) This is clear from \eqref{eqn:fctexp}.  

(2) The first claim is clear.
By \eqref{eqn:proflgth-a} and \eqref{eqn:ctprpoflg}, $\hatwo(\rsysp)=-\rsysp$,
so $\hatwo(\rsysbase)=-\rsysbase$.
Then the second claim follows from the first claim and \eqref{eqn:svsqur-b}.
The third claim follows from \eqref{eqn:proflgth-a}, 
\eqref{eqn:ctprpoflg}, \eqref{eqn:fctexp} and the fact that
$h$ is even and
all the exponents are odd, see 
\cite[Tables~3.1 and 3.2,\,Theorem~3.19]{Hum}.

(3) Let $\hatb_i:=\underbrace{\hats_1\hats_2\cdots\hats_i}_i$
($i\in\fkI$). Let $\hatwo^\prime:=\underbrace{\hatb_\rkN\hatb_{\rkN-1}\cdots\hatb_1}_\rkN$.
Then $\hatb_i(\vece_x)=\vece_{x+1}$ ($x\in\fkJ_{1,i}$),
$\hatb_i(\vece_{i+1})=\vece_1$, and 
$\hatb_i(\vece_y)=\vece_y$ ($y\in\fkJ_{i+1,\hatN}$).
Hence 
$\hatwo^\prime(\vece_x)=\vece_{\hatN-x+1}$ 
($x\in\fkJ_{1,\hatN}$). In particular, $\hatwo^\prime(\rsysbase)=-\rsysbase$.
By \eqref{eqn:str-A}, $|\rsysp|={\frac {\rkN(\rkN-1)} 2}$.
Hence the claim~(3) follows from \eqref{eqn:ctprpoflg}.

(4) Let $r\in\fkJ_{1,\rkN-1}$, and 
$\hatb:=\underbrace{\hats_r\hats_{r+1}\cdots\hats_\rkN}_{\rkN-r+1}$.
Then $\hatb(\vece_x)=\vece_x$ ($x\in\fkJ_{1,r-2}$),
$\hatb(\vece_y)=\vece_{y+1}$ ($y\in\fkJ_{1,\rkN-2}$),
$\hatb(\vece_{\hatN-1})=-\vece_r$
and $\hatb(\vece_\hatN)=-\vece_\hatN$.
Then we obtain the claim~(4) in way similar to that for the claim~(3).

(5) This can be proved directly. 
\hfill $\Box$

\subsection{On longest elements of type-A classical Weyl group}

\begin{proposition} \label{proposition:ANNlgest}
Let $\orderrsysbase=(\hatal_1,\ldots,\hatal_\rkN)$ be the
${\mathrm{A}}_\rkN$-data, and $\rsysbase:=\{\,\al_i\,|\,i\in\fkI\,\}$.
Let $\hatwo$ be the longest element of the Coxeter system
$(\hatWPi,\hatSPi)$. Let $\hats_i:=\hats_{\hatal_i}\in\hatSPi$ 
$(i\in\fkI)$.
Let $n:={\frac {\rkN(\rkN+1)} 2}$.
Let $m\in\fkJ_{1,\rkN}$ and $r:=n-m(\rkN-m+1)$.
Then 
\begin{equation} \label{eqn:ANNlgest-0}
\hatlng(\hatwo)=n.
\end{equation}
Moreover there exists $f\in\funcI_n$
such that $\{f(t)|t\in\fkJ_{1,r}\}=\fkI\setminus\{m\}$, $f(r+1)=m$
and $\underbrace{\hats_{f(1)}\hats_{f(2)}\cdots\hats_{f(l)}}_n$
is 
a
reduce expression 
of $\hatwo$.
Furthermore for such $f$, we have
\begin{equation} \label{eqn:ANNlgest-a}
\begin{array}{l}
\{\,\hats_{f(1)}\hats_{f(2)}\cdots\hats_{f(k-1)}(\hatal_{f(k)})\,|\,
k\in\fkJ_{1,r}\,\} \\
\quad =\{\,\vece_x-\vece_{x^\prime}\,|\,
x,\,x^\prime\in\fkJ_{1,m},\,x<x^\prime\,\} \\
\quad \quad \,\,
\cup\{\,\vece_y-\vece_{y^\prime}\,|\,
y,\,y^\prime\in\fkJ_{m+1,\rkN+1},\,y<y^\prime\,\},
\end{array}
\end{equation} and
\begin{equation} \label{eqn:ANNlgest-b}
\begin{array}{l}
\{\,\hats_{f(1)}\hats_{f(2)}\cdots\hats_{f(t-1)}(\hatal_{f(t)})\,|\,
t\in\fkJ_{r+1,\rkN}\,\} \\
\quad =\{\,\vece_z-\vece_{z^\prime}\,|\,
z\in\fkJ_{1,m},\,z^\prime\in\fkJ_{m+1,\rkN+1}\,\}.
\end{array}
\end{equation}
\end{proposition}
{\it{Proof.}}
The claim \eqref{eqn:ANNlgest-0} follows from \eqref{eqn:ctprpoflg} and \eqref{eqn:str-A}.
Let $\rsystem^\prime$ be the crystallographic root system in $\vecRN$
defined as
the RHS of \eqref{eqn:ANNlgest-a}.
Let $\rsysbase^\prime:=\rsysbase\setminus\{\hatal_m\}$.
Then $\rsysbase^\prime$ is a root basis of $\rsystem^\prime$,
Note that the Coxeter system
$({\hat{W}}(\rsysbase^\prime),{\hat{S}}(\rsysbase^\prime))$ 
is isomorphic to the product of
the Coxeter systems of types ${\mathrm{A}}_{m-1}$ and 
${\mathrm{A}}_{\rkN-m}$ (resp. the Coxeter system of type ${\mathrm{A}}_{\rkN-1}$)
if $m\in\fkJ_{2,\rkN-1}$ (resp. $m\in\{1,\rkN\}$).
Note that $r$ equals the length of the longest element of 
$({\hat{W}}(\rsysbase^\prime),{\hat{S}}(\rsysbase^\prime))$.
Then the remaining claims
follow from these facts,
\eqref{eqn:str-A}, \eqref{eqn:progrouplg} and \eqref{eqn:ctprpoflg-b}.
\hfill $\Box$
\newline\par
\begin{remark}
A reduced expression as in Proposition~\ref{proposition:ANNlgest}
is given by 
\begin{equation}\label{eqn:srAm}
\begin{array}{lcl}
\hatwo&=&(\underbrace{\hats_1\hats_2\cdots\hats_m}_m)(\underbrace{\hats_1\hats_2\cdots\hats_{m-1}}_{m-1})
\cdots(\underbrace{\hats_1\hats_2}_2)\underbrace{\hats_1}_1  \\
& &\cdot(\underbrace{\hats_{m+2}\hats_{m+3}\cdots\hats_\rkN}_{\rkN-m-1})
(\underbrace{\hats_{m+2}\hats_{m+3}\cdots\hats_{\rkN-1}}_{\rkN-m-2})
\cdots(\underbrace{\hats_{m+2}\hats_{m+3}}_2)\underbrace{\hats_{m+2}}_1 \\
& &\cdot(\underbrace{\hats_{m+1}\hats_{m+2}\cdots\hats_\rkN}_{\rkN-m})
(\underbrace{\hats_m\hats_{m+1}\cdots\hats_{\rkN-1}}_{\rkN-m})
\cdots(\underbrace{\hats_1\hats_2\cdots\hats_{\rkN-m}}_{\rkN-m}).  \\
\end{array}
\end{equation} This can be proved in a way similar
to that for Proposition~\ref{proposition:concretewo}~(1).
\end{remark}

\subsection{A longest element of type-B classical Weyl group}
\label{subsection:LEofBWG}

\begin{lemma}\label{lemma:BNNlgest}
Let $\orderrsysbase=(\hatal_1,\ldots,\hatal_\rkN)$ be the
${\mathrm{B}}_\rkN$-data, and $\rsysbase:=\{\,\hatal_i\,|\,i\in\fkI\,\}$.
Let $\hatwo$ be the longest element of
$(\hatWPi,\hatSPi)$. Let $\hats_i:=\hats_{\hatal_i}\in\hatSPi$ 
$(i\in\fkI)$.

{\rm{(1)}} Let $k$, $r\in\fkJ_{1,\rkN}$
with $k\leq r$. 
Let $\hatb:=-\stProj_{\fkJ_{k,r}}
+\stProj_{\fkJ_{1,\hatN}\setminus\fkJ_{k,r}}\in\rmGL_N(\bR)$.
Then $\hatb\in\hatWPi$ and 
\begin{equation}\label{eqn:EXShatb}
(\underbrace{\hats_k\hats_{k+1}\cdots\hats_{\rkN-1}\hats_\rkN\hats_{\rkN-1}\cdots\hats_{r+1}
\hats_r}_{2\rkN-k-r+1})^{r-k+1}
=\hatb.
\end{equation} 
Moreover the LHS of \eqref{eqn:EXShatb}
is a reduced expression of $\hatb$.

{\rm{(2)}}
Let $k$, $t$, $r\in\fkJ_{1,\hatN}$ be such that $k\leq t<r$.
Let $\hatb\in\hatWPi$ be as in the claim {\rm{(1)}}.
Define $\hatb_1\in\hatWPi$ 
{\rm{(}}resp. $\hatb_2\in\hatWPi${\rm{)}} be the one define in the same way as that for
$\hatb$ with $k$ and $t$ {\rm{(}}resp. $t+1$ and $r${\rm{)}} in place of 
$k$ and $r$ respectively. 
Then $\hatb=\hatb_1\hatb_2=\hatb_2\hatb_1$ 
and $\hatlng(\hatb)=\hatlng(\hatb_1)+\hatlng(\hatb_2)$.

{\rm{(3)}}
Let $m\in\fkJ_{1,\rkN-1}$.
Then
\begin{equation}\label{eqn:srBm}
\begin{array}{lcl}
\hatwo&=&(\underbrace{\hats_{\rkN-m+1}\hats_{\rkN-m+2}\cdots\hats_\rkN}_m)^m \\
& & \cdot (\underbrace{\hats_1\hats_2\cdots\hats_{\rkN-1}\hats_\rkN\hats_{\rkN-1}\cdots\hats_{\rkN-m+1}
\hats_{\rkN-m}}_{\rkN+m})^{\rkN-m}.
\end{array}
\end{equation}
Moreover the RHS of \eqref{eqn:srBm}
is a reduced expression of $\hatwo$. In particular,
\begin{equation}\label{eqn:srBm-2}
\hatlng(\hatwo)=\rkN^2.
\end{equation}
\end{lemma}
\proof
(1) 
Let $\hatb^\prime\in\hatWPi$ be the LHS of \eqref{eqn:EXShatb}.
By \eqref{eqn:conjsv}, we have
\begin{equation*}
\underbrace{\hats_r\hats_{r+1}\cdots\hats_{\rkN-1}\hats_\rkN\hats_{\rkN-1}\cdots\hats_{r+1}
\hats_r}_{2\rkN-2r+1}=\hats_{\brre_r}.
\end{equation*} Hence
\begin{equation*}
\hatb^\prime=(\underbrace{\hats_k\hats_{k+1}\cdots\hats_{r-1}}_{r-k}\hats_{\brre_r})^{r-k+1}.
\end{equation*} Then
by the same claim as Proposition~\ref{proposition:concretewo}~(2)
for the $\mathrm{B}_\rkN$-data
with $r-k$ in place of $\hatN$, we have $\hatb=\hatb^\prime\in\hatWPi$.
We see that
\begin{equation*}
\whLng(\hatb)=\{\,\vece_t\,|\,t\in\fkJ_{k,r}\}\cup
\{\,\vece_t+c\vece_{t^\prime}\,|\,c\in\{-1,1\},\,t\in\fkJ_{k,r},\,t^\prime\in\fkJ_{t^\prime,\rkN}\,\}.
\end{equation*} Hence by \eqref{eqn:proflgth-a}, we have
\begin{equation}\label{eqn:cphtbkr}
\begin{array}{l}
\hatlng(\hatb)=(r-k+1)+2\sum_{t=k}^r(\rkN-t) \\
\quad=(r-k+1)+2\rkN(r-k+1)-2({\frac {r(r+1)} 2}-{\frac {k(k-1)} 2}) \\
\quad=(r-k+1)(1+2\rkN-(r+k)) \\
\quad=(2\rkN-k-r+1)(r-k+1). 
\end{array}
\end{equation} Hence we obtain the last statement.

(2)
The first statement is clear.
The second statement follows from the claim (1)
and \eqref{eqn:cphatbktr-b} below.
\begin{equation}\label{eqn:cphatbktr-b}
\begin{array}{l}
\hatlng(\hatb_1)+\hatlng(\hatb_2) \\
\quad=(2\rkN-k-t+1)(t-k+1)+(2\rkN-t-r)(r-t) \\
\quad=2\rkN(r-k+1)-(k+t-1)(t-k+1)-(t+r)(r-t) \\
\quad=2\rkN(r-k+1)-(-k^2+t^2+2k-1)-(r^2-t^2) \\
\quad=2\rkN(r-k+1)+(k^2-r^2-2k+1) \\
\quad=2\rkN(r-k+1)+(k-1+r)(k-1-r) \\
\quad=(2\rkN-r-k-1)(r-k+1) \\
\quad=\hatlng(\hatb).
\end{array}
\end{equation}

(3) 
This follows immediately from the claims (1), (2) and 
\eqref{eqn:ctprpoflg}. 
\qed

\section{Longest elements of Weyl groupoids of the simple Lie superalgebra of type ABCD}
\label{section:LGEsuperCD}
\subsection{Super-data}\label{subsection:Superdatas}

Let $\barN\in\bN$.
Let $\{\,\brre_i\,|\,i\in\fkJ_{1,\barN}\,\}$ be the standard 
$\bR$-basis of $\bR^\barN$.
Let $m\in\fkJ_{0,\barN}$.
Let $\mcA_{m|\barN-m}$ be the set of all maps $p:\fkJ_{1,\barN}\to\fkJ_{0,1}$
with $\sum_{i=1}^\barN p(i)=m$.
For $p\in\mcA_{m|\barN-m}$, define the $\bR$-bilinear map $\breta^p:\bR^\barN\times\bR^\barN\to\bR$ by
$\breta^p(\vece_i,\vece_j)=\delta(p(i),p(j))(-1)^{p(i)}$.
Define $p^+_{m|\barN-m}\in\mcA_{m|\barN-m}$ by
$p^+_{m|\barN-m}(i)=0$ ($i\in\fkJ_{1,m}$) and $p^+_{m|\barN-m}(j)=1$ ($j\in\fkJ_{m+1,\barN}$).
Define $p^-_{m|\barN-m}\in\mcA_{m|\barN-m}$ by
$p^-_{m|\barN-m}(i)=1$ ($i\in\fkJ_{1,m}$) and $p^-_{m|\barN-m}(j)=0$ ($j\in\fkJ_{m+1,\barN}$).

The sets $\barR$ given in Definition~\ref{definition:defSuperbhm} below are
really the sets of the roots of finite dimensional contragredient Lie superalgebras,
see Theorem~\ref{theorem:multisuper}. 

\begin{definition} \label{definition:defSuperbhm}
{\rm{
Keep the notation as above. We also use the terminology in 
Definition~\ref{definition:strtbsis}.
Let 
\begin{equation*}
\orderbarPi=(\bal_i|i\in\fkI)=(\,\bal_1,\bal_2,\ldots,\bal_\rkN)\in
\underbrace{\bR^\barN\times\cdots\times\bR^\barN}_\rkN.
\end{equation*} 
Let $\barR$ be a subset of $\bR^\barN$.
Let $\parity:\fkI\to\fkJ_{0,1}$ be a map.

(1) Assume that $\barN-1=\rkN\geq 2$ and $m\in\fkJ_{1,\rkN}$.
We call $(\breta^{p^+_{m|\rkN+1-m}},\orderbarPi)$ {\it{the ${\mathrm{A}}(m-1,\rkN-m)$-data}}
if $\orderbarPi$ is the ${\mathrm{A}}_\rkN$-data.
We call $\barR$ {\it{the ${\mathrm{A}}(m-1,\rkN-m)$-type standard root system}} if
$\barR$ is the ${\mathrm{A}}_\rkN$-type standard root system, see \eqref{eqn:str-A}.
We call $\parity$ {\it{the ${\mathrm{A}}(m-1,\rkN-m)$-type parity map}}
if $\parity(m):=1$ and $\parity(i):=0$ ($i\in\fkI\setminus\{m\}$).

(2) Assume that $\barN=\rkN\geq 1$ and $m\in\fkJ_{0,\rkN-1}$.
We call $(\breta^{p^-_{\rkN-m|m}},\orderbarPi)$ {\it{the ${\mathrm{B}}(m,\rkN-m)$-data}}
if $\orderbarPi$ is the ${\mathrm{B}}_\rkN$-data. 
Let $\rsystem$ be the ${\mathrm{B}}_\rkN$-type standard root system, see \eqref{eqn:str-B}.
Assume that $\barR=\rsystem\cup\{2ce_i|i\in\fkJ_{1,\rkN-m},c\in\{-1,1\}\}$.
We call $\barR$ {\it{the ${\mathrm{B}}(m,\rkN-m)$-type standard root system}}.
Note that $\barR\setminus{\frac 1 2}\barR=\rsystem$.
We call $\parity$ {\it{the ${\mathrm{B}}(m,\rkN-m)$-type parity map}}
if $\parity(\rkN-m):=1$ and $\parity(i):=0$ ($i\in\fkI\setminus\{\rkN-m\}$).

(3) Assume that $\barN=\rkN\geq 3$. 
We call $(\breta^{p^-_{1|\rkN-1}},\orderbarPi)$ {\it{the ${\mathrm{C}}(\rkN)$-data}}
if $\orderbarPi$ is the ${\mathrm{C}}_\rkN$-data.. 
Let $\rsystem$ be the ${\mathrm{C}}_\rkN$-type standard root system,
see \eqref{eqn:str-C}.
Assume that $\barR=\rsystem\setminus\{2e_1,-2e_1\}$.
We call $\barR$ {\it{the ${\mathrm{C}}(\rkN)$-type standard root system}}.
We call $\parity$ {\it{the ${\mathrm{C}}(\rkN)$-type parity map}}
if $\parity(1):=1$ and $\parity(i):=0$ ($i\in\fkI\setminus\{1\}$).

(4) Assume that $\barN=\rkN\geq 3$ and $m\in\fkJ_{2,\rkN-1}$. 
We call $(\breta^{p^-_{\barN-m|m}},\orderbarPi)$ {\it{the ${\mathrm{D}}(m,\rkN-m)$-data}} if
$\orderbarPi$ is the ${\mathrm{D}}_\rkN$-data. 
Let $\rsystem$ be the ${\mathrm{D}}_\rkN$-type standard root system,
see \eqref{eqn:str-D}.
Assume that $\barR=\rsystem\cup\{2ce_i|i\in\fkJ_{1,\rkN-m},c\in\{-1,1\}\}$.
We call $\barR$ {\it{the ${\mathrm{D}}(m,\rkN-m)$-type standard root system}}.
We call $\parity$ {\it{the ${\mathrm{D}}(m,\rkN-m)$-type parity map}}
if $\parity(\rkN-m):=1$ and $\parity(i):=0$ ($i\in\fkI\setminus\{\rkN-m\}$).

(5) Assume $\barN=\rkN=4$. 
We call $(\breta^{p^-_{1|3}},\orderbarPi)$ {\it{the ${\mathrm{F}}(4)$-data}}
if $\bal_1={\frac 1 {\sqrt{2}}}({\sqrt{3}}\brre_1-\brre_2-\brre_3-\brre_4)$,
$\bal_2={\sqrt{2}}\brre_2$, $\bal_3={\sqrt{2}}(-\brre_2+\brre_3)$, 
and $\bal_4={\sqrt{2}}(-\brre_3+\brre_4)$. 
Assume that $\barR=\barR^+\cup-\barR^+$, where
$\barR^+:=\{\,
\bal_4,\bal_3+\bal_4,\bal_2+\bal_3+\bal_4,
2\bal_2+\bal_3+\bal_4,\bal_1+\bal_2+\bal_3+\bal_4,2\bal_2+2\bal_3+\bal_4,\bal_1+2\bal_2+\bal_3+\bal_4,
\bal_1+2\bal_2+2\bal_3+\bal_4,\bal_1+3\bal_2+2\bal_3+\bal_4,2\bal_1+3\bal_2+2\bal_3+\bal_4,
\bal_1,\bal_1+\bal_2,\bal_1+\bal_2+\bal_3,\bal_1+2\bal_2+\bal_3,\bal_2,2\bal_2+\bal_3,\bal_2+\bal_3,
\bal_3\,\}$.
We call $\barR$ {\it{the ${\mathrm{F}}(4)$-type standard root system}}.
We call $\parity$ {\it{the ${\mathrm{F}}(4)$-type parity map}}
if $\parity(1):=1$ and $\parity(i):=0$ ($i\in\fkI\setminus\{1\}$).

(6) Assume that $\rkN=3$ and $\barN=4$. 
We call $(\breta^{p^-_{1|2}},\orderbarPi)$ {\it{the ${\mathrm{G}}(3)$-data}}
if $\bal_1={\sqrt{2}}\brre_1+\brre_3-\brre_4$,
$\bal_2=\brre_2-\brre_3$, $\bal_3=-2\brre_2+\brre_3+\brre_4$. 
Assume that $\barR=\barR^+\cup-\barR^+$, where
$\barR^+:=\{\,\bal_1,\,\,\bal_1+\bal_2,\,\bal_1+\bal_2+\bal_3,\,
\bal_1+2\bal_2+\bal_3,\,\bal_1+3\bal_2+\bal_3,\,\bal_1+3\bal_2+2\bal_3,\,
\bal_1+4\bal_2+2\bal_3,\,
\bal_2,\,3\bal_2+\bal_3,\,3\bal_2+2\bal_3,\,2\bal_2+\bal_3,\,\bal_2+\bal_3,\,\bal_3,\,
2\bal_1+4\bal_2+2\bal_3\,\}$.
We call $\barR$ {\it{the ${\mathrm{G}}(3)$-type standard root system}}.
Note that $\barR\cap 2\barR=\{2\bal_1+4\bal_2+2\bal_3\}$.
We call $\parity$ {\it{the ${\mathrm{G}}(3)$-type parity map}}
if $\parity(1):=1$ and $\parity(i):=0$ ($i\in\fkI\setminus\{1\}$).

(7) Assume that $\rkN=3$ and $\barN=4$. Let $x\in\bK\setminus\{\,0,\,-1\,\}$.
Fix $\sqrt{-x^2-2x}\in\bK$ so that $({\sqrt{-x^2-2x}})^2=-x^2-2x$.
We call $(\breta^{p^-_{2|2}},\orderbarPi)$ {\it{the ${\mathrm{D}}(2,1;x)$-data}}
if $\bal_1=\brre_1-\brre_2$,
$\bal_2=\brre_2-\brre_3$, $\bal_3=-x\brre_3+\sqrt{-x^2-2x}\brre_4$. 
Assume that $\barR=\barR^+\cup-\barR^+$, where
$\barR^+:=\{\,\bal_1,\,\,\bal_1+\bal_2,\,2\bal_1+\bal_2+\bal_3,\,
\bal_1+\bal_3,\,\bal_3\,\}$.
We call $\barR$ {\it{the ${\mathrm{D}}(2,1;x)$-type standard root system}}.
We call $\parity$ {\it{the ${\mathrm{D}}(2,1;x)$-type parity map}}
if $\parity(1):=1$ and $\parity(i):=0$ ($i\in\fkI\setminus\{1\}$).

(8) Let $(\breta,\orderbarPi)$ be a data as in (1)-(7) above, i.e.,
$(\breta,\orderbarPi)$ is the ${\mathrm{X}}$-data with
${\mathrm{X}}$ being 
${\mathrm{A}}(m-1,\rkN-m)$ (for some $m\in\fkJ_{1,\rkN}$ if $\rkN\geq 2$), 
${\mathrm{B}}(m,\rkN-m)$ (for some $m\in\fkJ_{0,\rkN-1}$ if $\rkN\geq 1$), 
${\mathrm{C}}(\rkN)$ (if $\rkN\geq 3$),
${\mathrm{D}}(m,\rkN-m)$ (for some $m\in\fkJ_{2,\rkN-1}$ if $\rkN\geq 3$), ${\mathrm{F}}(4)$
(if $\rkN =4$), ${\mathrm{G}}(3)$ (if $\rkN =3$), or ${\mathrm{D}}(2,1;x)$
(for some $x\in\bK\setminus\{\,0,\,-1\,\}$ if $\rkN =3$).
(Note that ${\mathrm{X}}\ne{\mathrm{A}}(0,0)$.)
We call $(\breta,\orderbarPi)$
{\it{a rank-$\rkN$ standard super-data}}.
If $\barR$ is the ${\mathrm{X}}$-type standard root system,
we call $\barR$ {\it{the standard root system associated with $(\breta,\orderbarPi)$}}.
If $\parity$ is the ${\mathrm{X}}$-type parity map,
we call $\parity$ {\it{the parity map associated with $(\breta,\orderbarPi)$}}.
}}
\end{definition}

We can directly see
\begin{lemma}\label{lemma:rtparity}
Let $(\breta,\orderbarPi)$ be 
a rank-$\rkN$ standard super-data.
Let $\barR$ be the standard root system associated with $(\breta,\orderbarPi)$.
Let $\parity$ be the parity map associated with $(\breta,\orderbarPi)$.
Define the $\bZ$-module homomorphism 
$\barparity:\oplus_{\i\in\fkI}\bZ\bal_i\to\bZ$ by
$\barparity(\bal_i):=\parity(i)$.
Let $\bal\in\barR$. Then
$\barparity(\bal)\in 2\bZ+1$ if and only if 
$\breta(\bal,\bal)=0$ or $2\bal\in\barR$.
\end{lemma}

\begin{definition}\label{definition:DefLieSup} {\rm{(1)
Let $\mfkg$ be a $\bK$-linear space.
Let $\mfkg(t)$ ($t\in\fkJ_{0,1}$) be two subspaces of $\mfkg$
with $\mfkg=\mfkg(0)\oplus\mfkg(1)$. 
Define the subspaces $\mfkg(t)$ of $\mfkg$
($t\in\bZ\setminus\fkJ_{0,1}$)
by $\mfkg(t)=\mfkg(t\pm 2)$. Let
$[\,,\,]:\mfkg\times\mfkg\to\mfkg$
be a $\bK$-bilinear map.
We call $\mfkg$ {\it{a Lie superalgebra}} if the following
$({\rm{Su}}1)$ and $({\rm{Su}}2)$ are satisfied.
\newline\par
$({\rm{Su}}1)$ $[x,y]\in\mfkg(t+t^\prime)$
and $[x,y]=-(-1)^{tt^\prime}[y,x]$ \quad ($t$, $t^\prime\in\bZ$, 
$x\in\mfkg(t)$, $y\in\mfkg(t^\prime)$).
\par
$({\rm{Su}}3)$
$[x,[y,z]]=[[x,y],z]+(-1)^{tt^\prime}[y,[x,z]]$ \quad ($t$, $t^\prime\in\bZ$, 
$x\in\mfkg(t)$, $y\in\mfkg(t^\prime)$,
$z\in\mfkg$).
\newline\par (2) Let $\mfkg$ be a Lie superalgebra.
Let $Y$ be a non-empty subset of $\mfkg$.
Let $\langle Y\rangle^{(1)}_\mfkg:=\rmSpan_\bK(Y)$.
For $t\in\fkJ_{2,\infty}$, 
let $\langle Y\rangle^{(t)}_\mfkg:=\rmSpan_\bK(\{[y,z]|y\in Y,z\in
\langle Y\rangle^{(t-1)}_\mfkg\})$.
Let $\langle Y\rangle_\mfkg:=\rmSpan_\bK(\cup_{t=1}^\infty \langle Y\rangle^{(t)}_\mfkg)$.
We call $\langle Y\rangle_\mfkg$ is {\it{the sub Lie superalgebra of $\mfkg$
generated by $Y$}}. 
\par
(3) Let $(\breta,\orderbarPi)$
be a rank-$\rkN$ standard super-data.
Let $\parity$ be the parity map associated with $(\breta,\orderbarPi)$.
In a standard way, we have the Lie superalgebra $\mfkg(\breta,\orderbarPi,\parity)$
(unique up to isomorphism)
satisfying the following conditions.
\newline\par
$({\rm{CoSu1}})$ There exist linearly independent $3\rkN$-elements
$\barH_i$, $\barE_i$, $\barF_i$ ($i\in\fkI$) of $\mfkg(\breta,\orderbarPi,\parity)$
such that $\langle \{\barH_i,\barE_i,\barF_i|i\in\fkI\}\rangle_{\mfkg(\breta,\orderbarPi,\parity)}
=\mfkg(\breta,\orderbarPi,\parity)$.
\par
$({\rm{CoSu2}})$ \quad $\barH_i\in\mfkg(\breta,\orderbarPi,\parity)(0)$,
$\barE_i$, $\barF_i\in\mfkg(\breta,\orderbarPi,\parity)(\parity(i))$
\quad ($i\in\fkI$).
\par
$({\rm{CoSu3}})$ \quad $[\barH_i,\barH_j]=0$, $[\barH_i,\barE_j]=\breta(\bal_i,\bal_j)\barE_j$, 
$[\barH_i,\barF_j]=-\breta(\bal_i,\bal_j)\barF_j$, $[\barE_i,\barF_j]=\delta_{ij}\barH_i$
\quad ($i,j\in\fkI$).
\par
$({\rm{CoSu4}})$ There exist subspaces $\mfkg(\breta,\orderbarPi,\parity)_\bal$
($\bal\in\bZ\barPi$) of $\mfkg(\breta,\orderbarPi,\parity)$
such that $\mfkg(\breta,\orderbarPi,\parity)=\oplus_{\bal\in\bZ\barPi}\mfkg(\breta,\orderbarPi,\parity)_\bal$,
$[\mfkg(\breta,\orderbarPi,\parity)_\bbeta,\mfkg(\breta,\orderbarPi,\parity)_{\bbeta^\prime}]
\subseteq\mfkg(\breta,\orderbarPi,\parity)_{\bbeta+\bbeta^\prime}$
($\bbeta,\bbeta^\prime\in\bZ\barPi$), 
$\mfkg(\breta,\orderbarPi,\parity)_0=\oplus_{i\in\fkI}\bK\barH_i$, and
$\mfkg(\breta,\orderbarPi,\parity)_{\bal_j}=\bK\barE_j$, 
$\mfkg(\breta,\orderbarPi,\parity)_{-\bal_j}=\bK\barF_j$
($j\in\fkI$).
\par
$({\rm{CoSu5}})$ Let $\bal\in\bZgeqo\barPi\setminus(\barPi\cup\{0\})$.
For $X\in\mfkg(\breta,\orderbarPi,\parity)_\bal$, if 
$[X,\barF_i]=0$ for all $i\in\fkI$, then $X=0$.
For $Y\in\mfkg(\breta,\orderbarPi,\parity)_{-\bal}$, if 
$[\barE_i,Y]=0$ for all $i\in\fkI$, then $Y=0$.
(Note that the conditions $({\rm{CoSu1}})$-$({\rm{CoSu4}})$ imply that 
$\mfkg(\breta,\orderbarPi,\parity)_\bbeta=\{0\}$
for $\bbeta\in\bZ\barPi\setminus(\bZgeqo\barPi\cup -\bZgeqo\barPi)$.)
}}
\end{definition}

We call $\mfkg(\breta,\orderbarPi,\parity)$ the
{\it{contragredient Lie superalgebra}}.

It is well-known that 
\begin{theorem}\label{theorem:multisuper}
{\rm{(}}c.f. {\rm{\cite[Proposition~2.5.5]{Kac77}}}{\rm{)}} 
Let $(\breta,\orderbarPi)$
be a rank-$\rkN$ standard super-data.
Let $\barR$ and $\parity$ be 
the standard root system and the parity map associated with $(\breta,\orderbarPi)$
respectively.
Then for $\bbeta\in\bZ\barPi$,
$\dim\mfkg(\breta,\orderbarPi,\parity)_\bbeta\geq 1$
if and only if $\bbeta\in\barR\cup\{0\}$.
Moreover, for $\bbeta\in\barR$, we have $\dim\mfkg(\breta,\orderbarPi,\parity)_\bbeta =1$.
\end{theorem}

Let $(\breta,\orderbarPi)$
be a rank-$\rkN$ standard super-data. Let $\parity$ be 
the parity map associated with $(\breta,\orderbarPi)$.
Let $\mfkg:=\mfkg(\breta,\orderbarPi,\parity)$.
Following \cite{Kac77}, we give names as follows.
Assume that $(\breta,\orderbarPi)$ is the ${\mathrm{X}}$-data.
If ${\mathrm{X}}={\mathrm{A}}(m-1,\rkN-m)$ and $2m-1=\rkN\geq 2$,
then $\mfkg$ has the one-dimensional center ${\mathfrak{c}}$,
$\mfkg/{\mathfrak{c}}$ is a simple Lie superalgebra
and $\mfkg/{\mathfrak{c}}$ is called ${\mathrm{A}}(m-1,m-1)$.
Otherwise $\mfkg$ is a simple Lie superalgebra
and called ${\mathrm{X}}$.
If ${\mathrm{X}}$ is ${\mathrm{A}}(m-1,\rkN-m)$
(resp. ${\mathrm{B}}(m,\rkN-m)$, resp. ${\mathrm{C}}(\rkN)$,
resp. ${\mathrm{D}}(m,\rkN-m)$), then
$\mfkg$ is also called 
${\mathrm{sl}}(m|\rkN-m+1)$
(resp. ${\mathrm{osp}}(2m+1|2(\rkN-m))$, resp. ${\mathrm{osp}}(2|2(\rkN-1))$,
resp. ${\mathrm{osp}}(2m|2(\rkN-m))$).

Let $\barN\in\fkJ_{2,\infty}$.
For $i\in\fkJ_{1,\barN-1}$, define the bijection $\perm^{(\barN)}_i:\fkJ_{1,\barN}\to\fkJ_{1,\barN}$ by 
\begin{equation}\label{eqn:defsbarNi}
\perm^{(\barN)}_i(j):=
\left\{\begin{array}{ll}
i+1 & \quad\mbox{if $j=i$}, \\
i & \quad\mbox{if $j=i+1$}, \\
j & \quad\mbox{if $j\in\fkJ_{1,\barN}\setminus\fkJ_{i,i+1}$}.
\end{array}\right.
\end{equation}

\subsection{A Weyl groupoid of the simple Lie superalgebra 
${\mathrm{A}}(m-1,\rkN-m)$}
\label{subsection:LGEsuperA}

Let $m\in\fkJ_{1,\rkN}$.
For $i\in\fkI$, define the bijection 
$\superAnewtr_i:\mcA_{m|\rkN+1-m}\to\mcA_{m|\rkN+1-m}$ by
$\superAnewtr_i(p):=p\circ \perm^{(\rkN+1)}_i$.
Let $\orderrsysbase=(\hatal_i|i\in\fkI)$  and
$\rsystem$ 
be the ${\mathrm{A}}_\rkN$-data
and the ${\mathrm{A}}_\rkN$-type
standard root system respectively.
Assume that $\orderedPi=\orderrsysbase$,
so $\bV=\oplus_{i\in\fkI}\bR\hatal_i\subset\bR^{\rkN+1}$.
For any $p\in\mcA_{m|\rkN+1-m}$,
let $\superAR(p):=\rsystem$, and
define the generalized Cartan Matrix $\superAC^p=(\superAc^p_{ij})_{i,j\in\fkI}$ by
$\superAc^p_{ij}={\frac {2\hateta(\hatal_i,\hatal_j)} {\hateta(\hatal_i,\hatal_i)}}$.
We easily see that $\superAmcC_{m|\rkN+1-m}:=\mcC(\fkI,\mcA_{m|\rkN+1-m},(\superAnewtr_i)_{i\in\fkI},
(\superAC^p)_{p\in\mcA_{m|\rkN+1-m}})$ is a Cartan scheme,
and that $\superAmcR_{m|\rkN+1-m}:=\mcR(\superAmcC,(\superAR(p))_{p\in\mcA_{m|\rkN+1-m}})$
is a generalized root system of type $\superAmcC_{m|\rkN+1-m}$. 
In particular, we have
\begin{lemma} \label{lemma:pre-longest-superA-id}
Let $\orderrsysbase=(\hatal_i|i\in\fkI)$ be the ${\mathrm{A}}_\rkN$-data.
Then
for any $m\in\fkJ_{1,\rkN}$ and any $p\in\mcA_{m|\rkN+1-m}$,
$(\superAmcR_{m|\rkN+1-m},p)$ and $(\mcR_\orderrsysbase,a)$
are quasi-isomorphic,
where $a$ is the only element of $\mcA_\orderrsysbase$.
\end{lemma}

Note that
for the generalized root system $\superAmcR_{m|\rkN+1-m}$, we have
\begin{equation}\label{eqn:inveta-superA}
\breta^p(1^p\tils_{f,r}(x),1^p\tils_{f,r}(y))=\breta^{p_{f,r}}(x,y)\quad(x,y\in\bV)
\end{equation} for $p\in\mcA_{m|\rkN+1-m}$, $f\in\funcI_\infty$ and $r\in\bN$.

Recall that if $p:=p^+_{m|\rkN+1-m}\in\mcA_{m|\rkN+1-m}$,
then $(\breta^p,\orderrsysbase)$ is the ${\mathrm{A}}(m-1,\rkN-m)$-data
and $\rsystem$ is the standard root system associated with $(\breta^p,\orderrsysbase)$.

We have
\begin{lemma}\label{lemma:longest-superA-id}
Let $m\in\fkI$, and treat the Weyl groupoid $\mcW(\superAmcR_{m|\rkN+1-m})$.
Let $p:=p^+_{m|\rkN+1-m}\in\mcA_{m|\rkN+1-m}$.
Let $n:={\frac {\rkN(\rkN+1)} 2}$ and $r:=n-m(\rkN-m+1)$.
Then $\tilell^p(1^pw_0)=n$.
Let $f\in\funcI_n$ be as in Proposition~{\rm{\ref{proposition:ANNlgest}}}.
Let $\beta_x:=1^p\tils_{f,x-1}(\al_x)$
$(x\in\fkJ_{1,n})$.
Then $1^pw_0=1^p\tils_{f,n}$
and $\superARp(p)=\{\beta_x|x\in\fkJ_{1,n}\}$. Moreover
\begin{equation}\label{eqn:lsA-id-1}
\breta^p(\beta_k,\beta_k)\subseteq\{-2,2\}\,\,
(k\in\fkJ_{1,r})\quad\mbox{and}\quad 
\breta^p(\beta_t,\beta_t)=0\,\,
(t\in\fkJ_{r+1,n}).
\end{equation}
\end{lemma}
{\it{Proof}}. This follows from Lemmas~\ref{lemma:techlm} and \ref{lemma:pre-longest-superA-id}
and \eqref{eqn:ellRp}, \eqref{eqn:proflgth-d}, \eqref{eqn:ctprpoflg-b}, 
\eqref{eqn:ANNlgest-a} and \eqref{eqn:ANNlgest-b}. 
\hfill $\Box$

\subsection{The longest element of a Weyl groupoid of the simple Lie superalgebra 
${\mathrm{B}}(m,\rkN-m)$ with $m\geq 1$}
\label{subsection:LGEsuperB}

Let $m\in\fkJ_{1,\rkN}$.
For $i\in\fkI$, define the bijection 
$\superBnewtr_i:\mcA_{m|\rkN-m}\to\mcA_{m|\rkN-m}$ by
\begin{equation*}
\superBnewtr_i(p):=
\left\{\begin{array}{ll}
p\circ \perm^{(\rkN)}_i & \quad\mbox{if $i\in\fkJ_{1,\rkN-1}$}, \\
p & \quad\mbox{if $i=\rkN$}.
\end{array}\right.
\end{equation*} Let $\orderrsysbase=(\hatal_i|i\in\fkI)$  and
$\rsystem$ 
be the ${\mathrm{B}}_\rkN$-data
and the ${\mathrm{B}}_\rkN$-type
standard root system respectively.
Assume that $\orderedPi=\orderrsysbase$,
so $\bV=\oplus_{i\in\fkI}\bR\hatal_i=\truebRrkN$.
For any $p\in\mcA_{m|\rkN-m}$,
let $\superBR(p):=\rsystem$, and
define the generalized Cartan Matrix $\superBC^p=(\superBc^p_{ij})_{i,j\in\fkI}$ by
$\superBc^p_{ij}={\frac {2\hateta(\hatal_i,\hatal_j)} {\hateta(\hatal_i,\hatal_i)}}$.
We see that $\superBmcC_{m|\rkN-m}:=\mcC(\fkI,\mcA_{m|\rkN-m},(\superBnewtr_i)_{i\in\fkI},
(\superBC^p)_{p\in\mcA_{m|\rkN-m}})$ is a Cartan scheme,
and that $\superBmcR_{m|\rkN-m}:=\mcR(\superBmcC,(\superBR(p))_{p\in\mcA_{m|\rkN-m}})$
is a generalized root system of type $\superBmcC_{m|\rkN-m}$. 
In particular, we have
\begin{lemma} \label{lemma:pre-longest-superB-id}
Let $\orderrsysbase=(\hatal_i|i\in\fkI)$ be the ${\mathrm{B}}_\rkN$-data.
Then
for any $m\in\fkJ_{1,\rkN}$ and any $p\in\mcA_{m|\rkN-m}$,
$(\superBmcR_{m|\rkN-m},p)$ and $(\mcR_\orderrsysbase,a)$
are quasi-isomorphic,
where $a$ is the only element of $\mcA_\orderrsysbase$.
\end{lemma}

Note that
for the generalized root system $\superBmcR_{m|\rkN+1-m}$, we have
\begin{equation}\label{eqn:inveta-superB}
\breta^p(1^p\tils_{f,r}(x),1^p\tils_{f,r}(y))=\breta^{p_{f,r}}(x,y)\quad(x,y\in\bV)
\end{equation} for $p\in\mcA_{m|\rkN-m}$, $f\in\funcI_\infty$ and $r\in\bN$.

Define $\superBf_{m|\rkN-m}\in\funcI_{\rkN^2}$ by
\begin{equation}\label{eqn:superBf-def}
\superBf_{m|\rkN-m}(t):=
\left\{\begin{array}{ll}
\rkN-m+t & \quad\mbox{if $t\in\fkJ_{1,m}$}, \\
\superBf_{m|\rkN-m}(t-m) & \quad\mbox{if $t\in\fkJ_{m+1,m^2}$}, \\
t-m^2 & \quad\mbox{if $t\in\fkJ_{m^2+1,m^2+\rkN}$}, \\
\rkN-(t-(m^2+\rkN)) & \quad\mbox{if $t\in\fkJ_{m^2+\rkN+1,m^2+\rkN+m}$}, \\
\superBf_{m|\rkN-m}(t-(\rkN+m)) & \quad\mbox{if $t\in\fkJ_{m^2+\rkN+m+1,\rkN^2}$}.
\end{array}\right.
\end{equation}

Let $\superBp_{m|\rkN-m}:=p^-_{\rkN-m|m}\in\mcA_{m|\rkN-m}$.
Note that if $p=\superBp_{m|\rkN-m}$, then
$(\breta^p,\orderrsysbase)$ is the ${\mathrm{B}}(m,\rkN-m)$-data
and that if $\barR$ is the standard root system associated with $(\breta^p,\orderrsysbase)$,
then $\rsystem=\barR\setminus 2\barR$.

Let $\superBp_{m|\rkN-m}^{(0)}:=\superBp_{m|\rkN-m}$,
and for $k\in\fkJ_{1,m}$,
let $\superBp_{m|\rkN-m}^{(k)}:=\superBp_{m|\rkN-m}^{(k-1)}\circ\perm^{(\rkN)}_{\rkN-m+k-1}$.
Then for $t\in\fkJ_{1,\rkN^2}$, we have
\begin{equation}\label{eqn:superBp-fom}
\begin{array}{l}
(\superBp_{m|\rkN-m})_{\superBf_{m|\rkN-m},t} \\ \quad
=\left\{\begin{array}{ll}
\superBp_{m|\rkN-m} & \quad\mbox{if $t\in\fkJ_{1,m^2+\rkN-m-1}$}, \\
\superBp_{m|\rkN-m}^{(t-(m^2+\rkN-m-1))} & \quad\mbox{if $t\in\fkJ_{m^2+\rkN-m,m^2+\rkN-1}$}, \\
\superBp_{m|\rkN-m}^{(m-(t-(m^2+\rkN)))} & \quad\mbox{if $t\in\fkJ_{m^2+\rkN,m^2+\rkN+m}$}, \\
\superBp_{m|\rkN-m}^{(t-(\rkN+m))} & \quad\mbox{if $t\in\fkJ_{m^2+\rkN+m+1,\rkN^2}$}. 
\end{array}\right.
\end{array}
\end{equation}

We have
\begin{lemma}\label{lemma:longest-superB-id}
Let $m\in\fkI$, and treat the Weyl groupoid $\mcW(\superBmcR_{m|\rkN-m})$.
Let $n:=\rkN^2$, $f:=\superBf_{m|\rkN-m}\in\funcI_n$
and $p:=\superBp_{m|\rkN-m}\in\mcA_{m|\rkN-m}$.
Then 
\begin{equation}\label{eqn:longest-superB-id-1}
\tilell^p(1^pw_0)=n\quad\mbox{and}\quad 1^pw_0=1^p\tils_{f,n}.
\end{equation}
Moreover for $t\in\fkJ_{1,n}$, letting
$\beta_t:=1^p\tils_{f,t-1}(\al_t)$, we have
$\superBRp=\{\beta_t|t\in\fkJ_{1,n}\}$ and
\begin{equation}\label{eqn:longest-superB-id-2}
\breta^p(\beta_t,\beta_t)
= 
\left\{\begin{array}{ll}
2 & \quad\mbox{if $t\in\fkJ_{1,m-1}$}, \\
1 & \quad\mbox{if $t=m$}, \\
\breta^p(\beta_{t-m},\beta_{t-m}) &
\quad\mbox{if $t\in\fkJ_{m+1,m^2}$}, \\
-2 & \quad\mbox{if $t\in\fkJ_{m^2+1,m^2+\rkN-m-1}$}, \\
0 & \quad\mbox{if $t\in\fkJ_{m^2+\rkN-m,m^2+\rkN-1}$}, \\
-1 & \quad\mbox{if $t=m^2+\rkN$}, \\
0 & \quad\mbox{if $t\in\fkJ_{m^2+\rkN+1,m^2+\rkN+m}$}, \\
\breta^p(\beta_{t-(\rkN+m)},\beta_{t-(\rkN+m)}) &
\quad\mbox{if $t\in\fkJ_{m^2+\rkN+m+1,\rkN^2}$}.
\end{array}\right.
\end{equation}
\end{lemma}
{\it{Proof}}. We have \eqref{eqn:longest-superB-id-1} 
by Lemmas~\ref{lemma:techlm} and \ref{lemma:pre-longest-superB-id}
and
\eqref{eqn:proflgth-d}, \eqref{eqn:ctprpoflg-b}, 
\eqref{eqn:srBm},
\eqref{eqn:srBm-2}. 
Then we can directly see \eqref{eqn:longest-superB-id-2} by 
\eqref{eqn:longest-superB-id-1},
\eqref{eqn:ellRp},
\eqref{eqn:inveta-superB},
\eqref{eqn:superBf-def} and \eqref{eqn:superBp-fom}.
\hfill $\Box$

\subsection{Longest elements of Weyl groupoids of the simple Lie superalgebra ${\mathrm{D}}(m,\rkN-m)$}
\label{subsection:LGEsuperD}
In this subsection, assume $\rkN\geq 3$.
Let $\ospefkI:=\fkI\cup\{\rkN+1\}=\fkJ_{1,\rkN+1}$.

Let $m\in\fkJ_{1,\rkN-1}$.
Let $\ospemfD_{m|\rkN-m}$ be the set of
all maps $\ospep:\ospefkI\to\fkJ_{0,1}$ 
satisfying the following conditions $(p1)$ and $(p2)$.
\newline\par
$(p1)$ $\sum_{i\in\fkI}\ospep(i)=m$. \par
$(p2)$ If $\ospep(\rkN)=0$, then $\ospep(\rkN+1)=0$. 
\newline\par

For $i\in\fkI$, define the bijection 
\begin{equation*}
\ospenewtr_i:\ospemfD_{m|\rkN-m}\to\ospemfD_{m|\rkN-m}
\end{equation*}
by the following $(\ospenewtr\mbox{-}1)$-$(\ospenewtr\mbox{-}3)$,
where see \eqref{eqn:defsbarNi} for $\perm^{(\rkN+1)}_i$.
\newline\par
$(\ospenewtr\mbox{-}1)$ For
$i\in\fkJ_{1,\rkN-2}$, let $\ospenewtr_i(\ospep):=\ospep\circ \perm^{(\rkN+1)}_i$. \par
$(\ospenewtr\mbox{-}2)$ Assume $i=\rkN-1$.
If $\ospep(\rkN+1)=0$, let $\ospenewtr_i(\ospep):=\ospep\circ \perm^{(\rkN+1)}_{\rkN-1}$.
If $\ospep(\rkN+1)=1$, let $\ospenewtr_i(\ospep):=\ospep$. \par 
$(\ospenewtr\mbox{-}3)$ Assume $i=\rkN$.
If $\ospep(\rkN-1)=\ospep(\rkN)$, let $\ospenewtr_i(\ospep):=\ospep$.
If $(\ospep(\rkN-1),\ospep(\rkN),\ospep(\rkN+1))$ equals $(0,1,0)$
(resp. $(0,1,1)$, resp. $(1,0,0)$), let $\ospenewtr_i(\ospep)$
be such that $\ospenewtr_i(\ospep)_{|\fkJ_{1,\rkN-2}}:=\ospep_{|\fkJ_{1,\rkN-2}}$ and
$(\ospenewtr_i(\ospep(\rkN-1),\ospenewtr_i(\ospep(\rkN)),\ospenewtr_i(\ospep(\rkN+1)))$
 equals $(0,1,0)$
(resp. $(1,0,0)$, resp. $(0,1,1)$).
\newline\newline
All $\ospep\in\ospemfD_{m|\rkN-m}$ with $\rkN=4$ and $m=2$ are given in 
Figure~\ref{fig:DynkinDall} below.
\newline\par
Let $\ospep\in\ospemfD_{m|\rkN-m}$.
Define the $\bR$-bilinear map $\ospeeta^\ospep:\truebRrkN\times\truebRrkN\to\bR$
by $\ospeeta^\ospep(\vece_i,\vece_j):=\delta(\ospep(i),\ospep(j))(-1)^{\ospep(i)}$
($i,j\in\fkI$). 
Define the subset $\ospeR(\ospep)$ of $\truebRrkN$ by
\begin{equation} \label{eqn:defoseR}
\begin{array}{lcl}
\ospeR(\ospep)&:=&\{\,c\vece_i+c^\prime\vece_j\,|\,i,j\in\fkI,\,i\ne j,\,
c,\,c^\prime\in\{\,1,\,-1\,\}\,\} \\
& & \quad\cup\{\,c^{\prime\prime}\vece_i\,|\,i\in\fkI,\,\ospep(i)=1,\,
c^{\prime\prime}\in\{\,1,\,-1\,\}\,\}.
\end{array}
\end{equation}
For
$i\in\fkI$, define $\ospeal^\ospep_i\in\truebRrkN$ by
\begin{equation} \label{eqn:defospeal}
\ospeal^\ospep_i:=
\left\{
\begin{array}{ll}
\vece_i-\vece_{i+1} & \mbox{if $i\in\fkJ_{1,\rkN-2}$}, \\
\vece_{\rkN-1}-\vece_\rkN  & \mbox{if $i=\rkN-1$ and $\ospep(\rkN+1)=0$}, \\
-2\vece_\rkN & \mbox{if $i=\rkN-1$ and $\ospep(\rkN+1)=1$}, \\
\vece_{\rkN-1}+\vece_\rkN  & \mbox{if $i=\rkN$ and $\ospep(\rkN)=\ospep(\rkN+1)$}, \\
2\vece_\rkN & \mbox{if $i=\rkN$ and $\ospep(\rkN)=1$, $\ospep(\rkN+1)=0$}. 
\end{array}\right.
\end{equation} We can directly see that
\begin{equation} \label{eqn:defospeeta}
\ospeeta^\ospep(x,y)=\ospeeta^{\ospenewtr_i\ospep}(\hats_{\ospeal^\ospep_i}(x),\hats_{\ospeal^\ospep_i}(y))
\quad(i\in\fkI,\,x,y\in\truebRrkN).
\end{equation}
Define the generalized Cartan matrix
$\ospeC^\ospep=[\ospec^\ospep_{ij}]_{i,j\in\fkI}
\in{\mathrm{M}}_\rkN(\bZ)$ by
\begin{equation*}
\ospec^\ospep_{ij}:=
\left\{
\begin{array}{ll}
2 & \quad\mbox{if $i=j$}, \\ 
{\frac {2\ospeeta^\ospep(\ospeal^\ospep_i,\ospeal^\ospep_j)} {\ospeeta^\ospep(\ospeal^\ospep_i,\ospeal^\ospep_i)}} 
& \quad\mbox{if $i\ne j$ and $\ospeeta^\ospep(\ospeal^\ospep_i,\ospeal^\ospep_i)\ne 0$}, \\
0 & \quad\mbox{if $i\ne j$ and $\ospeeta^\ospep(\ospeal^\ospep_i,\ospeal^\ospep_i)=\ospeeta^\ospep(\ospeal^\ospep_i,\ospeal^\ospep_j)=0$},\\
-1 & \quad\mbox{if $i\ne j$, $\ospeeta^\ospep(\ospeal^\ospep_i,\ospeal^\ospep_i)=0$
and $\ospeeta^\ospep(\ospeal^\ospep_i,\ospeal^\ospep_j)\ne 0$}.
\end{array}\right.
\end{equation*} Define the $\bR$-linear homomorphism 
$\ospexi^\ospep:\truebRrkN\to\bV$
by $\ospexi^\ospep(\ospeal^\ospep_i):=\al_i$
($i\in\fkI$).

We can directly see Proposition~\ref{proposition:evenosp-1} below.

\begin{proposition} \label{proposition:evenosp-1}
{\rm{(1)}} $\ospemcC_{m|\rkN-m}:=
{\mcC}(\fkI,\ospemfD_{m|\rkN-m},(\ospenewtr_i)_{i\in\fkI},(\ospeC^\ospep)_{\ospep\in\ospemfD_{m|\rkN-m}})$
is a Cartan scheme. \par
{\rm{(2)}} There exists a unique 
generalized root system 
\begin{equation*}
\ospemcR_{m|\rkN-m}:=\mcR(\ospemcC_{m|\rkN-m},(\tR(\ospep))_{\ospep\in\ospemfD_{m|\rkN-m}})
\end{equation*} 
of type $\ospemcC_{m|\rkN-m}$
such that 
\begin{equation} \label{eqn:evenosp-1-1}
\ospexi^\ospep(\ospeR(\ospep))=\tR(\ospep)
\quad (\ospep\in\ospemfD_{m|\rkN-m}). 
\end{equation} Moreover
\begin{equation} \label{eqn:evenosp-1-2}
\ospexi^{\ospenewtr_i\ospep}\circ \hats_{\ospeal^\ospep_i}=\tils^\ospep_i \circ \ospexi^\ospep
\quad(\ospep\in\ospemfD_{m|\rkN-m},\,i\in\fkI).
\end{equation} 
\end{proposition}

Define 
$\ospep_{m|\rkN-m}\in\ospemfD_{m|\rkN-m}$ by
\begin{equation*}
\ospep_{m|\rkN-m}(i):=
\left\{
\begin{array}{ll}
1 & \quad\mbox{if $i\in\fkJ_{1,\rkN-m}$}, \\
0 & \quad\mbox{if $i\in\fkJ_{\rkN-m+1,\rkN+1}$}. 
\end{array}\right.
\end{equation*} 
Note that
$(\ospeeta^{\ospep_{m|\rkN-m}},(\ospeal^{\ospep_{m|\rkN-m}}_i|i\in\fkI))$
is the ${\mathrm{D}}(m,\rkN-m)$-data
and $\ospeR(\ospeeta^{\ospep_{m|\rkN-m}})$ is the
${\mathrm{D}}(m,\rkN-m)$ standard root system
associated with $(\ospeeta^{\ospep_{m|\rkN-m}},(\ospeal^{\ospep_{m|\rkN-m}}_i|i\in\fkI))$. 

For $k\in\fkJ_{0,m}$, define 
$\ospep^{(k)}_{m|\rkN-m}\in\ospemfD_{m|\rkN-m}$ by
\begin{equation*}
\ospep^{(k)}_{m|\rkN-m}(i):=
\left\{
\begin{array}{ll}
0 & \quad\mbox{if $i\in\fkJ_{1,\rkN-m-1}\cup\{\rkN-m+k\}\cup\{\rkN+1\}$}, \\
1 & \quad\mbox{if $i\in\fkJ_{\rkN-m,\rkN-m+k-1}\cup\fkJ_{\rkN-m+k+1,\rkN}$}, 
\end{array}\right.
\end{equation*} so 
\begin{equation}\label{eqn:ospepk}
\ospep_{m|\rkN-m}=\ospep^{(0)}_{m|\rkN-m}
\quad\mbox{and}\quad\ospep^{(k)}_{m|\rkN-m}=\ospep^{(k-1)}_{m|\rkN-m}\circ \perm^{(\rkN+1)}_{\rkN-m-1+k} 
\,\,(k\in\fkJ_{1,m}).
\end{equation}
By Lemma~\ref{lemma:lgslm}~(1), \eqref{eqn:defoseR} and \eqref{eqn:evenosp-1-1},
\begin{equation} \label{eqn:evenosp-1-3}
\tilell_\ospep(1^\ospep w_0)=|\tRp(\ospep)|=\rkN^2-m\quad(\ospep\in\ospemfD_{m|\rkN-m}).
\end{equation}
Define $\ospef_{m|\rkN-m}\in\funcI_{\rkN^2-m}$ by 
\begin{equation}\label{eqn:defddospef}
\begin{array}{l}
\ospef_{m|\rkN-m}(t) \\
\,\,:=
\left\{
\begin{array}{ll}
\rkN-m+t & \quad\mbox{if $t\in\fkJ_{1,m}$}, \\
\ospef_{m|\rkN-m}(t-m) & \quad\mbox{if $t\in\fkJ_{m+1,m(m-1)}$}, \\
t-m(m-1) & \quad\mbox{if $t\in\fkJ_{m(m-1)+1,m(m-1)+\rkN}$}, \\
\rkN-(t-(m(m-1)+\rkN)) & \quad\mbox{if $t\in\fkJ_{m(m-1)+\rkN+1,m(m-1)+\rkN+m}$}, \\
\ospef_{m|\rkN-m}(t-(\rkN+m)) & \quad\mbox{if $t\in\fkJ_{\rkN+m^2+1,\rkN^2-m}$}.
\end{array}\right.
\end{array}
\end{equation}
Using \eqref{eqn:ospepk}, we can see that
for $t\in\fkJ_{1,\rkN^2-m}$, 
\begin{equation} \label{eqn:evenosp-1-4}
\begin{array}{l}
(\ospep_{m|\rkN-m})_{\ospef_{m|\rkN-m},t}
= \\ \,\, \left\{
\begin{array}{ll}
\ospep_{m|\rkN-m} & \,\,\mbox{if $t\in\fkJ_{1,m(m-1)+\rkN-m-1}$}, \\
\ospep_{m|\rkN-m}^{(t-(m(m-1)+\rkN-m-1))} & \,\,\mbox{if $t\in\fkJ_{m(m-1)+\rkN-m,m(m-1)+\rkN-1}$}, \\
\ospep_{m|\rkN-m}^{(m-(t-(m(m-1)+\rkN))} & \,\,\mbox{if $t\in\fkJ_{m(m-1)+\rkN,m(m-1)+\rkN+m}$}, \\
\ospep_{m|\rkN-m}^{(t-(\rkN+m))} & \,\,\mbox{if $t\in\fkJ_{\rkN+m^2+1,\rkN^2-m}$}.
\end{array}\right.
\end{array}
\end{equation} 

Note that if $m\in\fkJ_{2,\rkN-1}$, $p=\ospep_{m|\rkN-m}$
and $\superordDbase=(\ospeal^p_i|i\in\fkI)$, then
$(\ospeeta^p,\superordDbase)$
is the ${\mathrm{D}}(m,\rkN-m)$-data and 
$\ospeR(p)$ is the standard root system associated with 
$(\ospeeta^p,\superordDbase)$.
\begin{proposition}\label{proposition:EVENosp-lngst} 
Let $p:=\ospep_{m|\rkN-m}\in\ospemfD_{m|\rkN-m}$.
Let $n:=\rkN^2-m$ and $f:=\ospef_{m|\rkN-m}\in\funcI_n$.
Then we have
\begin{equation} \label{eqn:evenosp-lngst}
\tilell_\ospep(1^p w_0)=\rkN^2-m,\quad 1^pw_0=1^p\tils_{f,n}
\end{equation}
and
\begin{equation}\label{eqn:evenosp-lngstd}
p_{f,n}=p.
\end{equation} Moreover, letting $r:=m(m-1)$, we have
\begin{equation}\label{eqn:longest-superD-id}
\begin{array}{l}
\ospeeta^{p_{f,t-1}}(\ospeal_t,\ospeal_t)
= \\
\quad
\left\{\begin{array}{ll}
2 & \quad\mbox{if $t\in\fkJ_{1,m}$}, \\
\ospeeta^{p_{f,t-m-1}}(\ospeal_{t-m},\ospeal_{t-m}) &
\quad\mbox{if $t\in\fkJ_{m+1,r}$}, \\
2 & \quad\mbox{if $t\in\fkJ_{r+1,r+\rkN-m-1}$}, \\
0 & \quad\mbox{if $t\in\fkJ_{r+\rkN-m,r+\rkN-1}$}, \\
-4 & \quad\mbox{if $t=m^2+\rkN$}, \\
0 & \quad\mbox{if $t\in\fkJ_{r+\rkN+1,r+\rkN+m}$}, \\
\ospeeta^{p_{f,t-(\rkN+m)-1}}(\ospeal_{t-(\rkN+m)},\ospeal_{t-(\rkN+m)}) &
\quad\mbox{if $t\in\fkJ_{r+\rkN+m+1,n}$}.
\end{array}\right.
\end{array}
\end{equation}
\end{proposition}
{\it{Proof.}} Let $r:=m(m-1)$ and $r^\prime:=n-r$.
Define $f^\prime\in\funcI_{r^\prime}$ by $f^\prime(y):=f(y+r)$
($y\in\fkJ_{1,r^\prime}$).
By \eqref{eqn:evenosp-1-4}, $p=p_{f,r}=p_{f,n}$, since $r^\prime=(\rkN+m)(\rkN-m)$. 
Hence we have \eqref{eqn:evenosp-lngstd} and 
\begin{equation}\label{eqn:EVENosp-lngst-(-1)}
1^p\tils_{f,n}=1^p\tils_{f,r}1^p\tils_{f^\prime,r^\prime}.
\end{equation}
For $t\in\fkJ_{1,r}$, since $f(t)\in\fkJ_{\rkN-m+1,\rkN}$ for $t\in\fkJ_{1,r}$, 
by \eqref{eqn:defospeal} and \eqref{eqn:evenosp-1-4}, we have
\begin{equation*}
\ospeal^{p_{f,t}}_{f(t)}=
\left\{\begin{array}{ll}
\vece_{f(t)}-\vece_{f(t)+1} & \mbox{if $f(t)\in\fkJ_{\rkN-m+1,\rkN-1}$}, \\
\vece_{f(t)-1}+\vece_{f(t)} & \mbox{if $f(t)=\rkN$}.
\end{array}\right.
\end{equation*}
Hence
by \eqref{eqn:hatsrsN} and \eqref{eqn:evenosp-1-2},
we have 
\begin{equation}\label{eqn:EVENosp-lngst-1}
(\ospexi^p)^{-1}\circ 1^p\tils_{f,r}\circ\ospexi^p=
\stProj_{\fkJ_{1,\rkN-m}}-\stProj_{\fkJ_{\rkN-m+1,\rkN-1}}-(-1)^m\stProj_{\fkJ_{\rkN,\rkN}}.
\end{equation}

For $t^\prime\in\fkJ_{1,r^\prime}$, we can directly see
\begin{equation}\label{eqn:EVENosp-lngst-2}
\ospeal^{p_{f^\prime,t^\prime}}_{f^\prime(t^\prime)}=
\left\{\begin{array}{ll} 
\vece_{f^\prime(t^\prime)}-\vece_{f^\prime(t^\prime)+1} 
&\quad\mbox{if $f^\prime(t^\prime)\in\fkJ_{1,\rkN-1}$}, \\
2\vece_\rkN &\quad\mbox{if $f^\prime(t^\prime)=\rkN$}. 
\end{array}\right.
\end{equation}
By \eqref{eqn:EXShatb},
\eqref{eqn:evenosp-1-2} and \eqref{eqn:EVENosp-lngst-2}, we have 
\begin{equation}\label{eqn:EVENosp-lngst-3}
(\ospexi^p)^{-1}\circ 1^p\tils_{f^\prime,r^\prime}\circ\ospexi^p=
-\stProj_{\fkJ_{1,\rkN-m}}+\stProj_{\fkJ_{\rkN-m+1,\rkN}}.
\end{equation}
By \eqref {eqn:EVENosp-lngst-(-1)}, \eqref{eqn:EVENosp-lngst-1}
and \eqref{eqn:EVENosp-lngst-3}, we have 
\begin{equation}\label{eqn:EVENosp-lngst-4}
(\ospexi^p)^{-1}\circ 1^p\tils_{f,n}\circ\ospexi^p=
-\stProj_{\fkJ_{1,\rkN-1}}-(-1)^m\stProj_{\fkJ_{\rkN,\rkN}}.
\end{equation}
By \eqref{eqn:defospeal}
and \eqref{eqn:EVENosp-lngst-4}, we have
\begin{equation*}
1^p\tils_{f,n}(\al_i)=
\left\{
\begin{array}{ll}
-\al_i & \quad\mbox{if $i\in\fkJ_{1,\rkN-2}$}, \\
-\al_i & \quad\mbox{if $m\in2\bN$ and $i\in\fkJ_{\rkN-1,\rkN}$}, \\
-\al_\rkN & \quad\mbox{if $m\in2\bN-1$ and $i=\rkN-1$}, \\
-\al_{\rkN-1} & \quad\mbox{if $m\in2\bN-1$ and $i=\rkN$}.
\end{array}\right.
\end{equation*}
Hence
$1^p\tils_{f,n}
(\Pi)=-\Pi$. By \eqref{eqn:lgsta} and \eqref{eqn:evenosp-1-3}, 
we have \eqref{eqn:evenosp-lngst}.
We can directly see \eqref{eqn:longest-superD-id}.
This completes the proof. \hfill $\Box$

\subsection{A longest elements of a Weyl groupoid of the simple Lie superalgebra ${\mathrm{C}}(\rkN)$}
\label{subsection:LGEsuperC}
In this subsection, assume $\rkN\geq 3$.
Define $\ospepCN\in\ospemfD_{1|\rkN-1}$ by 
$\ospepCN(i):=1-\delta_{i,1}-\delta_{i,\rkN+1}$
($i\in\ospefkI$).

Note that if $p=\ospepCN$
and $\superordCbase=(\ospeal^p_i|i\in\fkI)$, then
$(-\ospeeta^p,\superordDbase)$
is the ${\mathrm{C}}(\rkN)$-data and 
$\ospeR(p)$ is the standard root system associated with 
$(-\ospeeta^p,\superordCbase)$.

Define $\ospefCN\in\funcI_{\rkN^2-1}$ by.
\begin{equation*}
\ospefCN(t):=
\left\{\begin{array}{ll}
t+1 & \quad\mbox{if $t\in\fkJ_{1,\rkN-1}$}, \\
\ospefCN(t-(\rkN-1)) & \quad\mbox{if $t\in\fkJ_{\rkN,(\rkN-1)^2}$}, \\
t-(\rkN-1)^2 & \quad\mbox{if $t\in\fkJ_{(\rkN-1)^2+1,(\rkN-1)^2+\rkN}$}, \\
\rkN-1-(t-((\rkN-1)^2+\rkN)) & \quad\mbox{if $t\in\fkJ_{(\rkN-1)^2+\rkN+1,\rkN^2-1}$}.
\end{array}\right.
\end{equation*}

\begin{proposition}\label{proposition:ospCNlgnst}
We have
\begin{equation} \label{eqn:ospCN-lngst}
\tilell_\ospepCN(1^\ospepCN w_0)=\rkN^2-1,\quad 1^\ospepCN w_0=1^\ospepCN\tils_{\ospefCN,\rkN^2-1}
\end{equation} and
\begin{equation}\label{eqn:ospCN-lngstd}
(\ospepCN)_{\ospefCN,\rkN^2-1}(i)=1-\delta_{i,1}\quad(i\in\ospefkI).
\end{equation} Moreover letting 
$m_\al:=-\ospeeta^\ospep((\ospexi^\ospepCN)^{-1}(\al),(\ospexi^\ospepCN)^{-1}(\al))$
$(\al\in\tRp(\ospepCN))$, we have
\begin{equation}\label{eqn:ospCN-lngstdd}
\{\,1^\ospepCN\tils_{\ospefCN,t-1}(\al_{\ospefCN(t)})\,|\,t\in\fkJ_{1,(\rkN-1)^2}\,\}=
\{\al\in\tRp(\ospepCN)|m_\al\in\{2,4\}\},
\end{equation} and 
\begin{equation}\label{eqn:ospCN-lngstddd}
\{\,1^\ospepCN\tils_{\ospefCN,t-1}(\al_{\ospefCN(t)})\,|\,t\in\fkJ_{(\rkN-1)^2+1,\rkN^2-1}\,\}=
\{\al\in\tRp(\ospepCN)|m_\al=0\}.
\end{equation} 
\end{proposition}
{\it{Proof}}. Let $p:=\ospepCN$, $f:=\ospefCN$ and $r:=(\rkN-1)^2$.
For $t\in\fkJ_{1,r}$, we have $p_{f,t}=p$, so \eqref{eqn:defospeal} implies
\begin{equation}\label{eqn:ospCN-pf-0}
\ospeal^{p_{f,t}}_{f(t)}=
\left\{\begin{array}{ll} 
\vece_{f(t)}-\vece_{f(t)+1} 
&\quad\mbox{if $f(t)\in\fkJ_{2,\rkN-1}$}, \\
2\vece_\rkN &\quad\mbox{if $f(t)=\rkN$}. 
\end{array}\right.
\end{equation} In particular,
\begin{equation}\label{eqn:ospCN-pf-1}
p_{f,r}=p.
\end{equation} 

By \eqref{eqn:EXShatb},
\eqref{eqn:evenosp-1-2} and \eqref{eqn:ospCN-pf-0}, we have 
\begin{equation}\label{eqn:EVENospCN-pf-2}
(\ospexi^p)^{-1}\circ 1^p\tils_{f,r}\circ\ospexi^p=
\stProj_{\fkJ_{1,1}}-\stProj_{\fkJ_{2,\rkN}}.
\end{equation}

By \eqref{eqn:ospCN-pf-1}, we directly have \eqref{eqn:ospCN-lngstd} and
\begin{equation}\label{eqn:EVENospCN-pf-3}
\ospeal^{p_{f,r+t}}_{f(r+t)}=
\left\{\begin{array}{ll} 
\vece_t-\vece_{t+1} &\quad\mbox{if $t\in\fkJ_{1,\rkN-1}$}, \\
\vece_{\rkN-1}+\vece_\rkN &\quad\mbox{if $t=\rkN$}, \\
\vece_{2\rkN-1-t}-\vece_{2\rkN-t} &\quad\mbox{if $t\in\fkJ_{\rkN+1,2\rkN-2}$}. 
\end{array}\right.
\end{equation} for $t\in\fkJ_{1,2\rkN-2}$.
Since $\hats_{\vece_{\rkN-1}-\vece_\rkN}\hats_{\vece_{\rkN-1}+\vece_\rkN}=\hats_{\vece_{\rkN-1}}\hats_{\vece_{\rkN}}$,
by \eqref{eqn:EVENospCN-pf-3},
we have 
\begin{equation}\label{eqn:EVENospCN-pf-4}
\underbrace{\hats_{f(r+1)}\cdots\hats_{f(\rkN^2-1)}}_{\rkN^2-1-r}=\hats_{\vece_1}\hats_{\vece_\rkN}
=-\stProj_{\fkJ_{1,1}}+\stProj_{\fkJ_{2,\rkN-1}}-\stProj_{\fkJ_{\rkN,\rkN}}.
\end{equation}
By \eqref{eqn:EXShatb}, \eqref{eqn:EVENospCN-pf-2} and \eqref{eqn:EVENospCN-pf-4}, we have 
\begin{equation} \label{eqn:EVENospCN-pf-5}
(\ospexi^p)^{-1}\circ 1^p\tils_{f,\rkN^2-1}\circ\ospexi^p=
-\stProj_{\fkJ_{1,\rkN-1}}+\stProj_{\fkJ_{\rkN,\rkN}}.
\end{equation} By \eqref{eqn:defospeal}
and \eqref{eqn:EVENospCN-pf-5}, we have
\begin{equation*}
1^p\tils_{f,\rkN^2-1}(\al_i)=
\left\{
\begin{array}{ll}
-\al_i & \quad\mbox{if $i\in\fkJ_{1,\rkN-2}$}, \\
-\al_\rkN & \quad\mbox{if $i=\rkN-1$}, \\
-\al_{\rkN-1} & \quad\mbox{if $i=\rkN$}.
\end{array}\right.
\end{equation*}
Hence
$1^p\tils_{f,\rkN^2-1}
(\Pi)=-\Pi$.
Hence by \eqref{eqn:lgsta} 
and \eqref{eqn:evenosp-1-3},
we have \eqref{eqn:ospCN-lngst}.
By 
\eqref{eqn:ellwf}, \eqref{eqn:defoseR}, 
\eqref{eqn:defospeal}, \eqref{eqn:evenosp-1-1}
and \eqref{eqn:EVENospCN-pf-2},
we see that LHS of \eqref{eqn:ospCN-lngstdd}
is $\tRp(\ospepCN)\cap\ospexi^\ospepCN(\oplus_{x=2}^\rkN\bR\vece_x)$.
Hence we have \eqref{eqn:ospCN-lngstdd}.
This completes the proof.
\hfill $\Box$

\section{{Generalized quantum groups}} \label{section:GQG}
\subsection{Bi-homomorphism $\bhm$ and Dynkin diagram of $\bhm$}\label{subsection:BhDd}
We say that a map $\bhm:\bZPi\times \bZPi\to \bKt$
is a {\it bi-homomorphism on $\Pi$}
if
\begin{equation*}
\bhm(\al,\beta+\gamma)=\bhm(\al,\beta)\bhm(\al,\gamma),\quad
\bhm(\al+\beta,\gamma)=
\bhm(\al,\gamma)\bhm(\beta,\gamma)
\end{equation*} hold
for all $\al$, $\beta$, $\gamma\in\bZPi$.
Let $\bigchiPi$ be the set of bi-homomorphisms on $\Pi$.

Let $\bhm\in\bigchiPi$ and let $q_{ij}:=\bhm(\al_i,\al_j)$
for $i$, $j\in \fkI$.
By {\it{the Dynkin diagram}} of $\bhm$,
we mean the un-oriented graph with
$\rkN$-dots such that each $i$-th dot is labeled
$\al_i$ and $q_{ii}$,
each two $j$-th and $k$-th dots
with $j\ne k$ and  $q_{jk}q_{kj}\ne 1$
are joined by a single line
labeled $q_{jk}q_{kj}$.
For example, if $\rkN=3$ and
$q_{11}=-1$, $q_{22}=\bq^2$, $q_{33}=\bq^6$,
$q_{12}q_{21}=\bq^{-2}$,
$q_{23}q_{32}=\bq^{-6}$ and
$q_{13}q_{31}=1$ for some $\bq\in \bKtinf$,
then the Dynkin diagram of $\bhm$ is given by the leftmost one of Figure~\ref{fig:DynkinG}.
Note that the Dynkin diagram of $\bhm$
does not recover $\bhm$. In fact,
it recovers an equivalent class $\equiv$ of $\bhm$, which will be introduced in
\eqref{eqn:eqvchi} below.

\subsection{Quantum group $U=U(\bhm)$ associated with $\bhm\in\bigchiPi$}
\label{subsection:QgUaw}

From now on until the end of Subsection~\ref{subsection:Irrhigh},
we fix $\bhm\in\bigchiPi$, let $q_{ij}:=\bhm(\al_i,\al_j)$
for $i$, $j\in \fkI$.

Let $\tU:=\tU(\bhm)$ be the unital associative $\bK$-algebra defined
by generators
\begin{equation*}
\tK_\al,\,\tL_\al\, (\al\in\bZPi),\quad
\tE_i, \tF_i\, (i\in \fkI)
\end{equation*} and
relations
\begin{equation}\label{eqn:relone}
\begin{array}{l}
\tK_0=\tL_0=1,\,
\tK_\al \tK_\beta=\tK_{\al+\beta},\,
\tL_\al \tL_\beta=\tL_{\al+\beta},\,
\tK_\al \tL_\beta=\tL_\beta \tK_\al, \\
\tK_\al \tE_i =\bhm(\al,\al_i)\tE_i \tK_\al,\,
\tL_\al \tE_i  =\bhm(-\al_i,\al)\tE_i \tL_\al,\\
\tK_\al \tF_i  =\bhm(\al,-\al_i)\tF_i \tK_\al,\,
\tL_\al \tF_i  =\bhm(\al_i,\al)\tF_i\tL_\al,\\
\tE_i\tF_j-\tF_j\tE_i=\delta_{ij}(-\tK_{\al_i}+\tL_{\al_i})
\end{array}
\end{equation}
for all
$\al$, $\beta\in\bZPi$ and all $i$, $j\in \fkI$.

Define the $\bK$-algebra automorphism
$\tOm: \tU\to\tU$
by $\tOm(\tK_\al):=\tK_{-\al}$,
$\tOm(\tL_\al):=\tL_{-\al}$,
$\tOm(\tE_i):=\tF_i\tL_{-\al_i}$, and
$\tOm(\tF_i):=\tK_{-\al_i}\tE_i$.
Define $\bhmop\in\bigchiPi$ by
$\bhmop(\al,\beta):=\bhm(\beta,\al)$
($\al$, $\beta\in\bZPi$).
Define the $\bK$-algebra automorphism
$\tUp: \tU(\bhm^{\rm{op}})\to\tU(\bhm)$
by $\tUp(\tK_\al):=\tL_\al$, $\tUp(\tL_\al):=\tK_\al$,
$\tUp(\tE_i):=\tF_i$, and
$\tUp(\tF_i):=\tE_i$.

Let $\tU^0:=\tU^0(\bhm)$ (resp. $\tU^+:=\tU^+(\bhm)$, resp. $\tU^-:=\tU^-(\bhm)$)
be the unital subalgebra of
$\tU$ generated by $\tK_\al$,  $\tL_\al$ ($\al\in\bZPi$)
(resp. $\tE_i$ ($i\in \fkI$), resp. $\tF_i$ ($i\in \fkI$)).

\begin{lemma}\label{lemma:tUtr}
The elements
\begin{equation}\label{eqn:tUtreq}
\tF_{j_1}\cdots \tF_{j_m}\tK_\al \tL_\beta\tE_{i_1}\cdots \tE_{i_r}
\end{equation} with
$m$, $r\in\bZgeqo$, $i_x\in\fkI$
$(x\in\fkJ_{1,r})$, $j_y\in\fkI$
$(y\in\fkJ_{1,m})$, $\al$, $\beta\in\bZPi$ form
a $\bK$-basis of $\tU$,
where we use the convention that
if $r=0$ {\rm{(}}resp. $m=0${\rm{)}}, then $\tE_{i_1}\cdots \tE_{i_r}$
{\rm{(}}resp. $\tF_{j_1}\cdots \tF_{j_m}${\rm{)}}
means $1$.
\end{lemma}
{\it{Proof.}} This can be proved in a standard way as in \cite[Lemma~2.2]{Y89}.
\hfill $\Box$
\newline\par
Define the $\bZPi$-grading structure
$\tU=\oplus_{\al\in\bZPi}\tU_\al$
on $\tU$
by $\tK_\al\in \tU_0$, $\tL_\al\in \tU_0$, $\tE_i\in \tU_{\al_i}$,
$\tF_i\in \tU_{-\al_i}$, and $\tU_\al \tU_\beta\subseteq \tU_{\al+\beta}$.
For $\al\in\bZPi$, let
$\tU^\pm_\al:=\tU^\pm\cap\tU_\al$.
Then
$\tU^\pm=\oplus_{\al\in\pm\bZgeqoPi}\tU^\pm_\al$.

For $m\in\bZgeqo$, and $t_1$, $t_2\in\bKt$, let
\begin{equation*}
(m;t_1,t_2):=1-t_1^{m-1}t_2,
\quad\mbox{and}\quad
(m;t_1,t_2)!:=\prod_{j\in \fkJ_{1,m}}(j;t_1,t_2).
\end{equation*}
For $m\in\bZgeqo$ and $i$, $j\in \fkI$ with $i\ne j$, define
$\tEp_{m,\al_i,\al_j}\in\tU^+_{m\al_i+\al_j}$,
and
$\tFp_{m,\al_i,\al_j}\in\tU^-_{-m\al_i-\al_j}$
inductively by
$\tEp_{0,\al_i,\al_j}:=\tE_j$,
$\tFp_{0,\al_i,\al_j}:=\tF_j$, and
\begin{equation}\label{eqn:EFp}
\begin{array}{l}
\tEp_{m+1,\al_i,\al_j}:=\tE_i\tEp_{m,\al_i,\al_j}
-q_{ii}^mq_{ij}\tEp_{m,\al_i,\al_j}\tE_i,
\\
\tFp_{m+1,\al_i,\al_j}:=\tF_i\tFp_{m,\al_i,\al_j}
-q_{ii}^mq_{ji}\tFp_{m,\al_i,\al_j}\tF_i.
\end{array}
\end{equation}
\begin{remark}
The elements in \eqref{eqn:EFp} appear naturally when considering the
twisting of the Hopf algebra structure of quantum groups.
For example, see \cite[Section~7]{Y99}.
Here the term {\it{twisting}} means the method given by \cite[Proposition~7.2.3]{Y99}
in order to obtain another Hopf algebra structure.
\end{remark}
We have
\begin{equation*}
\tUp(\tFp_{m,\al_i,\al_j})=\tEp_{m,\al_i,\al_j}\quad\mbox{and}\quad\tUp(\tEp_{m,\al_i,\al_j})=\tFp_{m,\al_i,\al_j}.
\end{equation*}
\begin{lemma}\label{lemma:tUbidd}
{\rm{(1)}} For $m\in \bN$ and $i$, $j\in\fkI$ with
$i\ne j$, we have
\begin{equation}\label{eqn:tUbiddone}
[\tE_i,\tF_i^m]=(m)_{q_{ii}}(-\tK_{\al_i}+q_{ii}^{-m+1}\tL_{\al_i})\tF_i^{m-1}.
\end{equation}

{\rm{(2)}} For $m\in \bN$ and $i$, $j\in\fkI$ with
$i\ne j$, we have
\begin{equation}\label{eqn:tUbiddtwo}
[\tE_i,\tFp_{m,\al_i,\al_j}]=-(m)_{q_{ii}}
(m;q_{ii},q_{ij}q_{ji})\tK_{\al_i}\tFp_{m-1,\al_i,\al_j}.
\end{equation}

{\rm{(3)}} Let $n$, $m\in \bZgeqo$ with
$n<m$ and $i$, $j\in\fkI$ with
$i\ne j$. Then we have
\begin{equation*}
[\tEp_{n,\al_i,\al_j},\tFp_{m,\al_i,\al_j}]=(n)_{q_{ii}}!
{m \choose n}_{q_{ii}}
(m;q_{ii},q_{ij}q_{ji})!\tF_i^{m-n}\tL_{n\al_i+\al_j}.
\end{equation*} In particular, we have
\begin{equation}\label{eqn:tUbiddfour}
[\tEp_j,\tFp_{m,\al_i,\al_j}]=
(m;q_{ii},q_{ij}q_{ji})!\tF_i^m\tL_{\al_j}.
\end{equation}

{\rm{(4)}} For $m\in \bZgeqo$
and $i$, $j\in\fkI$ with
$i\ne j$, we have
\begin{equation}\label{eqn:tUbiddfive}
[\tEp_{m,\al_i,\al_j},\tFp_{m,\al_i,\al_j}]=(m)_{q_{ii}}!
(m;q_{ii},q_{ij}q_{ji})!(-\tK_{m\al_i+\al_j}+\tL_{m\al_i+\al_j}).
\end{equation}

{\rm{(5)}} For $m$, $n\in \bZgeqo$
and $i$, $j$, $k\in\fkI$ with
$i\ne j\ne k\ne i$, we have
\begin{equation*}
[\tEp_{m,\al_i,\al_j},\tFp_{n,\al_i,\al_k}]=0.
\end{equation*}

\end{lemma}
{\it{Proof.}} These equations are obtained
in a direct way as in \cite[Corollary~4.25, Lemma~4.27]{Hec10}.
\hfill $\Box$
\newline\par
For $\al=\sum_{i\in\fkI}n_i\al_i\in\bZgeqoPi$ with $n_i\in\bZgeqo$,
using induction on $\sum_{i\in\fkI}n_i$,
we define
a $\bK$-subspace $\tcI^-_{-\al}$ of
$\tU^-_{-\al}$ as follows.
Let $\tcI^-_0:=\{0\}$.
For $\al\in\bZgeqoPi\setminus\{0\}$,
let $\tcI^-_{-\al}$ be
the $\bK$-subspace of $\tU^-_{-\al}$
formed by the elements $\tY\in\tU^-_{-\al}$
with $[\tE_i,\tY]\in
\tcI^-_{-\al+\al_i}\tK_{\al_i}+\tcI^-_{-\al+\al_i}\tL_{\al_i}$
for all $i\in\fkI$
with $\al-\al_i\in\bZgeqoPi$.
Note that $\tcI^-_{-\al_i}=\{0\}$ for all $i\in\fkI$.

Let $\tcI^-:=\oplus_{\al\in\bZgeqoPi}\tcI^-_{-\al}$,
and set $\tcJ^-:=\rmSpan_\bK (\tcI^-\tU^0\tU^+)$. Then $\tcJ^-$ is an ideal of $\tU$.
We define the unital $\bK$-algebra $U=U(\bhm)$ by
\begin{equation}\label{eqn:dfU}
U:=\tU/(\tcJ^- +\tOm(\tcJ^-))=\tU/(\tcJ^- +\tUp(\tcJ^-))
\end{equation} (the quotient algebra).

The $\bK$-algebra $U(\bhm)$ defined by \eqref{eqn:dfU} is isomorphic to the one given by \cite[(3.14)]{HY10}, which is defined in a way similar to that given by Lusztig~\cite[3.1.1 (a)-(e)]{b-Lusztig93}.

\begin{remark}\label{remark:rYamane}
As for the defining relations of $U(\bhm)$
for $\bhm\in\dotXNSuper$ of Definition~\ref{definition:defSuperbhm-d} below, see \cite{A11}, \cite{AAY10},
\cite{Y94}, \cite{Y99}, \cite{Y01}, \cite{YHome}.
\end{remark}

Let $\pi_\bhm:=\pi:\tU\to U$ be the canonical map.
We denote $\pi(\tK_\al)$, $\pi(\tL_\al)$,
$\pi(\tE_i)$, $\pi(\tF_i)$, $\pi(\tEp_{m,\al_i,\al_j})$, $\pi(\tFp_{m,\al_i,\al_j})$
by $K_\al$, $L_\al$, $E_i$, $F_i$, $E_{m,\al_i,\al_j}$, and
$F_{m,\al_i,\al_j}$
respectively. Let $U^0:=U^0(\bhm):=\pi(\tU^0)$, and $U^\pm:=U^\pm(\bhm):=\pi(\tU^\pm)$.
For $\al\in\bZPi$, let $U_\al:=U(\bhm)_\al:=\pi(\tU_\al)$,
and $U^\pm_\al:=U^\pm(\bhm)_\al:=\pi(\tU^\pm_\al)$.

We have the $\bK$-algebra automorphism
$\Omega: U\to U$ with $\Omega\circ\pi=\pi\circ\tOm$.
We have the $\bK$-algebra isomorphism
$\Upsilon: U(\bhmop)\to U(\bhm)$ with $\Upsilon\circ\pi_\bhmop=\pi_\bhm\circ\tUp$.

We can easily see
\begin{lemma}\label{lemma:Utr}
There exists a unique $\bK$-linear isomorphism
from
$U^-\otimes U^0 \otimes U^+$
to $U$
sending
$Y\otimes Z\otimes X$
to
$YZX$
$(X\in U^+,\,Z\in U^0,\,Y\in U^-)$.
Moreover, we have
$U^0=\oplus_{\al,\beta\in\bZPi}\bK K_\al L_\beta$,
$U=\oplus_{\al\in\bZPi}U_\al$,
$U^\pm=\oplus_{\al\in\bZgeqoPi}U^\pm_{\pm\al}$, and
$\dim U^+_\al=\dim U^-_{-\al}$ for all $\al\in\bZgeqoPi$.
\end{lemma}

\begin{remark}\label{remark:idUzero}
By Lemma~\ref{lemma:tUtr}, we have the same results as Lemma~\ref{lemma:Utr}
with $\tU$ in place of $U$.
By Lemma~\ref{lemma:Utr}, we see that
the structure of $U^0$ as an unital $\bK$-algebra is independent of the choice of $\bhm\in\bigchiPi$.
\end{remark}

\subsection{Kharchenko-PBW theorem}
\label{subsection:KharPBW}

Define the map $h^\bhm:\bZPi\to\bN\cup\{\infty\}$ by
\begin{equation*}
h^\bhm(\al):=\left\{
\begin{array}{ll}
\infty & \mbox{if $(m)_{\bhm(\al,\al)}!\ne 0$ for all $m\in\bN$,} \\
\rmMax\{m\in\bN\,|\,(m)_{\bhm(\al,\al)}!\ne 0\} & \mbox{otherwise.}
\end{array}\right.
\end{equation*} For $i$, $j\in \fkI$, define
$c_{ij}:=c^\bhm_{ij}\in\{2\}\cup \fkJ_{-\infty,0}\cup\{-\infty\}$ by
\begin{equation}\label{eqn:cxij}
c_{ij}:=c^\bhm_{ij}:=\left\{
\begin{array}{l}
2 \,\,\mbox{if $i=j$,} \\
-\infty \,\,\mbox{if $i\ne j$ and $(m)_{q_{ii}}!(m;q_{ii},q_{ij}q_{ji})!\ne 0$
for all $m\in\bN$,} \\
-\rmMax\{m\in\bZgeqo\,|\,(m)_{q_{ii}}!(m;q_{ii},q_{ij}q_{ji})! \ne 0\}\,\,
\mbox{otherwise.}
\end{array}\right.
\end{equation}

By \eqref{eqn:tUbiddone} and
Lemma~\ref{lemma:Utr}, the following lemma holds.
\begin{lemma} \label{lemma:easyone}
Let $i\in\fkI$ and $m\in\bZgeqo$. Then
$F_i^m=0$ if and only if $m>h^\bhm(\al_i)$.
In particular, if $m\in \fkJ_{0,h^\bhm(\al_i)}$,
then $\dim U^-_{-m\al_i}=1$, and if $m>h^\bhm(\al_i)$, then
$\dim U^-_{-m\al_i}=0$.
\end{lemma}

We also have the following result.
\begin{lemma} \label{lemma:easytwo} Let $i$, $j\in\fkI$ with $i\ne j$, and $m\in\bZgeqo$.

{\rm{(1)}} $F_{m,\al_i,\al_j}\ne 0$ if and only if
$m\in\fkJ_{0,-c_{ij}}$.

{\rm{(2)}} If $m>h^\bhm(\al_i)-c_{ij}$, $\dim U^-_{-m\al_i-\al_j}=0$.
Also if $m\leq h^\bhm(\al_i)-c_{ij}$,
the elements $F_{r,\al_i,\al_j}F_i^{m-r}$
with $r\in\fkJ_{\rmMax\{0,m-h^\bhm(\al_i)\},{\rm{Min}}\{m,-c_{ij}\}}$ form
a $\bK$-basis of $U^-_{-m\al_i-\al_j}$.

{\rm{(3)}} $\dim U^-_{-\al_i-\al_j}\in\fkJ_{1,2}$.
Moreover,
$q_{ij}q_{ji}=1$ $\Leftrightarrow$ $c_{ij}=0$ $\Leftrightarrow$
$\dim U^-_{-\al_i-\al_j}= 1$ $\Leftrightarrow$ $F_iF_j=q_{ij}F_jF_i$.
\end{lemma}
{\it{Proof.}} {\rm{(1)}} This follows from
Lemmas~\ref{lemma:Utr} and \ref{lemma:easyone}, together with
\eqref{eqn:tUbiddtwo}, \eqref{eqn:tUbiddfour}, and
\eqref{eqn:tUbiddfive}.

{\rm{(2)}} If $m=0$, this follows from Lemma~\ref{lemma:easyone}.
Assume $m\geq 1$. Note that $-c_{ij}\leq h^\bhm(\al_i)$.
By (1) and \eqref{eqn:EFp}, we see that $U^-_{-m\al_i-\al_j}$
is spanned by elements as in the statement.

Assume $m\leq h^\bhm(\al_i)-c_{ij}$.
Let $Z:=\fkJ_{\rmMax\{0,m-h^\bhm(\al_i)\},{\rm{Min}}\{m,-c_{ij}\}}$,
and $X:=\sum_{r\in Z}y_rF_{r,\al_i,\al_j}F_i^{m-r}$
with $y_r\in\bK$.
Assume $X=0$.
Observing coefficients of $F_{r,\al_i,\al_j}F_i^{m-r-1}L_{\al_i}$
($r\in Z\cap\fkJ_{0,m-1}$)
of $[E_i,X]$, by Lemma~\ref{lemma:Utr}, and \eqref{eqn:tUbiddone}, \eqref{eqn:tUbiddtwo} and induction,
we see $y_r=0$ for $r\in Z\cap\fkJ_{0,m-1}$.
In a similar way, we see $y_m=0$ (if $m\leq -c_{ij}$) by \eqref{eqn:tUbiddtwo}.

{\rm{(3)}} This follows from (1) and (2). \hfill $\Box$
\newline\par
By the celebrated Kharchenko's PBW theorem, we have

\begin{theorem}\label{theorem:Khth} {\rm{(}}{\rm{\cite{Kha99}}}, see also 
{\rm{\cite[{\it{Section}}~3, (P)]{Hec06}}},
and {\rm{\cite[(2.15)]{HY10}{\rm{)}}}}

{\rm{(1)}} There exists a unique pair
\begin{equation*}
(R^+=R^+(\bhm),\varphi=\varphi_\bhm),
\end{equation*} where $R^+$ is a subset of $\bZgeqoPi\setminus\{0\}$, and
$\varphi$ is a map $\varphi:R^+\to\bN$,
satisfying the condition that there exist $k\in\bN\cup\{\infty\}$,
a surjective map $\psi:\fkJ_{1,k}\to R^+$
with $|\psi^{-1}(\{\al\})|=\varphi (\al)$
$(\al\in R^+)$, and $F[r]\in U^-_{-\psi(r)}\setminus\{0\}$
$(r\in \fkJ_{1,k})$ such that the elements
\begin{equation}\label{eqn:KhPBWbs}
F[1]^{x_1}\cdots F[m]^{x_m}\quad(m\in \fkJ_{1,k},\,
x_y\in \fkJ_{0,h^\bhm(\psi(y))}\,(y\in \fkJ_{1,m}))
\end{equation} form a $\bK$-basis of $U^-$
{\rm{(}}where we mean that for $m\leq {\bar m}$,
$F[1]^{x_1}\cdots F[m]^{x_m}=F[1]^{{\bar x}_1}\cdots F[{\bar m}]^{{\bar x}_{\bar m}}$
if and only if
$x_y={\bar x}_y$ for all $y\in \fkJ_{1,m}$ and ${\bar x}_{\bar y}=0$
for all ${\bar y}\in \fkJ_{m+1,{\bar m}}${\rm{)}}.

{\rm{(2)}} 
If $|R^+|<\infty$, then $\varphi(R^+)=\{1\}$.
{\rm{(}}See also Remark~{\rm{\ref{remark:prprtsys}~(1)}} below.{\rm{)}}
\end{theorem}

\begin{remark}\label{remark:prprtsys}
(1) 
Using $\Omega$ and $\Upsilon$,
we see that $(R^+(\bhmop),\varphi_\bhmop)=(R^+(\bhm),\varphi_\bhm)$.
It is clear that $\varphi_\bhm(\al_i)=1$ for $i\in\fkI$.
Then Theorem~\ref{theorem:Khth}~{(2)} can also be proved by using
Theorem~\ref{theorem:hecprone}~{(1)} and \eqref{eqn:ellRp} below.

(2) The following facts are well-known,
see \cite[33.3, Corollary~33.1.5]{b-Lusztig93} and \cite[Section~2.4]{HK07} for example.
Let $\bq\in\bKtinf$.
Let $A=[a_{ij}]_{i,j\in\fkI}$ be a symmetrizable generalized Cartan matrix.
Let $d_i\in\bN$ ($i\in\fkI$) be such that $d_ia_{ij}=d_ja_{ji}$
($i$, $j\in\fkI$). Let $\mfkg$ be the Kac-Moody Lie algebra defined for $A$.
Let ${\mathfrak{n}}^-$ be the negative part of $\mfkg$,
and ${\mathcal{U}}({\mathfrak{n}}^-)$ be the universal enveloping algebra of ${\mathfrak{n}}^-$.
Let $\bhm\in\bigchiPi$ be such that
$\bhm(\al_i,\al_j)=\bq^{d_ia_{ij}}$ ($i$, $j\in\fkI$).
Then the ideal $\tcI^-$ of $\tU^-(\bhm)$ is generated by
$\tF_{1-a_{ij},\al_i,\al_j}$ ($i\ne j$),
and $\dim U^-(\bhm)_{-\al}=\dim{\mathcal{U}}({\mathfrak{n}}^-)_{-\al}$
for all $\al\in\bZgeqoPi$, where ${\mathcal{U}}({\mathfrak{n}}^-)_{-\al}$ is the weight subspace of
${\mathcal{U}}({\mathfrak{n}}^-)$ corresponding to $-\al$.
In particular, $R^+(\bhm)$ can be identified with the set of positive roots of ${\frak{g}}$
and $\varphi_\bhm(\al)=\dim\mathfrak{g}_\al(=\dim\mathfrak{g}_{-\al})$
for all $\al\in R^+(\bhm)$.
\end{remark}

Once we know Theorem~\ref{theorem:Khth},
the following lemma is clear from Lemmas~\ref{lemma:easyone} and \ref{lemma:easytwo}.

\begin{lemma}\label{lemma:elprofR}
{\rm{(1)}} $\Pi\subseteq R^+$.

{\rm{(2)}} For $i$, $j\in \fkI$ with $i\ne j$,
we have
\begin{equation}\label{eqn:elprofRd}
R^+\cap (\al_j+\bZgeqo\al_i)=\{\al_j+n\al_i\,|\,n\in\fkJ_{0,-c^\bhm_{ij}}\},
\end{equation}
and $\varphi(\al_j+r\al_i)=1$ for all $r\in \fkJ_{0,-c_{ij}}$.

{\rm{(3)}} Let $\fkI^\prime$
be a non-empty proper subset of $\fkI$.
Let $\fkI^{\prime\prime}=\fkI\setminus\fkI^\prime$.
Then $R^+=(R^+\cap\oplus_{i\in \fkI^\prime}\bZgeqo\al_i)
\uplus (R^+\cap\oplus_{j\in \fkI^{\prime\prime}}\bZgeqo\al_j)$
if and only if $q_{ij}q_{ji}=1$ for all $i\in \fkI^\prime$
and all $j\in \fkI^{\prime\prime}$.

\end{lemma}

For $\bhm\in\bigchiPi$, we say that $\bhm$ is {\it{irreducible}}
if
its Dynkin diagram is connected, that is,
for any two $i$, $j\in \fkI$ with $i\ne j$, there exist
$m\in\bZgeqo$, $i_r\in \fkI$ ($r\in \fkJ_{1,m}$)
such that $q_{i_t,i_{t+1}}q_{i_{t+1},i_t}\ne 1$
for all $t\in \fkJ_{0,m}$, where we let $i_0:=i$
and $i_{m+1}:=j$.
Let $\bigchiPiirr :=\{\bhm\in \bigchiPi\,|\,\mbox{ $\bhm$ is irreducible}\}$.

\subsection{Irreducible highest weight modules}
\label{subsection:Irrhigh}
Let $\Lambda\in \rmCh(U^0)$.
By Lemma~\ref{lemma:Utr}, we have a unique left $U$-module
$\mclM_\bhm(\Lambda)$ satisfying the following conditions.
\newline\par
(i) There exists $\tv_\lambda\in \mclM_\bhm(\Lambda)\setminus\{0\}$
such that $Z\tv_\Lambda=\Lambda(Z)\tv_\Lambda$ for all $Z\in U^0$
and $E_i\tv_\Lambda=0$ for all $i\in \fkI$.
\par
(ii) The $\bK$-linear map $U^-\to\mclM_\bhm(\Lambda)$, $Y\mapsto Y\tv_\Lambda$,
is bijective.
\newline\par For $i\in \fkI$ and $m\in\bZgeqo$, by \eqref{eqn:tUbiddone}, we have
\begin{equation}\label{eqn:esyf}
E_iF_i^m\tv_\Lambda=
\left\{\begin{array}{ll}
0 &\quad\mbox{if $m=0$,} \\
(m)_{q_{ii}}
(-q_{ii}^{1-m}\Lambda(K_{\al_i})+\Lambda(L_{\al_i}))
F_i^{m-1}\tv_\Lambda & \quad\mbox{otherwise.}
\end{array}\right.
\end{equation}
For $\al\in\bZPi$,
let $\mclM_\bhm(\Lambda)_\al:=U^-_\al \tv_\Lambda$.
We say that a $\bK$-subspace $\mathcal{V}$ of $\mclM_\bhm(\Lambda)$
is {\it{$\bZPi$-graded}} if
$\mathcal{V}=\oplus_{\al\in\bZgeqoPi}(\mathcal{V}\cap \mclM_\bhm(\Lambda)_\al)$.
If $\mathcal{V}$ is a $\bZPi$-graded $U(\bhm)$-submodule  of
$\mclM_\bhm(\Lambda)$, then $\mathcal{V}\ne\mclM_\bhm(\Lambda)$
if and only if $\tv_\Lambda\notin\mathcal{V}$.

Let $\cN:=\cN_\bhm(\Lambda)$ be the maximal
proper $\bZPi$-graded $U(\bhm)$-submodule of
$\mclM_\bhm(\Lambda)$. Note $\cN\cap\bK\tv_\Lambda=\{0\}$.
Let $\mclL_\bhm(\Lambda)$ be
the quotient left $U$-module
defined by
\begin{equation*}
\mclL_\bhm(\Lambda):=\mclM_\bhm(\Lambda)/\cN.
\end{equation*}
We denote the element $\tv_\Lambda+\cN$ of
$\mclL_\bhm(\Lambda)$ by $v_\Lambda$.
For $\al\in\bZPi$, let $\mclL_\bhm(\Lambda)_\al:=U^-_\al v_\Lambda$.
Then $\mclL_\bhm(\Lambda)=\oplus_{\al\in\bZgeqoPi}\mclL_\bhm(\Lambda)_{-\al}$.
We also have $\mclL_\bhm(\Lambda)_0=\bK v_\Lambda$, and $\dim\mclL_\bhm(\Lambda)_0=1$.
We say that a $\bK$-subspace $\mathcal{V}^\prime$ of $\mclL_\bhm(\Lambda)$
is {\it{$\bZPi$-graded}} if
$\mathcal{V}^\prime=\oplus_{\al\in\bZgeqoPi}(\mathcal{V}^\prime\cap\mclL_\bhm(\Lambda)_\al)$.
Then there exists no non-zero proper $\bZPi$-graded $U(\bhm)$-submodule
of $\mclL_\bhm(\Lambda)$.

By \eqref{eqn:esyf}, for $m\in\bN$, we have
\begin{equation}\label{eqn:fvo}
F_i^m v_\Lambda=0\quad\Longleftrightarrow\quad (m)_{q_{ii}}!(m;q_{ii}^{-1},\Lambda(K_{\al_i}L_{-\al_i}))!=0.
\end{equation} Let
\begin{equation}\label{eqn:hixl}
h_{\bhm,\Lambda,i}:=\left\{
\begin{array}{l}
\infty \quad\mbox{if $(m)_{q_{ii}}!(m;q_{ii}^{-1},\Lambda(K_{\al_i}L_{-\al_i}))!\ne 0$
for all $m\in\bN$,} \\
\rmMax\{m\in\bZgeqo|(m)_{q_{ii}}!(m;q_{ii}^{-1},\Lambda(K_{\al_i}L_{-\al_i}))!\ne 0\} \\
\quad\mbox{otherwise.}
\end{array}\right.
\end{equation} By \eqref{eqn:fvo}, since $\mclL_\bhm(\Lambda)$ is $\bZPi$-graded, we have
\begin{equation}\label{eqn:fvod}
\dim\mclL_\bhm(\Lambda)< \infty\quad\Longrightarrow\quad\forall i\in \fkI,\,
h_{\bhm,\Lambda,i}< \infty.
\end{equation}

\subsection{Notation $\equiv$}\label{subsection:NTequiv}

\begin{notation}
Let $\bhm$, $\bhmp\in\bigchiPi$.
Let $q_{ij}:=\bhm(\al_i,\al_j)$ and $q^\prime_{ij}:=\bhmp(\al_i,\al_j)$.
We write
\begin{equation}\label{eqn:eqvchi}
\bhm\equiv \bhmp
\end{equation} if $q_{ii}=q^\prime_{ii}$ for all $i\in\fkI$
and $q_{jk}q_{kj}=q^\prime_{jk}q^\prime_{kj}$
for all $j$, $k\in\fkI$. By \eqref{eqn:cxij},
\begin{equation}\label{eqn:cijccp}
\bhm\equiv \bhmp\quad \Longrightarrow\quad
\forall i,\,\forall j\in \fkI,\,\,c^\bhm_{ij}=c^{\bhmp}_{ij}.
\end{equation}
Let $\Lambda\in \rmCh(U^0(\bhm))$ and $\Lambda^\prime\in \rmCh(U^0(\bhmp))$.
We write
\begin{equation*}
(\bhm,\Lambda)\equiv (\bhmp,\Lambda^\prime)
\end{equation*} if $\bhm\equiv \bhmp$ and
$\Lambda(K_{\al_i}L_{-\al_i})=\Lambda^\prime(K_{\al_i}L_{-\al_i})$
for all $i\in\fkI$.
\end{notation}

Note that
\begin{equation}\label{eqn:hchcp}
(\bhm,\Lambda)\equiv (\bhmp,\Lambda^\prime)
\quad\Longrightarrow\quad\forall i\in\fkI,\,h_{\bhm,\Lambda,i}=h_{\bhmp,\Lambda^\prime,i}.
\end{equation}

\subsection{Weyl groupoids of bi-homomorphisms}
\label{subsection:WGpoid-bih}

Let $\bigchiPifin:=\{\bhm\in\bigchiPi\,|\,|R^+(\bhm)|<\infty\}$.
If $\bhm\in\bigchiPifin\cap\bigchiPiirr$,  we say that $\bhm$ is {\it{of finite-type}}.

Let $i\in\fkI$,
and set $(\ppbigchiPifin)_i:=\{\bhm\in\bigchiPi\,|\,\forall j\in\fkI,\,c^\bhm_{ij}\ne -\infty\,\}$.
For $\bhm\in(\ppbigchiPifin)_i$, define
$\sbhm_i\in\rmGL_\rkN(\bR)$ by $\sbhm_i(\al_j)=\al_j-c^\bhm_{ij}\al_i$.
Note that $\sbhm_i(\al_i)=-\al_i$, and $(\sbhm_i)^2=\rmid_\bV$ ($i\in\fkI$). Define
$\newtr_i\bhm\in\bigchiPi$ 
by
\begin{equation}\label{eqn:newtribhm}
\newtr_i\bhm(\al,\beta)
:=\bhm(\sbhm_i(\al),\sbhm_i(\beta))
\end{equation} for all $\al$, $\beta\in\bZPi$.

\begin{lemma}\label{lemma:ichilm}
{\rm{(}}See also {\rm{\cite[(2.10)-(2.11)]{HY10}}}.{\rm{)}} Let $\bhm\in(\ppbigchiPifin)_i$.

{\rm{(1)}} If $\bhm(\al_i,\al_i)=1$, then
$\bhm(\al_i,\al_j)\bhm(\al_j,\al_i)=1$ for all $j\in\fkI$.

{\rm{(2)}}
We have
\begin{equation}\label{eqn:ichi}
c^{\newtr_i\bhm}_{ij}=c^\bhm_{ij}\,(j\in\fkI),\,\,\newtr_i\bhm\in(\ppbigchiPifin)_i,\,\,
s^{\newtr_i\bhm}_i=s^\bhm_i,\,\, \mbox{and}\,\, \newtr_i \newtr_i\bhm=\bhm.
\end{equation}

{\rm{(3)}}
If $\bhm\in\bigchiPiirr$, then $\newtr_i\bhm\in\bigchiPiirr$.
\end{lemma}
{\it{Proof.}}
(1) This is clear from \eqref{eqn:cxij}
since $\bhm\in(\ppbigchiPifin)_i$.

(2) Let $c_{ij}:=c^\bhm_{ij}$ ($i,j\in\fkI$).
Let $q_{xy}:=\bhm(\al_x,\al_y)$, and $q^\prime_{xy}:=\newtr_i\bhm(\al_x,\al_y)$ ($x$, $y\in\fkI$).
We have $q^\prime_{ii}=q_{ii}$, and
$(q^\prime_{ii})^mq^\prime_{ij}q^\prime_{ji}=(q_{ii}^{2(-c_{ij})-m}(q_{ij}q_{ji}))^{-1}$
for $j\in\fkI$ and $m\in\bZ$.
If $q_{ii}=1$, the statement is clear from (1) since
$c_{ij}=0$ for all $j\in\fkI\setminus\{i\}$.
Assume $q_{ii}\ne 1$. Let $j\in\fkI\setminus\{i\}$.
Assume $c_{ij}=0$. Then $q_{ij}q_{ji}=1$.
Hence $q^\prime_{ij}q^\prime_{ji}=1$, so $c^{\newtr_i\bhm}_{ij}=0$.
Assume $-c_{ij}\geq 1$.
Assume $q_{ii}^{-c_{ij}}q_{ij}q_{ji}=1$.
Then $(q^\prime_{ii})^{-c_{ij}}q^\prime_{ij}q^\prime_{ji}=1$. Moreover, we have
$(q^\prime_{ii})^{m+1}=q_{ii}^{m+1}\ne 1$ and
$(q^\prime_{ii})^mq^\prime_{ij}q^\prime_{ji}=q_{ii}^{-(-c_{ij}-m)}\ne 1$
for $m\in\fkJ_{0,-c_{ij}-1}$. Hence $c^{\newtr_i\bhm}_{ij}=c_{ij}$.
Assume $q_{ii}^{-c_{ij}}q_{ij}q_{ji}\ne 1$. Then
$q_{ii}$
is a primitive $(1-c_{ij})$-th root of unity.
Since $q_{ii}^mq_{ij}q_{ji}\ne 1$ for $m\in\fkJ_{0,-c_{ij}}$,
we have
$q_{ii}^nq_{ij}q_{ji}\ne 1$ for all $n\in\bZ$.
Hence $c^{\newtr_i\bhm}_{ij}=c_{ij}$.

(3) Let $q_{xy}$ and $q^\prime_{xy}$ be as above.
By \eqref{eqn:cxij}, 
for $x$, $y\in\fkI$
with $c_{ix}=c_{iy}=0$, we have $q_{xy}=q^\prime_{xy}$.
Then (3) follows from (2).
\hfill $\Box$
\newline\par
Let $\ppbigchiPifin:=\cap_{i\in\fkI}(\ppbigchiPifin)_i$.
By \eqref{eqn:elprofRd}, we have
\begin{equation}\label{eqn:subpp}
\bigchiPifin\subseteq\ppbigchiPifin.
\end{equation}

\begin{notation}
Let $\bhm,\,\bhmp\in\bigchiPi$.
We write $\bhm\sim\bhmp$ if
$\bhm=\bhmp$ or
there exist $m\in\bN$,
$i_t\in \fkI$  ($t\in \fkJ_{1,m}$) and
$\bhm_r\in\bigchiPi$ ($r\in \fkJ_{1,m+1}$)
such that $\bhm=\bhm_1$, $\bhmp=\bhm_{m+1}$ and
$\bhm_t\in(\ppbigchiPifin)_{i_t}$,
$\newtr_{i_t}\bhm_t=\bhm_{t+1}$ ($t\in \fkJ_{1,m}$).
\end{notation}

For $\bhm\in\bigchiPi$, let
$\Abhm:=\{\,\bhmp\in\bigchiPi\,|\,
\bhmp\sim\bhm\,\}$.

Let $\pbigchiPifin:=\{\,\bhm\in\ppbigchiPifin\,|\,
\newtr_i\bhmp\in\ppbigchiPifin\,(\bhmp\in\Abhm,\,i\in\fkI)\,\}$.

\begin{definition}
{\rm{Let $\bhm\in\pbigchiPifin$. 
For $i\in\fkI$, define the map
$\newtrAbhm_i:\Abhm\to\Abhm$ by
$\newtrAbhm_i(\bhmp):=\newtr_i\bhmp$.
For $\bhmp\in\Abhm$,
let $C^{\bhmp}$
be the $\rkN\times \rkN$-matrix $[c^{\bhmp}_{ij}]$
over $\bZ$, where note that by Lemma~\ref{lemma:easytwo}~(3), $C^{\bhmp}$
is a generalized Cartan matrix,
see $({\rm{M}}1)$, $({\rm{M}}2)$ of Subsection~\ref{eqn:deftilsa}. We call the quadruple
\begin{equation*}
\Cbhm=\Cbhm(\fkI,\Abhm,(\newtrAbhm_i)_{i\in\fkI},(C^{\bhmp})_{\bhmp\in\Abhm})
\end{equation*} {\it{the Cartan scheme associated with}} $\bhm$. 
Indeed, by \eqref{eqn:ichi}, $\Cbhm$ is a Cartan scheme, 
see $({\rm{C}}1)$, $({\rm{C}}2)$ of Subsection~\ref{eqn:deftilsa}.
}}
\end{definition} 
Let $\bhm\in\pbigchiPifin$.
Recall Notation~\ref{notation:Ab-chift-1}. Since $\Cbhm$ is a Cartan scheme, by \eqref{eqn:newtribhm}, we have
\begin{equation} \label{eqn:newtribhmd}
\begin{array}{l}
\bhm_{f,t-1}(\al_{f(t)},\al_{f(t)})=\bhm(\bhms_{f,t-1}(\al_{f(t)}),\bhms_{f,t-1}(\al_{f(t)})) \\
\quad (n\in\bZgeqo\cup\{\infty\},\,f\in\funcI_n,\,t\in\bN).
\end{array}
\end{equation}

For $\bhm\in\bigchiPi$,
let $R(\bhm):=R^+(\bhm)\cup (-R^+(\bhm))$,
and extend the initial domain of $\varphi_\bhm$
to $R(\bhm)$ by $\varphi_\bhm(-\al)=\varphi_\bhm(\al)$.

It is well-known that
\begin{theorem} \label{theorem:hecprone} {\rm{(\cite[Proposition~1]{Hec06},
{\it{see also}} \cite[Example~4]{HY08})}}
We have the followings.

{\rm{(1)}} For $i\in\fkI$ and $\bhm\in(\ppbigchiPifin)_i$, we have
\begin{equation}\label{eqn:hecpronedash} 
\begin{array}{l}
\newtr_i\bhm\in(\ppbigchiPifin)_i, \,\,
s^\bhm_i(R^+(\bhm)\setminus\{\al_i\})=R^+(\newtr_i\bhm)\setminus\{\al_i\}, \\
\varphi_{\newtr_i\bhm}(s^\bhm_i(\beta))=\varphi_\bhm(\beta) \,\,
(\beta\in R(\bhm)).
\end{array}
\end{equation} 
In particular,
$s^\bhm_i(R(\bhm))=R(\newtr_i\bhm)$.

{\rm{(2)}} Let $\bhm\in\bigchiPifin$.
Then $\newtr_i\bhm\in\bigchiPifin$
and $|R^+(\bhm)|=|R^+(\newtr_i\bhm)|$.

\end{theorem}
We have obtained the first property in \eqref{eqn:hecpronedash} by \eqref{eqn:ichi}.

The following theorem is also well-known.

\begin{theorem} {\rm{(\cite[{\it{Theorem}}~3.14]{Hec10})}}
\label{theorem:hecproneAX} 
Let $\bhm\in\pbigchiPifin$. 
Then the data $\Rbhm=\Rbhm(\Cbhm,(R(\bhmp))_{\bhmp\in\Abhm})$ is
{\it{a root system}} of type $\Cbhm$
{\rm{(}}see Definition~{\rm{\ref{definition:GnRtSys}}}{\rm{)}}.
In particular, for $\bhmp$, $\bhmpp\in\Abhm$ 
and $w\in\mcH(\bhmp,\bhmpp)$,
we have
\begin{equation}\label{eqn:prevR}
w(R(\bhmpp))=R(\bhmp).
\end{equation}
\end{theorem}
{\it{Proof.}} This can easily be shown by Theorem~\ref{theorem:hecprone}
and Definition~\ref{definition:GnRtSys}. \hfill $\Box$

\begin{corollary}\label{corollary:coince}
Let $\bhm$, $\bhmp\in\pbigchiPifin$ be such that $R(\bhm)=R(\bhmp)$
{\rm{(}}as a subset of $\bV${\rm{)}}. 
Then $\bhms_{f,t}=\bhmps_{f,t}$ 
{\rm{(}}as an element of $\rmGL_\rkN(\bR)${\rm{)}} for all
$f\in\funcI_\infty$ and 
$t\in\bZgeqo$.
\end{corollary}
{\it{Proof.}} This follows from Lemma~\ref{lemma:elprofR}
and \eqref{eqn:prevR}.
(See also Lemma~\ref{lemma:techlm}.)
\hfill $\Box$
\newline\par
By \eqref{eqn:subpp} and Theorem~\ref{theorem:hecprone}~(2), we have
\begin{equation}\label{eqn:subp}
\bigchiPifin\subseteq\pbigchiPifin.
\end{equation}

\begin{lemma} \label{lemma:crtyp}
Let $i\in \fkI$ and $\bhm\in(\ppbigchiPifin)_i$.
Assume
$\bhm(\al_i,\al_i)\in \bKtinf$. Then
$\newtr_i\bhm\equiv \bhm$.
\end{lemma}
{\it{Proof.}} Let $q_{xy}:=\bhm(\al_x,\al_y)$ and
$q^\prime_{xy}:=\newtr_i\bhm(\al_x,\al_y)$  for $x$, $y\in\fkI$.
Let $j$, $k\in\fkI$.
Assume $j\ne i\ne k$. By \eqref{eqn:cxij}, $q_{ii}^{-c_{ij}}q_{ij}q_{ji}=q_{ii}^{-c_{ik}}q_{ik}q_{ki}=1$
since $q_{ii}\in \bKtinf$.
Hence
\begin{equation*}
\begin{array}{ll}
q^\prime_{jk}q^\prime_{kj}&=\bhm(\al_j-c_{ij}\al_i,\al_k-c_{ik}\al_i)\bhm(\al_k-c_{ik}\al_i,\al_j-c_{ij}\al_i) \\
&= q_{jk}q_{kj}(q_{ij}q_{ji})^{-c_{ik}}(q_{ik}q_{ki})^{-c_{ij}}q_{ii}^{2c_{ij}c_{ik}} \\
&= q_{jk}q_{kj},
\end{array}
\end{equation*} as desired.
The other cases can be treated similarly. \hfill $\Box$
\newline\par
Let $i\in\fkI$.
Let $\bhm\in(\ppbigchiPifin)_i$,
and let $\bhmp\in\bigchiPi$.
By \eqref{eqn:cijccp}, we have
\begin{equation}\label{eqn:cijccpd}
\bhm\equiv\bhmp\,\,\Longrightarrow\,\,
\bhmp\in(\ppbigchiPifin)_i,\,\,s^\bhm_i=s^{\bhmp}_i,
\,\,\newtr_i\bhm\equiv \newtr_i\bhmp.
\end{equation}

\begin{lemma} \label{lemma:prlgelmmm-d}
Let $\bhm\in\bigchiPifin$.
Let $\bhmp\in\bigchiPi$ be such that
$\bhmp\equiv\bhm$. Then we have
$\bhmp\in\bigchiPifin$,
$\lgst =1^{\bhmp}w_0$
and $R(\bhmp)=R(\bhm)$.
\end{lemma}
{\it{Proof.}}
This follows from \eqref{eqn:cijccpd}
and Lemma~\ref{lemma:techlm}.
\hfill $\Box$

\section{FID-type bi-homomorphisms}

\subsection{Some bi-homomorphisms}\label{subsection:Somebihm}
\begin{definition}\label{definition:defCartanbhm}
{\rm{Recall Definition~\ref{definition:strtbsis}.
Let $\dotXNCartan$ be the subset of $\bigchiPi$ formed by 
$\bhm$ with
\begin{equation}\label{eqn:bhmCtn}
\bhm(\al_i,\al_j)=\bq^{\hateta(\hatal_i,\hatal_j)}\quad(i, j\in\fkI)
\end{equation} 
for some $\bq\in \bKtinf$ and some rank-$\rkN$ Cartan data
$\orderrsysbase=(\hatal_i|i\in\fkI)$.
Let $\dotXNCartan({\mathrm{X}}_\rkN)$ be the subset of $\dotXNCartan$ formed by $\bhm$ 
as in \eqref{eqn:bhmCtn} for which
$\orderrsysbase=(\hatal_i|i\in\fkI)$ are 
the ${\mathrm{X}}_\rkN$-data,
where ${\mathrm{X}}$ is one of 
${\mathrm{A}},\ldots,{\mathrm{G}}$.
}}
\end{definition}

Using Lemma~\ref{lemma:techlm}, we can easily see 
\begin{lemma}\label{lemma:CaTylm}
Let $\orderrsysbase$ be a
rank-$\rkN$ Cartan data.
Let $\bhm\in\dotXNCartan$ be as in \eqref{eqn:bhmCtn}. 
Then $(\Rbhm,\bhm)$ is isomorphic to $(\mcR_\orderrsysbase,a)$,
where $a\in\mcA_\orderrsysbase$ {\rm{(}}note $|\mcA_\orderrsysbase|=1${\rm{)}}.
In particular, $|\Abhm|=1$, and identifying $\orderedPi$ with $\orderrsysbase$, 
we have $R(\bhm)=\rsystem$, where $\rsystem$ is the rank-$\rkN$ root system
corresponding to $\orderrsysbase$.
\end{lemma}

\begin{definition}\label{definition:defSuperbhm-d}{\rm{
Recall Definition~\ref{definition:defSuperbhm}.
Let $\dotXNSuper$ be the subset of $\bigchiPi$ formed by $\bhm$ with
\begin{equation}\label{eqn:defSuperbhm-eq}
\bhm(\al_i,\al_j)=(-1)^{\parity(i)\parity(j)}\bq^{\breta(\bal_i,\bal_j)}\quad(\,i,\,j\in\fkI\,)
\end{equation}
for some $\bq\in\bKtinf$ and some rank-$\rkN$ standard super-data 
$(\breta,\orderbarPi=(\bal_i|i\in\fkI))$,
where $\parity$ is the parity map associated with $(\breta,\orderbarPi)$;
moreover, assume $x\in\bZ\setminus\{0,-1\}$
if $\rkN=3$ and $(\breta,\orderbarPi)$ is the ${\mathrm{D}}(2,1;x)$-data.
Let $\dotXNSuper({\mathrm{X}})$ be the subset of $\dotXNSuper$ formed by
$\bhm$ as in \eqref{eqn:defSuperbhm-eq} for which $(\breta,\orderbarPi=(\bal_i|i\in\fkI))$
are the ${\mathrm{X}}$-data. 
}}
\end{definition}

\begin{figure}
\begin{center}
\figureDall
\end{center}
\caption{
All $\ospep\in\ospemfD_{2|2}$ and the Dynkin diagrams of all $\bhm\in\dotXNSuper({\mathrm{D}}(2,2))$, where $(z_1,z_2,z_3,z_4,z_5)\in\fkJ_{0,1}^5$
means $\ospep\in\ospemfD_{2|2}$
with $\ospep(i)=z_i$}
\label{fig:DynkinDall}
\end{figure}

We can directly see
\begin{lemma}\label{lemma:SuCaTylm}
{\rm{(1)}} {\rm{(}}See also Subsection~{\rm{\ref{subsection:LGEsuperA}}}.{\rm{)}}
Assume $\rkN\geq 2$. Let $m\in\fkJ_{1,\rkN}$ 
and
$\bhm\in\dotXNSuper({\mathrm{A}}(m-1,\rkN-m))$. Let 
$p:=p^+_{m|\rkN+1-m}$.
Let $\bq:=\bhm(\al_1,\al_2)^{-1}$ 
{\rm{(}}resp. $\bq:=\bhm(\al_{\rkN-1},\al_\rkN)${\rm{)}} if $m\geq 2$
{\rm{(}}resp. $m\leq\rkN-1${\rm{)}}.
Then $\bq\in\bKtinf$, and for $f\in\funcI_\infty$ and $t\in\bZgeqo$, letting
$\parity_{f,t}(k):=\delta(\breta^{p_{f,t}}(\al_k,\al_k),0)$
$(k\in\fkI)$,
we have
\begin{equation}\label{eqn:SuCa-eq-1}
\bhm_{f,t}(\al_i,\al_j)=(-1)^{\parity_{f,t}(i)\parity_{f,t}(j)}\bq^{\breta^{p_{f,t}}(\al_i,\al_j)}
\quad(i,j\in\fkI),
\end{equation} where $p_{f,t}$ is the one defined for the Cartan scheme $\superAmcC_{m|\rkN+1-m}$.
In particular, $(\Rbhm,\bhm)$ is isomorphic to $(\superAmcR_{m|\rkN+1-m},p)$.

{\rm{(2)}} {\rm{(}}See also Subsection~{\rm{\ref{subsection:LGEsuperB}}}.{\rm{)}}
Assume $\rkN\geq 1$. Let $m\in\fkJ_{1,\rkN}$
and
$\bhm\in\dotXNSuper({\mathrm{B}}(m,\rkN-m))$.
Let $p:=p^-_{\rkN-m|m}$.
Let $\bq:=\bhm(\al_\rkN,\al_\rkN)$ 
{\rm{(}}resp. $\bq:=-\bhm(\al_\rkN,\al_\rkN)^{-1}${\rm{)}} if $m\leq 1$
{\rm{(}}resp. $m=0${\rm{)}}.
Then $\bq\in\bKtinf$, and for $f\in\funcI_\infty$ and $t\in\bZgeqo$, 
letting
\begin{equation*}
\parity_{f,t}(k):=
\left\{\begin{array}{ll}
\delta(\breta^{p_{f,t}}(\al_k,\al_k),0) & \quad\mbox{if $k\in\fkJ_{1,\rkN-1}$}, \\
{\frac {1-\breta^{p_{f,t}}(\al_\rkN,\al_\rkN)} 2}  & \quad\mbox{if $k=N$} 
\end{array}\right.
\end{equation*} $(k\in\fkI)$, we have
\begin{equation}\label{eqn:SuCa-eq-2}
\bhm_{f,t}(\al_i,\al_j)=(-1)^{\parity_{f,t}(i)\parity_{f,t}(j)}\bq^{\breta^{p_{f,t}}(\al_i,\al_j)}
\quad(i,j\in\fkI),
\end{equation} where $p_{f,t}$ is the one defined for $\superBmcC_{m|\rkN-m}$.
In particular, $(\Rbhm,\bhm)$ is isomorphic to $(\superBmcR_{m|\rkN-m},p)$.

{\rm{(3)}} {\rm{(}}See also Subsection~{\rm{\ref{subsection:LGEsuperC}}}.{\rm{)}}
Assume $\rkN\geq 3$.
Let $\bhm\in\dotXNSuper({\mathrm{C}}(\rkN))$
and $p:=\ospepCN$.
Let $\bq:=\bhm(\al_1,\al_2)^{-1}$.
Then $\bq\in\bKtinf$, and for $f\in\funcI_\infty$ and $t\in\bZgeqo$, letting
$\parity_{f,t}(k):=\delta(\ospeeta^{p_{f,t}}(\al_k,\al_k),0)$
$(k\in\fkI)$,
we have
\begin{equation}\label{eqn:SuCa-eq-3}
\bhm_{f,t}(\al_i,\al_j)=(-1)^{\parity_{f,t}(i)\parity_{f,t}(j)}\bq^{-\ospeeta^{p_{f,t}}(\al_i,\al_j)}
\quad(i,j\in\fkI),
\end{equation} where $p_{f,t}$ is the one defined for the Cartan scheme $\ospemcC_{1|\rkN-1}$.
In particular, $(\Rbhm,\bhm)$ is isomorphic to $(\ospemcR_{1|\rkN-1},p)$. 

{\rm{(4)}} {\rm{(}}See also Subsection~{\rm{\ref{subsection:LGEsuperD}}}.{\rm{)}}
Assume $\rkN\geq 3$. Let $m\in\fkJ_{2,\rkN-1}$
and $\bhm\in\dotXNSuper({\mathrm{D}}(m,\rkN-m))$.
Let $p:=\ospep_{m|\rkN-m}$ and $\bq:=\bhm(\al_{\rkN-2},\al_{\rkN-1})^{-1}$.
Then $\bq\in\bKtinf$, and for $f\in\funcI_\infty$ and $t\in\bZgeqo$, letting
$\parity_{f,t}(k):=\delta(\ospeeta^{p_{f,t}}(\al_k,\al_k),0)$
$(k\in\fkI)$,
we have
\begin{equation}\label{eqn:SuCa-eq-4}
\bhm_{f,t}(\al_i,\al_j)=(-1)^{\parity_{f,t}(i)\parity_{f,t}(j)}\bq^{\ospeeta^{p_{f,t}}(\al_i,\al_j)}
\quad(i,j\in\fkI),
\end{equation} where $p_{f,t}$ is the one defined for the Cartan scheme $\ospemcC_{m|\rkN-m}$.
In particular, $(\Rbhm,\bhm)$ is isomorphic to $(\ospemcR_{m|\rkN-m},p)$.  
\end{lemma}

\begin{remark}\label{remark:superPRS}
Keep the notation as in Definition~\ref{definition:defSuperbhm-d},
so let $\bhm\in\dotXNSuper$ be as in \eqref{eqn:defSuperbhm-eq}.
Let $\barR$ be the standard root system associated with $(\breta,\orderbarPi)$. 
Define the $\bR$-linear isomorphism 
$\bartrincl:\oplus_{i\in\fkI}\bR\bal_i\to\bV$ by
$\bartrincl(\bal_i):=\al_i$.
Then 
\begin{equation}\label{eqn:eqRK}
R(\bhm)=\bartrincl(\barR\setminus 2\barR).
\end{equation}
In particular, $\bhm\in\bigchiPifin$. 
We can prove \eqref{eqn:eqRK} in a way similar to that for  \cite[Proposition~10.4.1~(i)]{Y94}.
It can also be proved as follows. 
If $\bhm\in\dotXNSuper({\mathrm{A}}(m-1,\rkN-1))$
(resp. $\bhm\in\dotXNSuper({\mathrm{B}}(m,\rkN-m))$),
it follows from Lemmas~\ref{lemma:longest-superA-id} and \ref{lemma:SuCaTylm}~(1)
(resp.  Lemmas~\ref{lemma:longest-superB-id} and \ref{lemma:SuCaTylm}~(2)).
If $\bhm\in\dotXNSuper({\mathrm{C}}(\rkN))\cup
\dotXNSuper({\mathrm{D}}(m,\rkN-m))$,
it follows from \eqref{eqn:defoseR}, \eqref{eqn:evenosp-1-1}
and Lemma~\ref{lemma:SuCaTylm}~(3),(4).
If $\bhm\in\dotXNSuper({\mathrm{F}}(4))$ with $\rkN=4$
or $\bhm\in\dotXNSuper({\mathrm{G}}(3))$ with $\rkN=3$, we can directly prove it
using Lemma~\ref{lemma:lgelmmm} and \eqref{eqn:hatfwo} below.
\end{remark}

\begin{definition} \label{definition:defExtrabhm}
{\rm{
If $\rkN\in\bN\setminus\fkJ_{2,4}$, let $\dotXNExtra
:=\emptyset$.
If $\rkN=2$
(resp. $\rkN=3$, resp. $\rkN=4$), let $\dotXNExtra$ be the set of bi-homomorphisms 
$\bhm\in\bigchiPi$ satisfying the condition \eqref{eqn:deEXrk2} 
(resp. \eqref{eqn:deEXrk3}, resp. \eqref{eqn:deEXrk4}) below.
In the following, let $q_{ij}:=\bhm(\al_i,\al_j)$.
\begin{equation}\label{eqn:deEXrk2}
q_{11}^2+q_{11}+1=0,\,q_{12}=q_{21}\in\bKtinf,\,q_{12}q_{21}q_{22}=1.
\end{equation}
\begin{equation}\label{eqn:deEXrk3}
\begin{array}{l}
q_{12}=q_{21}\in\bKtinf,\,q_{11}q_{12}q_{21}=1,\,q_{22}=-1,\,q_{13}=q_{31}=1,\\
q_{23}=q_{32},\,q_{23}q_{32}q_{33}=1,\,q_{33}\ne 1,\,q_{11}q_{33}\ne 1.
\end{array}
\end{equation}
\begin{equation}\label{eqn:deEXrk4}
\begin{array}{l}
q_{12}=q_{21}=q_{23}=q_{32}\in\bKtinf,\,q_{34}=q_{43},\\
q_{13}=q_{31}=q_{14}=q_{41}=q_{24}=q_{42}=1, \\
q_{33}=q_{11}q_{44}=-1,\,q_{11}q_{12}q_{21}=q_{12}q_{21}q_{22}=1,\,q_{34}q_{43}q_{44}=1.
\end{array}
\end{equation}
}}
\end{definition}

\begin{remark}\label{remark:extraPRS}
Note that 
\begin{equation}\label{eqn:DtwonSbExt}
\mbox{if $\rkN=3$, then $\cup_{x\in\bZ\setminus\{0,-1\}}\dotXNSuper({\mathrm{D}}(2,1;x))\subset\dotXNExtra$.}
\end{equation}
Let $\bhm\in\dotXNExtra$.
If $\rkN=2$,
then $R^+(\bhm)=\{\al_2,\al_1+\al_2,2\al_1+\al_2,\al_1\}$.
If $\rkN=4$,
then $R^+(\bhm)=\{\al_1,\al_3,\al_1+\al_2+\al_3,\al_2+\al_3,
\al_1+2\al_2+\al_3,\al_1+\al_2,\al_2\}$.
If $\rkN=3$,
then $R^+(\bhm)=\{\al_1,\al_1+\al_2,\al_2,\al_4,
\al_1+\al_2+\al_3+\al_4,\al_1+\al_2+\al_3,
\al_1+2\al_2+2\al_3+\al_4,
\al_1+\al_2+2\al_3+\al_4,
\al_1+2\al_2+3\al_3+2\al_4,
\al_1+2\al_2+3\al_3+\al_4,
\al_2+\al_3+\al_4,
\al_3+\al_4,
\al_2+2\al_3+\al_4,
\al_2+\al_3,
\al_3\}$.
We can directly prove these facts
using Lemma~\ref{lemma:lgelmmm} and \eqref{eqn:hatfwo} below.
\end{remark}

\setlength{\unitlength}{1mm}
\begin{picture}(120,165)(-10,-60)

\put(0,90){\figzero}

\put(20,78){\vector(0,1){10}}
\put(20,85){\vector(0,-1){10}}
\put(15,81){$\newtr_3$}

\put(0,60){\figone}

\put(52,66){\vector(1,-1){10}}
\put(60,58){\vector(-1,1){10}}
\put(56,64){$\newtr_4$}

\put(55,45){\figthree}

\put(20,48){\vector(0,1){10}}
\put(20,55){\vector(0,-1){10}}
\put(15,51){$\newtr_2$}

\put(105,33){\vector(0,1){10}}
\put(105,40){\vector(0,-1){10}}
\put(100,36){$\newtr_2$}

\put(0,30){\figtwo}
\put(-5,25){\vector(1,1){5}}
\put(-5,25){\line(0,-1){38}}
\put(-5,-13){\vector(1,-1){5}}
\put(-3,-5){$\newtr_1$}
\put(0,-30){\figfour}

\put(20,18){\vector(0,1){10}}
\put(20,25){\vector(0,-1){10}}
\put(15,21){$\newtr_4$}

\put(105,3){\vector(0,1){10}}
\put(105,10){\vector(0,-1){10}}
\put(100,6){$\newtr_4$}

\put(0,0){\figfive}

\put(67,06){\vector(1,-1){10}}
\put(75,-2){\vector(-1,1){10}}
\put(71,04){$\newtr_2$}

\put(70,-15){\figsix}

\put(70,15){\figseven}

\put(55,75){\figeight}

\put(120,65){\vector(-1,1){5}}
\put(120,65){\line(0,-1){30}}
\put(120,35){\vector(-1,-1){5}}
\put(115,38){$\newtr_1$}

\put(55,-45){\fignine}

\put(90,-27){\vector(0,1){10}}
\put(90,-20){\vector(0,-1){10}}
\put(85,-24){$\newtr_3$}

\put(30,-55){{\rm{Figure~$1$: Dynkin diagrams of ${{{\dot {\mathcal{X}}}^{\mathrm{Extra}}_4}}$}}}

\end{picture}

\subsection{Classification of FID-type bi-homomorphisms}

\begin{lemma}\label{lemma:fin}
Let $\bhm\in\bigchiPiirr$.
Let $q_{ij}:=\bhm(\al_i,\al_j)$ for $i$, $j\in\fkI$.
Then
$\dim U^-(\bhm)<\infty$ if and only if $\bhm\in\bigchiPifin$
and either {\rm{(x)}} or {\rm{(y)}} below holds.

{\rm{(x)}} $\rkN=1$ and $q_{11}\notin \bKtinf\cup \bKt_1$.

{\rm{(y)}} $\rkN\geq 2$ and $q_{ij}q_{ji}\notin \bKtinf$
for all $i$, $j\in\fkI$.
\end{lemma}
{\it Proof.}
By Theorem~\ref{theorem:Khth},
we see that
\begin{equation}\label{eqn:cndfn}
\dim U^-(\bhm)<\infty\,\,\Longleftrightarrow\,\,
\bhm\in\bigchiPifin\,\,\mbox{and}\,\,
\forall \al\in R^+(\bhm),\,
\bhm(\al,\al)\notin \bKtinf\cup \bKt_1.
\end{equation}
Hence, if $\rkN=1$, the statement is true.

Assume $\rkN\geq 2$ and $\bhm\in\bigchiPifin$.
Since $\bhm\in\bigchiPiirr$, by \eqref{eqn:subpp}
and Lemma~\ref{lemma:ichilm}~{(1)},
we have $q_{ii}\ne 1$ for all $i\in\fkI$.

Assume that (y) is not true.
If $q_{ii}\in\bKtinf$ for some $i\in\fkI$,
then $\dim U^-(\bhm)=\infty$ by \eqref{eqn:cndfn}.
Assume that there exist $i$, $j\in\fkI$ with $i\ne j$
such that $q_{ii}$, $q_{jj}\notin\bKtinf$
and $q_{ij}q_{ji}\in\bKtinf$.
Then $-c^\bhm_{ij}\in\bN$ by \eqref{eqn:cxij}
and \eqref{eqn:subpp}.
Let $\beta:=\al_i+\al_j$.
By \eqref{eqn:elprofRd}, $\beta\in R^+(\bhm)$.
Since $\bhm(\beta,\beta)\in\bKtinf$,
we have $\dim U^-(\bhm)=\infty$ by \eqref{eqn:cndfn}.

Assume that (y) is true.
Let $\al\in R^+(\bhm)$.
It is clear that $\bhm(\al,\al)\notin\bKtinf$.
By \eqref{eqn:ellRp}, Theorem~\ref{theorem:hecprone}~(2)
and Lemma~\ref{lemma:ichilm}~(3), there
exist $\bhmp\in\bigchiPifin\cap\bigchiPiirr$
and $i\in\fkI$ such that
$\bhmp\sim\bhm$
and $\bhmp(\al_i,\al_i)=\bhm(\al,\al)$.
By \eqref{eqn:subpp} and Lemma~\ref{lemma:ichilm}~(1),
since $\bhmp\in\bigchiPifin\cap\bigchiPiirr$,
we have
$\bhm(\al,\al)\ne 1$.
By \eqref{eqn:cndfn}, we have $\dim U^-(\bhm)<\infty$,
as desired.
This completes the proof.
\hfill $\Box$
\newline\par
For $\bhm$, $\bhmp\in\bigchiPi$,
we write $\bhm\approx \bhmp$ if
there exist $\bhm_1$, $\bhm_2\in\bigchiPi$
and a bijection $f:\fkI\to\fkI$ 
such that $\bhm\sim\bhm_1\equiv\bhm_2$
and 
$\bhmp(\al_i,\al_j)=\bhm_2(\al_{f(i)},\al_{f(j)})$
($i$, $j\in \fkI$).

By the Heckenberger's classification \cite[Tables~1-4, Theorems~17,~22]{Hec09}
and Lemma~\ref{lemma:fin}, we have

\begin{theorem}{\rm{(\cite{Hec09})}} \label{theorem:cl}
{\rm{(1)}} Assume $\rkN=1$. Let $\bhm\in\bigchiPi$.
Then $\dim U^-(\bhm)=\infty$ if and only if 
$\bhm(\al_1,\al_1)\in\bKtinf$ or $\bhm(\al_1,\al_1)=1$.

{\rm{(2)}} Assume $\rkN\geq 2$.
Then we have
\begin{equation}\label{eqn:cl-eqn}
\begin{array}{l}
\{\,\bhm\in\bigchiPiirr\cap\bigchiPifin\,|\,\dim U^-(\bhm)=\infty \,\} \\
\quad =\{\,\bhm\in\bigchiPi\,|\,\exists\bhmp\in\dotXNCartan\cup
\dotXNSuper\cup\dotXNExtra,\,\bhm\approx\bhmp\,\}.
\end{array}
\end{equation}
\end{theorem}

\begin{remark}\label{remark:smbhms}
We have $\dotXNCartan({\mathrm{B}}_\rkN)=\dotXNSuper({\mathrm{B}}(0,\rkN))$.
We have the bijection
$\dotXNSuper({\mathrm{A}}(m-1,\rkN-m))\to\dotXNSuper({\mathrm{A}}(\rkN-m,m-1))$,
$\bhm\mapsto\bhmp$,
defined by
$\bhmp(\al_i,\al_j):=\bhm(\al_{\rkN-i+1},\al_{\rkN-j+1})$
($i,j\in\fkI$), 
where we also have $\bhm\approx\bhmp$, see Figure~\ref{fig:DynkinA}. 
If $\rkN=3$, 
then we have \eqref{eqn:DtwonSbExt} and
$\dotXNSuper({\mathrm{D}}(2,1))=\dotXNSuper({\mathrm{D}}(2,1;1))$,
and see that for $\bhm$, $\bhmp\in\dotXNExtra$,
$\bhm\approx\bhmp$ if and only if 
$\{q_{11},q_{33},(q_{11}q_{33})^{-1}\}
=\{q^\prime_{11},q^\prime_{33},(q^\prime_{11}q^\prime_{33})^{-1}\}$,
where $q_{ij}:=\bhm(\al_i,\al_j)$ and $q^\prime_{ij}:=\bhmp(\al_i,\al_j)$,
see also Figure~\ref{fig:DynkinDX}.
\end{remark}

\section{Lusztig isomorphisms} \label{section:Liso}
\subsection{Lusztig isomprphisms of the generalized quantum groups}
In this section, fix $i\in\fkI$
and $\bhm\in(\ppbigchiPifin)_i$,
and let $q_{ij}:=\bhm(\al_i,\al_j)$ ($j\in\fkI$).

\begin{theorem}{\rm{(\cite[Theorem~6.11]{Hec10})}}\label{theorem:Lusiso}
Assume $\bhm\in(\ppbigchiPifin)_i$.
Then there exists a unique $\bK$-algebra isomorphism
\begin{equation*}
T_i:=T^{\newtr_i\bhm}_i:U(\newtr_i\bhm)\to U(\bhm)
\end{equation*} such that
\begin{equation*}
\begin{array}{l}
T_i(K_\al)=K_{s^{\newtr_i\bhm}_i(\al)},\quad T_i(L_\al)=L_{s^{\newtr_i\bhm}_i(\al)},\\
T_i(E_i)=F_iL_{-\al_i},\quad T_i(F_i)=K_{-\al_i}E_i,\\
T_i(E_j)=E_{-c^\bhm_{ij},\al_i,\al_j}, \\
T_i(F_j)={\frac 1 {(-c^\bhm_{ij})_{q_{ii}}!(-c^\bhm_{ij};q_{ii},
q_{ij}q_{ji})!}}F_{-c^\bhm_{ij},\al_i,\al_j}
\end{array}
\end{equation*}
for $\al\in\bZPi$ and $j\in \fkI\setminus\{i\}$.
In particular,
\begin{equation}\label{eqn:Tidd}
T_i(U(\newtr_i\bhm)_\al)=U(\bhm)_{s^\bhm_i(\al)}\quad\quad (\al\in\bZPi).
\end{equation}
\end{theorem}

\subsection{Lusztig isomorphisms between irreducible modules}
\label{subsection:LisMo}

Let
$\Lambda\in\rmCh(U^0(\bhm))$, and $h:=h_{\bhm,\Lambda,i}$. Assume
that $\bhm\in(\ppbigchiPifin)_i$ and $h\ne\infty$.
Define $\newtr_i \Lambda :=\newtr_i^\bhm \Lambda\in \rmCh(U^0(\newtr_i\bhm))$ by
\begin{equation}\label{eqn:dfilm}
\newtr_i^\bhm \Lambda(K_\al L_\beta):=\Lambda(K_{s^{\newtr_i\bhm}_i(\al)}L_{s^{\newtr_i\bhm}_i(\beta)})
{\frac {\bhm(\al_i,s^{\newtr_i\bhm}_i(\beta))^h} {\bhm(s^{\newtr_i\bhm}_i(\al),\al_i)^h}}
\quad (\al,\,\beta\in\bZPi).
\end{equation} By \eqref{eqn:hixl} and \eqref{eqn:ichi}, we have
\begin{equation} \label{eqn:lgelmrd}
h_{\newtr_i \bhm,\newtr_i \Lambda,i}=h,\quad\mbox{and}\quad \newtr_i^{\newtr_i\bhm} \newtr_i^\bhm \Lambda=\Lambda.
\end{equation} By \eqref{eqn:hchcp}, \eqref{eqn:cijccpd}  and \eqref{eqn:dfilm},
for $\bhmp\in\bigchiPi$
and $\Lambda^\prime\in\rmCh(U^0(\bhmp))$, we have

\begin{equation} \label{eqn:lgelmrdd}
(\bhm,\Lambda)\equiv(\bhmp,\Lambda^\prime)\,\,\Longrightarrow\,\,
(\newtr_i \bhm,\newtr_i^\bhm \Lambda)\equiv(\newtr_i \bhmp,\newtr_i^{\bhmp} \Lambda^\prime),
\end{equation} where we have $\bhmp\in(\ppbigchiPifin)_i$.

\begin{lemma} \label{lemma:lgelmdd} Let $\Lambda$ and $h$ be as above.
Assume $h<\infty$.
There exists
a unique $\bK$-linear isomorphism
\begin{equation*}
\hT_i:=\hT^{\newtr_i \bhm,\newtr_i \Lambda}_i:\mclL_{\newtr_i\bhm}(\newtr_i^\bhm \Lambda)\to\mclL_\bhm(\Lambda)
\end{equation*} such that
\begin{equation} \label{eqn:lgelmrvl}
\hT_i (X v_{\newtr_i\bhm})=T_i(X)F_i^hv_\Lambda\quad
(X\in U(\newtr_i \bhm)).
\end{equation}

\end{lemma}
{\it {Proof.}} We can regard
$\mclL_\bhm(\Lambda)$ as a left $U(\newtr_i \bhm)$-module
defined by $X\cdot u:=T_i(X)u$ ($X\in U(\newtr_i\bhm)$, $u\in\mclL_\bhm(\Lambda)$).
Let $v^\prime:=F_i^hv_\Lambda\in\mclL_\bhm(\Lambda)$.
Note that a $U(\bhm)$-submodule of $\mclL_\bhm(\Lambda)$
is a $U(\newtr_i\bhm)$-submodule, and vice versa.
 By \eqref{eqn:fvo},
$v^\prime\ne 0$ and $E_i\cdot v^\prime=0$, so we also have
$E_j\cdot v^\prime=0$ for $j\in\fkI\setminus\{i\}$.
Then
we have a $U(\newtr_i \bhm)$-module homomorphism
$z:\mclM_{\newtr_i \bhm}(\newtr_i\Lambda)\to\mclL_\bhm(\Lambda)$ such that
$z(X\tv_{\newtr_i \Lambda})= X\cdot v^\prime$  for $X\in U(\newtr_i\bhm)$.
By \eqref{eqn:esyf}, $z$ is surjective,
and we also have $z(\mclM_{\newtr_i\bhm}(\newtr_i \Lambda)_\al)
=\mclL_\bhm(\Lambda)_{s_i(\al)-h\al_i}$
for $\al\in\bZgeqoPi$.
In particular, $\ker z$ is a proper
$\bZPi$-graded left $U(\newtr_i\bhm)$-submodule
of $\mclM_{\newtr_i \bhm}(\newtr_i\Lambda)$,
so $\ker z\subseteq\cN_{\newtr_i\bhm}(\newtr_i\Lambda)$.
By \eqref{eqn:fvo} and \eqref{eqn:lgelmrd},
we have $\mclM_{\newtr_i\bhm}(\newtr_i \Lambda)_{-h\al_i}
\cap \cN_{\newtr_i\bhm}(\newtr_i\Lambda)=\{0\}$
since $\mclM_{\newtr_i\bhm}(\newtr_i \Lambda)_{-h\al_i}=F_i^h\tv_{\newtr_i \Lambda}$.
Therefore, $z(\cN_{\newtr_i\bhm}(\newtr_i\Lambda))$ is a proper
$\bZPi$-graded $U(\bhm)$-submodule of
$\mclL_\bhm(\Lambda)$.
Hence $z(\cN_{\newtr_i\bhm}(\newtr_i\Lambda))=\{0\}$,
which implies $\cN_{\newtr_i\bhm}(\newtr_i\Lambda)=\ker z$.
Hence $z$ induces a $U(\newtr_i\bhm)$-module
isomorphism $\hT_i$,
as desired.
This completes the proof. \hfill $\Box$
\newline\par
Using  Lemma~\ref{lemma:crtyp},
together with \eqref{eqn:hchcp} and \eqref{eqn:subp},
the following lemme is an easy exercise for the reader.

\begin{lemma} \label{lemma:crtypmo}
Let
$\bhm\in\bigchiPifin$.
Assume that
$q_{ii}\in \bKtinf$
and $h_{\bhm,\Lambda,i}\ne\infty$.
Then
$(\newtr_i\bhm,\newtr_i^\bhm \Lambda)\equiv (\bhm,\Lambda)$,
and $h_{\newtr_i\bhm,\newtr_i^\bhm\Lambda,j}=h_{\bhm,\Lambda,j}$
for $j\in \fkI$.
\end{lemma}

\begin{definition}\label{definition:fchi}{\rm{
(Definition of  $H(\bhm,\Lambda,f)$)
Let $\bhm\in\bigchiPifin$ and
$n:=|R^+(\bhm)|$. Let $f\in\funcI_n$
be such that $\bhms_{f,n}=\lgst$,
where note that the existence 
of $\lgst$ follows from Lemma~{\rm{\ref{lemma:lgslm}~(2)}}.
Recall Notation~\ref{notation:Ab-chift-1}.
Let $\Lambda_{\bhm,f,0}:=\Lambda$.
If $t\in\fkJ_{1,n}$, we define
$\Lambda_{\bhm,f,t}:=\newtr_{f(t)}^{\bhm_{f,t-1}}\Lambda_{\bhm,f,t-1}$
if $\Lambda_{\bhm,f,t-1}$ can be defined and
$h_{\bhm_{f,t-1},\Lambda_{\bhm,f,t-1},f(t)}<\infty$.
Define
$H(\bhm,\Lambda,f)\in\fkJ_{0,n}$ as follows.
If there exists $t\in\fkJ_{0,n-1}$
such that $h_{\bhm_{f,k-1},\Lambda_{\bhm,f,k-1},f(k)}<\infty$
for all $k\in\fkJ_{1,t}$
and $h_{\bhm_{f,t},\Lambda_{\bhm,f,t},f(t+1)}=\infty$,
let $H(\bhm,\Lambda,f):=t$.
If $h_{\bhm_{f,k-1},\Lambda_{\bhm,f,k-1},f(k)}<\infty$
for all $k\in\fkJ_{1,n}$,
let $H(\bhm,\Lambda,f):=n$.
}}
\end{definition}

The following lemma is a crucial key in this paper.

\begin{lemma} \label{lemma:lgelmd}
Keep the notation as in Definition~{\rm{\ref{definition:fchi}}}.
Then
$\dim\mclL_\bhm(\Lambda)<\infty$
if and only if $H(\bhm,\Lambda,f)=n$.
\end{lemma}
{\it{Proof.}} Assume $\dim\mclL_\bhm(\Lambda)<\infty$.
Then $h_{\bhm_{f,0},\Lambda_{f,0},f(1)}<\infty$.
By Lemma~\ref{lemma:lgelmdd},
$\dim\mclL_{\bhm_{f,1}}(\Lambda_{f,1})<\infty$.
Repeating this argument, we have $H(\bhm,\Lambda,f)=n$.

Assume $H(\bhm,\Lambda,f)=n$.
Let $\gamma:=\sum_{t\in\fkJ_{1,n}}h_{\bhm_{f,t-1},\Lambda_{f,t-1},f(t)}\bhms_{f,t-1}(\al_{f(t)})$.
By \eqref{eqn:ellRp}, $\gamma\in\bZgeqoPi$.
By \eqref{eqn:lgelmrvl}, $\hT_{f(1)}\cdots\hT_{f(n)}(v_{\Lambda_{\bhm,f,n}})\in\mclL_\bhm(\Lambda)_{-\gamma}$.
Let 
\begin{equation*}
X:=\{\beta\in\bZgeqoPi\,|\,\gamma-\beta\in\bZgeqoPi\}.
\end{equation*} Then $|X|<\infty$.
We have
\begin{equation*}
\begin{array}{lll}
\mclL_\bhm(\Lambda)&=&\hT_{f(1)}\cdots\hT_{f(n)}(\mclL_{\bhm_{f,n}}(\Lambda_{\bhm,f,n})) \\
&=& \hT_{f(1)}\cdots\hT_{f(n)}(U^-(\bhm_{f,n})v_{\Lambda_{f,n}}) \\
&=& \oplus_{\al\in\bZgeqoPi}\hT_{f(1)}\cdots\hT_{f(n)}(U^-(\bhm_{f,n})_{-\al} v_{\Lambda_{\bhm,f,n}}) \\
&=& \oplus_{\al\in\bZgeqoPi}\hT_{f(1)}\cdots\hT_{f(n)}(U(\bhm_{f,n})_{-\al} v_{\Lambda_{\bhm,f,n}}) \\
&=& \oplus_{\al\in\bZgeqoPi}T_{f(1)}\cdots T_{f(n)}(U(\bhm_{f,n})_{-\al})\hT_{f(1)}\cdots\hT_{f(n)}(v_{\Lambda_{\bhm,f,n}}) \\
& & \quad\mbox{(by \eqref{eqn:lgelmrvl})} \\
&=& \oplus_{\beta \in\bZgeqoPi}U(\bhm)_\beta \hT_{f(1)}\cdots\hT_{f(n)}(v_{\Lambda_{\bhm,f,n}}) \\
& & \quad\mbox{(by \eqref{eqn:lgsta} and \eqref{eqn:Tidd})} \\
&=& \oplus_{\beta \in\bZgeqoPi}\mclL_\bhm(\Lambda)_{\beta-\gamma} \\
&=& \oplus_{\beta \in X}\mclL_\bhm(\Lambda)_{\beta-\gamma} \\
& & \quad\mbox{(since $\mclL_\bhm(\Lambda)=\oplus_{\al \in\bZgeqoPi}\mclL_\bhm(\Lambda)_{-\al}$)}.
\end{array}
\end{equation*} Hence
$\dim\mclL_\bhm(\Lambda)=\sum_{\beta \in X}\dim\mclL_\bhm(\Lambda)_{\beta-\gamma}<\infty$, as desired.
\hfill $\Box$
\newline

\begin{figure}
\begin{center}
\figHWM
\end{center}
\caption{${\mbox{$\mclL_\bhm(\Lambda)$ for 
$\bhmp\in\dotXNSuper({\mathrm{A}}(m-1,\rkN-m))$
with $\rkN=2$, $m=1$,}\quad\quad\quad\quad\quad\quad\quad\quad
\atop \mbox{$\lambda_1\ne -1$, $\lambda_2= \bq^8 $, $\lambda_1\lambda_2\bq^2\ne -1$,
where $\lambda_i:=\Lambda(K_{\al_i}L_{-\al_i})$}} $ }
\label{fig:HWM}
\end{figure}

The following lemma follows from \eqref{eqn:hchcp}, \eqref{eqn:lgelmrdd}, and
Lemmas~\ref{lemma:lgelmd} and \ref{lemma:prlgelmmm-d}.
\begin{lemma} \label{lemma:eqll}
Assume $\bhm\in\bigchiPifin$.
Let $\bhmp\in\bigchiPifin$ and
$\Lambda^\prime\in\rmCh(U^0)$
be such that $(\bhmp,\Lambda^\prime)\equiv (\bhm,\Lambda)$.
Then
$\dim\mclL_{\bhmp}(\Lambda^\prime)<\infty$
if and only if $\dim\mclL_\bhm(\Lambda)<\infty$.
\end{lemma}

\section{Main theorem}\label{section:MainTh}
\subsection{Irreducible weights for Cartan  and super-$AC$ cases}
\label{subsection:IrrMdCartan}
For $\bhm\in\bigchiPi$ and $i\in\fkI$, let
\begin{equation*}
\mbbS_i(\bhm):=\{\,\Lambda\in\rmCh(U^0(\bhm))\,|\,
\exists r\in\bZgeqo,\,\Lambda(K_{\al_i}L_{-\al_i})=\bhm(\al_i,\al_i)^r\,\}.
\end{equation*}

\begin{theorem}\label{theorem:IrrRepRkOne} 
Assume $\rkN=1$. Let $\bhm\in\bigchiPi$ be such that 
$\dim U^-(\bhm)=\infty$.
Let $\Lambda\in\rmCh(U^0(\bhm))$. 
Then $\dim\mclL_\bhm(\Lambda)<\infty$
if and only if $\Lambda\in\mbbS_1(\bhm)$.
{\rm{(}}See also Theorem~{\rm{\ref{theorem:cl}~(1)}}.{\rm{)}}
\end{theorem}
{\it{Proof.}} This follows from \eqref{eqn:fvo}.
\hfill $\Box$

\begin{theorem}\label{theorem:IrrRepCartan} 
Assume $\rkN\geq 2$. Let $\bhm\in\bigchiPi$ be such that $\bhm\approx\bhmp$
for some $\bhmp\in\dotXNCartan$.
Let $\Lambda\in\rmCh(U^0(\bhm))$. 
Then $\dim\mclL_\bhm(\Lambda)<\infty$
if and only if $\Lambda\in\cap_{i\in\fkI}\mbbS_i(\bhm)$.
\end{theorem}
{\it{Proof.}} By Theorem~\ref{theorem:cl}, 
$\bhm\in\bigchiPifin$.
Note that $\bhm(\al_i,\al_i)\in\bKtinf$ for all $i\in\fkI$.
Then this theorem follows from \eqref{eqn:hixl}, \eqref{eqn:fvod} and
Lemmas~\ref{lemma:crtypmo} and \ref{lemma:lgelmd}.
\hfill $\Box$

\begin{theorem}\label{theorem:IrrRepSuperAC} 
{\rm{(1)}}
Assume $\rkN\geq 2$. Let $m\in\fkI$.
Let $\bhm\in\bigchiPi$ be such that $\bhm\equiv\bhmp$
for some $\bhmp\in\dotXNSuper({\mathrm{A}}(m-1,\rkN-m))$.
Let $\Lambda\in\rmCh(U^0(\bhm))$. Then $\dim\mclL_\bhm(\Lambda)<\infty$
if and only if $\Lambda\in\cap_{i\in\fkI\setminus\{m\}}\mbbS_i(\bhm)$.

{\rm{(2)}} Assume $\rkN\geq 3$.
Let $\bhm\in\bigchiPi$ be such that $\bhm\equiv\bhmp$
for some $\bhmp\in\dotXNSuper({\mathrm{C}}(\rkN))$. 
Let $\Lambda\in\rmCh(U^0(\bhm))$.
Then $\dim\mclL_\bhm(\Lambda)<\infty$
if and only if $\Lambda\in\cap_{i\in\fkI\setminus\{1\}}\mbbS_i(\bhm)$.
\end{theorem}
{\it{Proof.}} (1) By Theorem~\ref{theorem:cl}, we see that $\bhm$, 
$\bhmp\in\bigchiPifin$.
By Lemma~\ref{lemma:eqll}, we may assume $\bhm=\bhmp$.
By \eqref{eqn:hixl} and \eqref{eqn:fvod},
we see that the `only-if' part holds.

We show the `if' part
Let $n$, $r\in\bN$ and $f\in\funcI_n$ be as in Proposition~\ref{proposition:ANNlgest}.
By Lemmas~\ref{lemma:techlm}, \ref{lemma:longest-superA-id} and 
\ref{lemma:SuCaTylm}~(1), we have $\bhms_{f,n}=\lgst$.
By Lemma~\ref{lemma:rtparity} and 
\ref{lemma:SuCaTylm}~(1) and \eqref{eqn:newtribhmd}, \eqref{eqn:lsA-id-1}, \eqref{eqn:defSuperbhm-eq}, 
we see that 
\begin{equation}\label{eqn:IRSuperA-1}
\bhm_{f,k-1}(\al_k,\al_k)\in\bKtinf\quad(k\in\fkJ_{1,r}), 
\end{equation} and
\begin{equation}\label{eqn:IRSuperA-2}
\bhm_{f,t-1}(\al_t,\al_t)=-1\quad(t\in\fkJ_{r+1,n}).
\end{equation}
By \eqref{eqn:hixl}, \eqref{eqn:IRSuperA-1} and Lemma~\ref{lemma:crtypmo}, we also see that
$\Lambda\in\cap_{i\in\fkI\setminus\{m\}}\mbbS_i(\bhm)$ implies $H(\bhm,\Lambda,f)\geq r$.
By \eqref{eqn:hixl}, \eqref{eqn:IRSuperA-2} and Lemma~\ref{lemma:crtypmo}, 
$H(\bhm,\Lambda,f)\geq r$ must be $H(\bhm,\Lambda,f)=n$.
Thus the `if' part follows from Lemma~\ref{lemma:lgelmd}. 
This completes the proof of the claim~(1).

(2) We can prove the claim~(2) in the same way as that for the claim~(1)
by using Propositions~\ref{proposition:ospCNlgnst}
and Lemma~\ref{lemma:SuCaTylm}~(3). 
\hfill $\Box$

\subsection{Some technical maps}\label{subsection:TechMap}
In Section~\ref{section:MainTh},
for $\lambda\in(\bKt)^\rkN$ and $i\in\fkI$,
let $\lambda_i$ mean the $i$-th component of $\lambda$, that is,
$\lambda=(\lambda_1,\ldots,\lambda_\rkN)$.

In Subsection~\ref{subsection:TechMap}, assume $\rkN\geq 2$ and 
let $\bq\in\bKtinf$ and $m\in\fkJ_{1,\rkN-1}$.
Let $\mcQ_\bq:=\{\bq^x|x\in\bZgeqo\}$.
Let
\begin{equation}\label{eqn:dfKmq}
\mcK^{(m)}_\bq:=\{\,\lambda\in(\bKt)^\rkN\,|\,
\lambda_i\in\mcQ_\bq\,(i\in\fkJ_{\rkN-m+1,\rkN})\,\}.
\end{equation}
Define the maps $\toolmap^{(\bq,m)}_k:\mcK^{(m)}_\bq\to\bKt$
($k\in\fkJ_{\rkN-m,\rkN}$) and $\tltoolmap^{(\bq,m)}_\rkN:\mcK^{(m)}_\bq\to\bKt$ 
by 
\begin{equation} \label{eqn:pretoolMain}
\begin{array}{l}
\toolmap^{(\bq,m)}_{\rkN-m}(\lambda):=\lambda_{\rkN-m}, \\
\toolmap^{(\bq,m)}_k(\lambda):=\lambda_k\toolmap^{(\bq,m)}_{k-1}(\lambda)\bq^{2(1-\delta(1,\toolmap^{(\bq,m)}_{k-1}(\lambda)))}
\,\,(k\in\fkJ_{\rkN-m+1,\rkN}), \\
\tltoolmap^{(\bq,m)}_\rkN(\lambda):={\frac {\lambda_\rkN} {\lambda_{\rkN-1}}}
\toolmap^{(\bq,m)}_{\rkN-1}(\lambda)^2\bq^{4(1-\delta(1,\toolmap^{(\bq,m)}_{\rkN-1}(\lambda)))}.
\end{array}
\end{equation}

The following lemma is used in the proofs of Theorems~\ref{theorem:IrrRepSuperB}
and \ref{theorem:IrrRepSuperD} below.

\begin{lemma} \label{lemma:toolMain}
{\rm{(1)}} 
Let $z\in\bKt\setminus\mcQ_\bq$. Then
\begin{equation}\label{eqn:tlMain-1}
(\toolmap^{(\bq,m)}_\rkN)^{-1}(\{z\})=\{\,\lambda\in\mcK^{(m)}_\bq\,|
\,\prod_{i=\rkN-m}^\rkN\lambda_i=
\bq^{-2m}z\,\}.
\end{equation}
 
{\rm{(2)}} 
Let $z\in\bKt\setminus\mcQ_\bq$. Then
\begin{equation}\label{eqn:tlMain-2}
\begin{array}{l}
(\tltoolmap^{(\bq,m)}_\rkN)^{-1}(\{z\}) \\
\quad=\{\,\lambda\in\mcK^{(m)}_\bq\,|\,
\prod_{j=\rkN-m}^{\rkN-1}\lambda_j\ne\bq^{-2(m-1)},\,
\lambda_{\rkN-1}\lambda_\rkN\prod_{i=\rkN-m}^{\rkN-2}\lambda_i^2=
\bq^{-4m}z\,\} \\
\quad\quad\cup
\{\,\lambda\in\mcK^{(m)}_\bq\,|\,
\prod_{j=\rkN-m}^{\rkN-1}\lambda_j=\bq^{-2(m-1)},\,
\lambda_\rkN=
\lambda_{\rkN-1}z\,\}.
\end{array}
\end{equation} 

{\rm{(3)}} 
We have 
\begin{equation}\label{eqn:tlMain-3}
\begin{array}{l}
(\toolmap^{(\bq,m)}_\rkN)^{-1}(\{1\}) \\
\,\,=\cup_{t=0}^m\{\,\lambda\in\mcK^{(m)}_\bq\,|\,\prod_{i=\rkN-m}^{\rkN-m+t}\lambda_i=
\bq^{-2t},\,
\lambda_j=1\,(j\in\fkJ_{\rkN-m+t+1,\rkN})\,\}.
\end{array}
\end{equation}

{\rm{(4)}} We have
\begin{equation}\label{eqn:tlMain-4}
\begin{array}{l}
(\tltoolmap^{(\bq,m)}_\rkN)^{-1}(\{1\}) \\
\,\,=\{\,\lambda\in\mcK^{(m)}_\bq\,|\,\lambda_i=1\,(i\in\fkJ_{\rkN-m,\rkN})\,\} \\
\quad\,\,\cup(\cup_{t=1}^{m-2}\{\,\lambda\in\mcK^{(m)}_\bq\,|\,\prod_{j=\rkN-m}^{\rkN-m+t}\lambda_j=
\bq^{-2t},\,
\lambda_i=1\,(i\in\fkJ_{\rkN-m+t+1,\rkN})\,\}) \\
\quad\,\,\cup \{\,\lambda\in\mcK^{(m)}_\bq\,|\,\prod_{j=\rkN-m}^{\rkN-1}\lambda_j=
\bq^{-2(m-1)},\,
\lambda_{\rkN-1}=\lambda_\rkN\,\} \\
\quad\,\,\cup \{\,\lambda\in\mcK^{(m)}_\bq\,|\,\prod_{j=\rkN-m}^{\rkN-1}\lambda_j\ne\bq^{-2(m-1)},\,
\lambda_{\rkN-1}\lambda_\rkN\prod_{i=\rkN-m}^{\rkN-2}\lambda_i^2=
\bq^{-4m}\,\}.
\end{array}
\end{equation}

\end{lemma}
{\it{Proof.}}
In this proof, we fix $\lambda\in\mcK^{(m)}_\bq$, and use the following notations.
For $k\in\fkJ_{\rkN-m-1,\rkN+1}$, let
\begin{equation*}
c_k(\lambda):=
\left\{\begin{array}{ll}
0 & \quad\mbox{if $k=\rkN-m-1$}, \\
|\{\,r\in\fkJ_{\rkN-m,k}\,|\,\toolmap^{(\bq,m)}_r(\lambda)\ne 1\,\}|
& \quad\mbox{if $k\in\fkJ_{\rkN-m,\rkN}$}, \\
c_\rkN(\lambda) & \quad\mbox{if $k=\rkN+1$}.
\end{array}\right.
\end{equation*}
Let
\begin{equation*}
g_k(\lambda):=\prod_{j=\rkN-m}^k\lambda_j\,\,(k\in\fkJ_{\rkN-m,\rkN}),
\quad\mbox{and}\quad 
{\tilde g}_N(\lambda):=\lambda_{\rkN-1}\lambda_\rkN\prod_{j=\rkN-m}^{\rkN-2}\lambda_j^2.
\end{equation*}
Let $t(\lambda):=\delta(g_{\rkN-1}(\lambda),\bq^{-2(m-1)})\in\fkJ_{0,1}$.
Let 
\begin{equation*}
r(\lambda):=
\left\{\begin{array}{l}
{\rm{Min}}\{\,x\in\fkJ_{\rkN-m,\rkN}\,|\,\toolmap^{(\bq,m)}_x(\lambda)= 1\,\} \\
\quad\mbox{if $\toolmap^{(\bq,m)}_y(\lambda)= 1$ for some $y\in\fkJ_{\rkN-m,\rkN}$}, \\
\rkN+1 \quad\mbox{otherwise}.
\end{array}\right.
\end{equation*}
Note that
\begin{equation}\label{eqn:toolpf-(-1)}
\begin{array}{l}
r(\lambda)-(\rkN-m)=c_{r(\lambda)-1}(\lambda)=c_{r(\lambda)}(\lambda),
\,\mbox{and that} \\
\mbox{for $y\in\fkJ_{r(\lambda)+1,\rkN+1}$},\,
c_{r(\lambda)}(\lambda)\leq c_y(\lambda),\,\mbox{and} \\
c_{r(\lambda)}(\lambda)=c_y(\lambda),\,\mbox{if and only if}\,
\lambda_z=1\,(z\in\fkJ_{r(\lambda)+1,y}).
\end{array}
\end{equation}

We have \eqref{eqn:toolpf-o-(i)} and \eqref{eqn:toolpf-o-(ii)} below.
\begin{equation}\label{eqn:toolpf-o-(i)}
\mbox{For $i\in\fkJ_{\rkN-m+1,\rkN}$, $\lambda_i=\bq^{l_i}$ for some $l_i\in\bZgeqo$.}
\end{equation}
\begin{equation}\label{eqn:toolpf-o-(ii)}
\mbox{For $k\in\fkJ_{\rkN-m+1,\rkN}$, $\toolmap^{(\bq,m)}_k(\lambda)=\bq^{2c_{k-1}(\lambda)}g_k(\lambda)$.}
\end{equation}
By \eqref{eqn:toolpf-o-(i)}, we can easily see
\eqref{eqn:toolpf-o-(iv)} below.
\begin{equation}\label{eqn:toolpf-o-(iv)}
\mbox{
For $k\in\fkJ_{\rkN-m,\rkN}$,
$g_k(\lambda)=\bq^{-2(k-(\rkN-m))}$ 
if and only if $r(\lambda)=k$.} 
\end{equation}
By \eqref{eqn:toolpf-o-(i)}, \eqref{eqn:toolpf-o-(ii)} 
and \eqref{eqn:toolpf-o-(iv)}, we easily see the claims (1) and (3).

Since $c_{\rkN-1}(\lambda)=c_{\rkN-2}(\lambda)+\delta(1,\toolmap^{(\bq,m)}_{\rkN-1}(\lambda))$, 
by \eqref{eqn:toolpf-o-(ii)},
we have
\begin{equation}\label{eqn:toolpf-o-(iii)}
\tltoolmap^{(\bq,m)}_\rkN(\lambda)=\bq^{4c_{\rkN-1}(\lambda)}{\tilde g}_\rkN(\lambda).
\end{equation}

By \eqref{eqn:toolpf-(-1)} and \eqref{eqn:toolpf-o-(iv)}, 
if $r(\lambda)\in\fkJ_{\rkN-m,\rkN-1}$, then
\begin{equation}\label{eqn:toolpf-iii-2}
{\tilde g}_\rkN(\lambda)= 
\left\{\begin{array}{ll}
\bq^{-4c_{r(\lambda)}(\lambda)}
\lambda_{\rkN-1}\lambda_\rkN\prod_{i=r(\lambda)+1}^{\rkN-2}\lambda_i^2
& \mbox{if $r(\lambda)\in\fkJ_{\rkN-m,\rkN-2}$}, \\
\bq^{-4c_{r(\lambda)}(\lambda)}
{\frac {\lambda_\rkN} {\lambda_{\rkN-1}}}
 &
\mbox{if $r(\lambda)=\rkN-1$}. 
\end{array}\right.
\end{equation}

We show the claim (2).
Let $Y_1$ and $Y_2$ be LHS and RHS of \eqref{eqn:tlMain-2}.

Let $\lambda\in Y_2$.
Then
\begin{equation}\label{eqn:toolpf-0}
z=\bq^{4(m-t(\lambda))}{\tilde g}_\rkN(\lambda).
\end{equation}
Assume 
$r(\lambda)\in\fkJ_{\rkN-m,\rkN-2}$.
By \eqref{eqn:toolpf-(-1)}, \eqref{eqn:toolpf-iii-2}
and \eqref{eqn:toolpf-0}, 
\begin{equation}\label{eqn:toolpf-0b}
z=\bq^{4(\rkN-t(\lambda)-r(\lambda))}
\lambda_{\rkN-1}\lambda_\rkN(\prod_{i=r(\lambda)+1}^{\rkN-2}\lambda_i^2).
\end{equation}
Since $t(\lambda)\in\fkJ_{0,1}$, by \eqref{eqn:toolpf-o-(i)} and \eqref{eqn:toolpf-0b}, 
$z\in\mcQ_\bq\setminus\{1\}$, contradiction.
Hence 
$r(\lambda)\in\fkJ_{\rkN-1,\rkN+1}$. Hence
$c_{\rkN-2}(\lambda)=m-1$.
By \eqref{eqn:toolpf-o-(iv)}, $c_{\rkN-1}(\lambda)=m-t(\lambda)$.
By \eqref{eqn:toolpf-o-(iii)} and \eqref{eqn:toolpf-0},
we have $\lambda\in Y_1$.
Hence $Y_2\subseteq Y_1$. 
Let $\lambda\in Y_1$, i.e., $\tltoolmap^{(\bq,m)}_\rkN(\lambda)=z$.
Assume $r(\lambda)\in\fkJ_{\rkN,\rkN+1}$. By \eqref{eqn:toolpf-o-(iii)},
${\tilde g}_\rkN(\lambda)=\bq^{-4m}z$.
Since $r(\lambda)\ne\rkN-1$, 
by \eqref{eqn:toolpf-o-(iv)},
$\bq^{2(m-1)}g_{\rkN-1}(\lambda)\ne 1$.
Hence $\lambda\in Y_2$. 
Assume $r(\lambda)\in\fkJ_{\rkN-m,\rkN-2}$. By  
\eqref{eqn:toolpf-o-(iii)} and \eqref{eqn:toolpf-iii-2}, we have 
\begin{equation*}
z = \bq^{4(c_{\rkN-1}(\lambda)-c_{r(\lambda)}(\lambda))}
\lambda_{\rkN-1}\lambda_\rkN\prod_{i=r(\lambda)+1}^{\rkN-2}
\lambda_i^2.
\end{equation*}
Hence by \eqref{eqn:toolpf-(-1)} and \eqref{eqn:toolpf-o-(i)}, 
$z\in\mcQ_\bq$,
contradiction. Assume $r(\lambda)=\rkN-1$. 
By \eqref{eqn:toolpf-o-(iv)}, $g_{\rkN-1}(\lambda)=\bq^{-2(m-1)}$.
By \eqref{eqn:toolpf-(-1)}, $c_{\rkN-1}(\lambda)= m-1$.
By \eqref{eqn:toolpf-iii-2} and \eqref{eqn:toolpf-o-(iii)}, $z={\frac {\lambda_\rkN} {\lambda_{\rkN-1}}}$.
Hence $\lambda\in Y_2$. Hence $Y_1\subseteq Y_2$.
Hence $Y_1= Y_2$. Thus we obtain the claim (2).

We show the claim (4).
Let $Y_3$ and $Y_4$ be LHS and RHS of 
\eqref{eqn:tlMain-4} respectively.
For $t\in\fkJ_{\rkN-m,\rkN}$,
let 
\begin{equation*}
\begin{array}{l}
Y_{4,t}:= \\
\quad\left\{
\begin{array}{ll}
\{\,\lambda\in\mcK^{(m)}_\bq\,|\,r(\lambda)=t,\,
\lambda_i=1\,(i\in\fkJ_{t+1,\rkN})\,\} & 
\mbox{if $t\in\fkJ_{\rkN-m,\rkN-2}$}, \\
\{\,\lambda\in\mcK^{(m)}_\bq\,|\,r(\lambda)=\rkN-1,\,
\lambda_{\rkN-1}=\lambda_\rkN\,\} & 
\mbox{if $t=\rkN-1$}, \\
\{\,\lambda\in\mcK^{(m)}_\bq\,|\,r(\lambda)\in\fkJ_{\rkN,\rkN+1},\,
{\tilde g}_\rkN(\lambda)=\bq^{-4m}\,\} &
\mbox{if $t=\rkN$}.
\end{array}\right.
\end{array}
\end{equation*}
By \eqref{eqn:toolpf-o-(i)}, \eqref{eqn:toolpf-o-(iv)} and \eqref{eqn:toolpf-iii-2}, we have
$Y_4=\cup_{t=\rkN-m}^\rkN Y_{4,t}$.
Then, by \eqref{eqn:toolpf-o-(i)}, \eqref{eqn:toolpf-o-(iii)} and \eqref{eqn:toolpf-iii-2}, we can easily see
$Y_4\subseteq Y_3$.
Let $\lambda\in Y_3$. By \eqref{eqn:toolpf-o-(iii)},  
\begin{equation}\label{eqn:toolpf-2}
{\tilde g}_\rkN(\lambda)=\bq^{-4c_{\rkN-1}(\lambda)}.
\end{equation}
Assume $r(\lambda)\in\fkJ_{\rkN,\rkN+1}$, that is, 
$c_{\rkN-1}(\lambda)=m$. By \eqref{eqn:toolpf-2}, $\lambda\in Y_{4,\rkN}$.
Assume $r(\lambda)=\rkN-1$. 
By \eqref{eqn:toolpf-iii-2} and \eqref{eqn:toolpf-2}, ${\frac {\lambda_\rkN} {\lambda_{\rkN-1}}}=1$.
Hence  $\lambda\in Y_{4,\rkN-1}$.
Assume $r(\lambda)\in\fkJ_{\rkN-m,\rkN-2}$.
By \eqref{eqn:toolpf-(-1)}, $c_{r(\lambda)}\leq c_{\rkN-1}$.
By \eqref{eqn:toolpf-o-(i)}, \eqref{eqn:toolpf-iii-2} and \eqref{eqn:toolpf-2},
we have $\lambda\in Y_{4,r(\lambda)}$. Thus we have $Y_3\subseteq Y_4$.
Hence $Y_3=Y_4$, as desired.
This completes the proof.
\hfill $\Box$

\subsection{Irreducible weights for super-$BD$ cases}
\label{subsection:IrrMdBD}

\begin{theorem}\label{theorem:IrrRepSuperB} 
Assume $\rkN\geq 2$. Let $m\in\fkJ_{1,\rkN-1}$.
Let $\bhm\in\bigchiPi$ be such that $\bhm\equiv\bhmp$
for some $\bhmp\in\dotXNSuper({\mathrm{B}}(m,\rkN-m))$.
Let $\bq:=\bhm(\al_\rkN,\al_\rkN)\in\bKtinf$.
Let $\Lambda\in\rmCh(U^0(\bhm))$,
and $g(\Lambda):=\prod_{i=\rkN-m}^\rkN\Lambda(K_{\al_i}L_{-\al_i})$. Then $\dim\mclL_\bhm(\Lambda)<\infty$
if and only if the following conditions $({\rm{irrB}}\mbox{-}m\mbox{-}1)$-$({\rm{irrB}}\mbox{-}m\mbox{-}3)$ are satisfied.
\newline\par
$({\rm{irrB}}\mbox{-}m\mbox{-}1)$ \quad $\Lambda\in\cap_{i\in\fkI\setminus\{\rkN-m\}}\mbbS_i(\bhm)$. \par
$({\rm{irrB}}\mbox{-}m\mbox{-}2)$ \quad $g(\Lambda)\in\{\bq^{-2x}|x\in\fkJ_{0,m-1}\}
\cup\{(-\bq)^{-(y+2m)}|y\in\bZgeqo\}$. \par
$({\rm{irrB}}\mbox{-}m\mbox{-}3)$ \quad If $g(\Lambda)=\bq^{-2x}$ for some $x\in\fkJ_{0,m-1}$,
then $\Lambda(K_{\al_i}L_{-\al_i})=1$ for all
$i\in\fkJ_{\rkN-m+x+1,\rkN}$.
\end{theorem}
{\it{Proof.}} See also Figure~\ref{fig:DynkinB}. 

By Theorem~\ref{theorem:cl}, we see that $\bhm$, 
$\bhmp\in\bigchiPifin$.
By Lemma~\ref{lemma:eqll}, we may assume $\bhm=\bhmp$.
Let $n:=\rkN^2$. 
Let $f\in\funcI_n$ be $\superBf_{m|\rkN-m}$ (see \eqref{eqn:superBf-def}).
By Lemmas~\ref{lemma:techlm}, \ref{lemma:longest-superB-id},
\ref{lemma:prlgelmmm-d} and 
\ref{lemma:SuCaTylm}~(2), we have 
$n=|R^+(\bhm)|$ and $\bhms_{f,n}=\lgst$.
By \eqref{eqn:inveta-superB}, \eqref{eqn:longest-superB-id-2} and \eqref{eqn:SuCa-eq-2}, 
we see the following. 
\begin{equation}\label{eqn:IRRsb-eq-1}
\bhm_{f,t-1}(\al_{f(t)},\al_{f(t)})\in\bKtinf
\quad (t\in\fkJ_{1,m^2+\rkN-m-1}).
\end{equation}
\begin{equation}\label{eqn:IRRsb-eq-2}
\bhm_{f,t-1}(\al_{f(t)},\al_{f(t)})=-1
\quad (t\in\fkJ_{m^2+\rkN-m,m^2+\rkN+m}\setminus\{m^2+\rkN\}).
\end{equation}
\begin{equation}\label{eqn:IRRsb-eq-3}
\bhm_{f,t-1}(\al_{f(t)},\al_{f(t)})=-\bq^{-1}
\quad\mbox{if $t=m^2+\rkN$.}
\end{equation}

By \eqref{eqn:superBf-def},
\begin{equation}\label{eqn:IRRsb-eq-4}
f(\fkJ_{1,m^2+\rkN-m-1})=\fkI\setminus\{\rkN-m\}.
\end{equation} 
By \eqref{eqn:hixl}, \eqref{eqn:hchcp}, \eqref{eqn:IRRsb-eq-1}, \eqref{eqn:IRRsb-eq-4}
and Lemma~\ref{lemma:crtypmo}, we see that 
\begin{equation}\label{eqn:IRRsb-eq-4a}
\mbox{$H(\bhm,\Lambda,f)\geq m^2+\rkN-m-1$ if and only if
$({\rm{irrB}}\mbox{-}m\mbox{-}1)$ holds.}
\end{equation}
By \eqref{eqn:hixl} and \eqref{eqn:IRRsb-eq-2},
$H(\bhm,\Lambda,f)\geq m^2+\rkN-m-1$ must be
$H(\bhm,\Lambda,f)\geq m^2+\rkN-1$. Hence
\begin{equation}\label{eqn:IRRsb-eq-4b}
\mbox{$H(\bhm,\Lambda,f)\geq m^2+\rkN-1$ if and only if 
$({\rm{irrB}}\mbox{-}m\mbox{-}1)$ holds.}
\end{equation}

Assume $H(\bhm,\Lambda,f)\geq m^2+\rkN-1$. 
(By \eqref{eqn:superBp-fom} and \eqref{eqn:SuCa-eq-2}, 
$\bhm_{f,m^2+\rkN-m-1}=\bhm$.)
By
\eqref{eqn:IRRsb-eq-1}
and Lemma~\ref{lemma:crtypmo}, we have
\begin{equation}\label{eqn:IRRsb-eq-5}
(\bhm_{f,m^2+\rkN-m-1},\Lambda_{\bhm,f,m^2+\rkN-m-1})\equiv
(\bhm,\Lambda).
\end{equation}
Let $\lambda:=(\Lambda(K_{\al_i}L_{-\al_i})|i\in\fkI)
\in(\bKt)^\rkN$.
By \eqref{eqn:IRRsb-eq-4b}, $\lambda\in\mcK^{(m)}_\bq$
(see also \eqref{eqn:dfKmq}).
Let $t\in\fkJ_{m^2+\rkN-m,m^2+\rkN-1}$. 
Let ${\bar h}_t:=h_{\bhm_{f,t-1},\Lambda_{\bhm,f,t-1},f(t)}$.
By \eqref{eqn:superBf-def}, 
\begin{equation}\label{eqn:IRRsb-eq-5a}
f(t)=t-m^2\,\,\mbox{and}\,\,f(t+1)=f(t)+1.
\end{equation}
By \eqref{eqn:IRRsb-eq-5a} and Lemmas~\ref{lemma:pre-longest-superB-id} and \ref{lemma:SuCaTylm}~(2),
\begin{equation}\label{eqn:IRRsb-eq-5b}
s^{\newtr_{f(t)}\bhm_{f,t-1}}_{f(t)}(\al_{f(t)+1})
=\al_{f(t)}+\al_{f(t)+1}.
\end{equation}
Using induction, we have \eqref{eqn:IRRsb-eq-6} and \eqref{eqn:IRRsb-eq-7} below.
\begin{equation}\label{eqn:IRRsb-eq-6}
\begin{array}{l}
\Lambda_{\bhm,f,t}(K_{\al_{f(t+1)}}L_{-\al_{f(t+1)}}) \\
\quad =\newtr_{f(t)}^{\bhm_{f,t-1}}\Lambda_{\bhm,f,t-1}(K_{\al_{f(t)+1}}L_{-\al_{f(t)+1}})
\quad(\mbox{by \eqref{eqn:IRRsb-eq-5a} and Definition~\ref{definition:fchi}}) \\
\quad =\Lambda_{\bhm,f,t-1}(K_{\al_{f(t)}+\al_{f(t)+1}}L_{-(\al_{f(t)}+\al_{f(t)+1})})
{\frac {\bhm_{f,t-1}(\al_{f(t)},-(\al_{f(t)}+\al_{f(t)+1}))^{{\bar h}_t}} 
{\bhm_{f,t-1}(\al_{f(t)}+\al_{f(t)+1},\al_{f(t)})^{{\bar h}_t}}} \\
\quad\quad\quad(\mbox{by \eqref{eqn:dfilm}
and \eqref{eqn:IRRsb-eq-5b}})
\\
\quad =\Lambda_{\bhm,f,t-1}(K_{\al_{f(t)}+\al_{f(t)+1}}L_{-(\al_{f(t)}+\al_{f(t)+1})})\bq^{2{\bar h}_t} \\
\quad\quad\quad(\mbox{by \eqref{eqn:IRRsb-eq-2},
since $\bhm_{f,t-1}(\al_{f(t)},\al_{f(t)+1})\bhm_{f,t-1}(\al_{f(t)+1},\al_{f(t)})=\bq^{-2}$} \\
\quad\quad\quad\,\,\mbox{by \eqref{eqn:superBp-fom} and \eqref{eqn:SuCa-eq-2}}) \\
\quad =\Lambda_{\bhm,f,t-1}(K_{\al_{f(t)}+\al_{f(t)+1}}L_{-(\al_{f(t)}+\al_{f(t)+1})})\bq^{2(1-\delta(
\Lambda_{\bhm,f,t-1}(K_{\al_{f(t)}}L_{-\al_{f(t)}}),1))} \\
\quad\quad\quad (\mbox{by \eqref{eqn:hixl} and \eqref{eqn:IRRsb-eq-2}}) \\
\quad =\toolmap^{(\bq,m)}_{f(t)+1}(\lambda)
\quad(\mbox{by \eqref{eqn:pretoolMain} and \eqref{eqn:IRRsb-eq-7} below}). \\
\end{array}
\end{equation} Similarly as above, we have
\begin{equation}\label{eqn:IRRsb-eq-7}
\Lambda_{\bhm,f,t-1}(K_{\al_i}L_{-\al_i})=\Lambda(K_{\al_i}L_{-\al_i})
\quad(i\in\fkJ_{f(t)+1,\rkN}),
\end{equation} where use \eqref{eqn:IRRsb-eq-5} if $t=m^2+\rkN-m$.

Since $\bq\in\bKtinf$,
by \eqref{eqn:hixl}, \eqref{eqn:IRRsb-eq-3}, \eqref{eqn:IRRsb-eq-5a} and \eqref{eqn:IRRsb-eq-6} for $t=m^2+\rkN-1$, we see that
\begin{equation}\label{eqn:IRRsb-eq-p8}
\begin{array}{l}
\mbox{$H(\bhm,\Lambda,f)\geq m^2+\rkN$  if and only if 
$({\rm{irrB}}\mbox{-}m\mbox{-}1)$ holds and there exists} \\
\mbox{$x\in\bZgeqo$ with
$\toolmap^{(\bq,m)}_\rkN(\lambda)=(-\bq^{-1})^x$.}
\end{array}
\end{equation}
By \eqref{eqn:tlMain-1}, \eqref{eqn:tlMain-3}
and \eqref{eqn:IRRsb-eq-p8}, we can see that
\begin{equation}\label{eqn:IRRsb-eq-8}
\mbox{$H(\bhm,\Lambda,f)\geq m^2+\rkN$  if and only if 
$({\rm{irrB}}\mbox{-}m\mbox{-}1)$-$({\rm{irrB}}\mbox{-}m\mbox{-}3)$ hold.}
\end{equation}

Assume $H(\bhm,\Lambda,f)\geq m^2+\rkN$.
By \eqref{eqn:hixl} and \eqref{eqn:IRRsb-eq-2}, we see that
\begin{equation}\label{eqn:IRRsb-eq-8b}
H(\bhm,\Lambda,f)\geq m^2+\rkN+m.
\end{equation}
(By \eqref{eqn:superBp-fom} and \eqref{eqn:SuCa-eq-2}, 
$\bhm_{f,m^2+\rkN-1}=\bhm_{f,m^2+\rkN}$.)
By \eqref{eqn:IRRsb-eq-2} and Lemma~\ref{lemma:crtypmo},
\begin{equation}\label{eqn:IRRsb-eq-9}
(\bhm_{f,m^2+\rkN-1},\Lambda_{\bhm,f,m^2+\rkN-1})\equiv
(\bhm_{f,m^2+\rkN},\Lambda_{\bhm,f,m^2+\rkN}).
\end{equation}
(By \eqref{eqn:superBp-fom} and \eqref{eqn:SuCa-eq-2}, 
$\bhm=\bhm_{f,m^2+\rkN+m}$.)
By \eqref{eqn:superBf-def}, 
\begin{equation}\label{eqn:IRRsb-eq-10}
f(t)=f(2(m^2+\rkN)-t)\quad(t\in\fkJ_{m^2+\rkN+1,m^2+\rkN+m}).
\end{equation}
For $t\in\fkJ_{m^2+\rkN+1,m^2+\rkN+m}$, we inductively see
\begin{equation*}
\begin{array}{l}
(\bhm_{f,t},\Lambda_{\bhm,f,t}) \\
\quad = (\newtr_{f(t)}\bhm_{f,t-1},\newtr^{\bhm_{f,t-1}}_{f(t)}\Lambda_{\bhm,f,t-1}) \\
\quad\quad\quad (\mbox{by Notation~\ref{notation:Ab-chift-1}~(1) and Definition~\ref{definition:fchi}}) \\
\quad \equiv (\newtr_{f(t)}\bhm_{f,2(m^2+\rkN)-t},\newtr^{\bhm_{f,2(m^2+\rkN)-t}}_{f(t)}\Lambda_{\bhm,f,2(m^2+\rkN)-t}) 
\\
\quad\quad\quad (\mbox{by induction and \eqref{eqn:lgelmrdd};
use \eqref{eqn:IRRsb-eq-9} if $t=m^2+\rkN+1$}) \\
\quad \equiv (\newtr_{f(t)}\newtr_{f(t)}\bhm_{f,2(m^2+\rkN)-t-1},
\newtr^{\bhm_{f,2(m^2+\rkN)-t}}_{f(t)}\newtr^{\bhm_{f,2(m^2+\rkN)-t-1}}_{f(t)}\Lambda_{\bhm,f,2(m^2+\rkN)-t-1}) \\
\quad\quad\quad (\mbox{by Notation~\ref{notation:Ab-chift-1}~(1), Definition~\ref{definition:fchi}
and \eqref{eqn:IRRsb-eq-10}}) \\
\quad = (\bhm_{f,2(m^2+\rkN)-t-1},
\Lambda_{\bhm,f,2(m^2+\rkN)-t-1}) \quad (\mbox{by \eqref{eqn:ichi} and \eqref{eqn:lgelmrd}}).
\end{array}
\end{equation*}
In particular, 
we have
\begin{equation}\label{eqn:IRRsb-eq-11}
(\bhm_{f,m^2+\rkN+m},\Lambda_{\bhm,f,m^2+\rkN+m})\equiv
(\bhm_{f,m^2+\rkN-m-1},\Lambda_{\bhm,f,m^2+\rkN-m-1}).
\end{equation}

By \eqref{eqn:IRRsb-eq-8} and \eqref{eqn:IRRsb-eq-8b}, we see that
\begin{equation}\label{eqn:IRRsb-eq-12}
\mbox{$H(\bhm,\Lambda,f)\geq m^2+\rkN+m$  if and only if 
$({\rm{irrB}}\mbox{-}m\mbox{-}1)$-$({\rm{irrB}}\mbox{-}m\mbox{-}3)$ hold.}
\end{equation}
By \eqref{eqn:superBf-def}, 
$f(t)=f(t-(\rkN+m))$ ($t\in\fkJ_{m^2+\rkN+m+1,n}$).
Hence by \eqref{eqn:IRRsb-eq-11} and \eqref{eqn:IRRsb-eq-12}, we see that
\begin{equation}\label{eqn:IRRsb-eq-13}
\mbox{$H(\bhm,\Lambda,f)=n$  if and only if 
$({\rm{irrB}}\mbox{-}m\mbox{-}1)$-$({\rm{irrB}}\mbox{-}m\mbox{-}3)$ hold.}
\end{equation}
By \eqref{eqn:IRRsb-eq-13} and Lemma~\ref{lemma:lgelmd}, this theorem holds.
\hfill $\Box$

\begin{figure}
\begin{center}
\figureA
\end{center}
\caption{Dynkin diagrams of $\bhm=\bhm_{f,0}\equiv\bhmp\in
\dotXNSuper({\mathrm{A}}(m-1,\rkN-m))$
with $\rkN=4$ and $m=1$, and $\bhm_{f,u}$ }
\label{fig:DynkinA}
\end{figure}

\begin{figure}
\begin{center}
\figureC
\end{center}
\caption{Dynkin diagrams of $\bhm=\bhm_{f,0}\equiv\bhmp\in
\dotXNSuper({\mathrm{C}}(\rkN))$  with $\rkN=4$, and $\bhm_{f,u}$
with $f=\ospefCN$}
\label{fig:DynkinC}
\end{figure}

\begin{theorem}\label{theorem:IrrRepSuperD} 
Assume $\rkN\geq 3$. Let $m\in\fkJ_{2,\rkN-1}$.
Let $\bhm\in\bigchiPi$ be such that $\bhm\equiv\bhmp$
for some $\bhmp\in\dotXNSuper({\mathrm{D}}(m,\rkN-m))$.
Let $\bq\in\bKtinf$ be such that 
$\bq^2=\bhm(\al_\rkN,\al_\rkN)$.
Let $\Lambda\in\rmCh(U^0(\bhm))$,
and 
\begin{equation*}
{\tilde g}(\Lambda):=(\prod_{i=\rkN-m}^{\rkN-2}\Lambda(K_{\al_i}L_{-\al_i}))^2
\Lambda(K_{\al_{\rkN-1}}L_{-\al_{\rkN-1}})\Lambda(K_{\al_\rkN}L_{-\al_\rkN}).
\end{equation*} Then $\dim\mclL_\bhm(\Lambda)<\infty$
if and only if the following conditions $({\rm{irrD}}\mbox{-}m\mbox{-}1)$-$({\rm{irrD}}\mbox{-}m\mbox{-}4)$ are satisfied.
\newline\par
$({\rm{irrD}}\mbox{-}m\mbox{-}1)$ \quad $\Lambda\in\cap_{i\in\fkI\setminus\{\rkN-m\}}\mbbS_i(\bhm)$. \par
$({\rm{irrD}}\mbox{-}m\mbox{-}2)$ \quad ${\tilde g}(\Lambda)=\bq^{-4x}$ for some $x\in\bZgeqo$. \par
$({\rm{irrD}}\mbox{-}m\mbox{-}3)$ \quad If ${\tilde g}(\Lambda)=\bq^{-4y}$ for some $y\in\fkJ_{0,m-2}$,
then $\prod_{i=\rkN-m}^{\rkN-m+y}\Lambda(K_{\al_i}L_{-\al_i})=\bq^{-2y}$ 
and $\Lambda(K_{\al_j}L_{-\al_j})=1$ for all
$j\in\fkJ_{\rkN-m+y+1,\rkN}$. \par
$({\rm{irrD}}\mbox{-}m\mbox{-}4)$ \quad If ${\tilde g}(\Lambda)=\bq^{-4(m-1)}$,
then $\prod_{i=\rkN-m}^{\rkN-1}\Lambda(K_{\al_i}L_{-\al_i})=\bq^{-2(m-1)}$.
\end{theorem}
{\it{Proof}}. See also Figure~\ref{fig:DynkinD}.
Let $n:=\rkN^2-m$.
Let $f\in\funcI_n$ be $\ospef_{m|\rkN-m}$ defined by \eqref{eqn:defddospef}.
Then, using this $f$,
this theorem can be proved in a way similar to that for Theorem~\ref{theorem:IrrRepSuperB}.
Here we only mention the facts as follows. 
By Lemmas~\ref{lemma:lgslm}~(1) and \ref{lemma:SuCaTylm}~(4) and \eqref{eqn:evenosp-lngst}, 
${\tilde \ell}_\bhm(1^\bhm w_0)=|R^+(\bhm)|=n$
and $\bhms_{f,n}=\lgst$.
Let $r:=m(m-1)+\rkN$.
By \eqref{eqn:longest-superD-id} and \eqref{eqn:SuCa-eq-4},
we have $\bhm_{f,{t_1}-1}(\al_{f(t_1)},\al_{f(t_1)})\in\bKtinf$
($t_1\in\fkJ_{1,r-m-1}$) and $\bhm_{f,{t_2}-1}(\al_{f(t_2)},\al_{f(t_2)})=-1$
($t_2\in\fkJ_{r-m,r-1}$). Then,
similarly to \eqref{eqn:IRRsb-eq-4b},
we can see that
$H(\bhm,\Lambda,f)\geq r-1$ if and only if 
$({\rm{irrD}}\mbox{-}m\mbox{-}1)$ holds.
Note that $f(r)=\rkN$.
By \eqref{eqn:longest-superD-id} and \eqref{eqn:SuCa-eq-4},
$\bhm_{f,r-1}(\al_\rkN,\al_\rkN)=\bq^{-4}\in\bKtinf$.
Similarly to \eqref{eqn:IRRsb-eq-6},
letting $\lambda:=(\Lambda(K_{\al_i}L_{-\al_i})|i\in\fkI)$,  we have
$\Lambda_{\bhm,f,r-1}(K_{\al_\rkN}L_{-\al_\rkN})
=\tltoolmap^{(\bq,m)}_\rkN(\lambda)$.
Similarly to \eqref{eqn:IRRsb-eq-6}, 
by \eqref{eqn:tlMain-2} and \eqref{eqn:tlMain-4},
we can see that
$H(\bhm,\Lambda,f)\geq r-1$ if and only if 
$({\rm{irrD}}\mbox{-}m\mbox{-}1)$-$({\rm{irrD}}\mbox{-}m\mbox{-}4)$ hold.
Similarly to \eqref{eqn:IRRsb-eq-13}, 
we can see that
$H(\bhm,\Lambda,f)=n$ if and only if 
$({\rm{irrD}}\mbox{-}m\mbox{-}1)$-$({\rm{irrD}}\mbox{-}m\mbox{-}4)$ hold.
Hence, by Lemma~\ref{lemma:lgelmd}, this theorem holds.
\hfill $\Box$

\begin{figure}
\begin{center}
\figureB
\end{center}
\caption{Dynkin diagrams of $\bhm=\bhm_{f,0}\equiv\bhmp\in\dotXNSuper({\mathrm{B}}(m,\rkN-m))$
with $\rkN=4$ and $m=2$, and $\bhm_{f,u}$ 
with $f=\superBf_{m|\rkN-m}$}
\label{fig:DynkinB}
\end{figure}

\begin{figure}
\begin{center}
\figureD
\end{center}
\caption{Dynkin diagrams of $\bhm=\bhm_{f,0}\equiv\bhmp\in\dotXNSuper({\mathrm{D}}(m,\rkN-m))$
with $\rkN=5$ and $m=3$, and $\bhm_{f,u}$ with $f=\ospef_{m|\rkN-m}$}
\label{fig:DynkinD}
\end{figure}

\subsection{Irreducible weights for super-$FG$ and extra cases}
\label{subsection:IrrMdFG}

\begin{theorem} \label{theorem:MainSec} Let $\bhm\in\bigchiPi$,
$\Lambda\in\rmCh(U^0(\bhm))$
and $\lambda_i:=\Lambda(K_{\al_i}L_{-\al_i})$
$(i\in\fkI)$.

{\rm{(1)}} Assume that $\rkN=4$
and $\bhm\equiv\bhmp$
for some $\bhmp\in\dotXNSuper({\mathrm{F}}(4))$.
Let $\bq\in\bKtinf$ be such that 
$\bq^2=\bhm(\al_2,\al_2)$. Then $\dim\mclL_\bhm(\Lambda)<\infty$
if and only if one of the following conditions 
$({\rm{irrF}}\mbox{-}1)$-$({\rm{irrF}}\mbox{-}4)$ holds.
\newline\par
$({\rm{irrF}}\mbox{-}1)$ \quad $\Lambda\in\cap_{i=2}^4\mbbS_i(\bhm)$ and 
$\lambda_1^2\lambda_2^3\lambda_3^2\lambda_4=\bq^{-6(x+4)}$
for some $x\in\bZgeqo$. \par
$({\rm{irrF}}\mbox{-}2)$ \quad $\lambda_1=\lambda_2=\lambda_3=\lambda_4=1$. \par
$({\rm{irrF}}\mbox{-}3)$ \quad $\Lambda\in\mbbS_3(\bhm)$, 
$\lambda_2=\lambda_4=1$ and $\lambda_1\lambda_3=\bq^{-6}$. \par
$({\rm{irrF}}\mbox{-}4)$ \quad $\Lambda\in\cap_{i=3}^4\mbbS_i(\bhm)$, $\lambda_1\lambda_3\lambda_4^2=\bq^{-12}$ and $\lambda_2=\bq^2\lambda_4$.
\newline\par
{\rm{(2)}} Assume that $\rkN=3$
and $\bhm\equiv\bhmp$
for some $\bhmp\in\dotXNSuper({\mathrm{G}}(3))$.
Let $\bq\in\bKtinf$ be such that 
$\bq^2=\bhm(\al_2,\al_2)$. Then $\dim\mclL_\bhm(\Lambda)<\infty$
if and only if one of the following conditions 
$({\rm{irrG}}\mbox{-}1)$-$({\rm{irrG}}\mbox{-}3)$ holds.
\newline\par
$({\rm{irrG}}\mbox{-}1)$ \quad $\Lambda\in\cap_{i=2}^3\mbbS_i(\bhm)$ and
$\lambda_1\lambda_2^2\lambda_3=(-\bq^{-2})^{x+6}$
for some $x\in\bZgeqo$.  \par
$({\rm{irrG}}\mbox{-}2)$ \quad $\lambda_1=\lambda_2=\lambda_3=1$. \par
$({\rm{irrG}}\mbox{-}3)$ \quad $\Lambda\in\mbbS_3(\bhm)$, $\lambda_1\lambda_3=\bq^{-8}$ and $\lambda_2=1$. 
\newline\par
{\rm{(3)}} Assume that $\rkN=3$
and $\bhm\equiv\bhmp$
for some $\bhmp\in\dotXNExtra$.
Let $\bq:=\bhm(\al_1,\al_1)\in\bKtinf$ and $\br:=\bhm(\al_3,\al_3)\in\bKt\setminus\{1,\bq^{-1}\}$.
If $\bq\br\notin\bKtinf$ and $\br\in\bKtinf$
{\rm{(}}resp.  $\br\notin\bKtinf${\rm{)}}, then $\dim\mclL_\bhm(\Lambda)<\infty$
if and only if $\Lambda\in\mbbS_1(\bhm)\cap\mbbS_3(\bhm)$
{\rm{(}}resp.  $\Lambda\in\mbbS_1(\bhm)${\rm{)}}. 
If $\bq\br\in\bKtinf$ and $\br\in\bKtinf$ 
{\rm{(}}resp.  $\br\notin\bKtinf${\rm{)}}, then $\dim\mclL_\bhm(\Lambda)<\infty$
if and only if $\Lambda\in\mbbS_1(\bhm)\cap\mbbS_3(\bhm)$ 
{\rm{(}}resp.  $\Lambda\in\mbbS_1(\bhm)${\rm{)}}
and one of the following conditions 
$({\rm{irrEx3}}\mbox{-}1)$-$({\rm{irrEx3}}\mbox{-}4)$ holds.
\newline\par
$({\rm{irrEx3}}\mbox{-}1)$ \quad 
$\lambda_1\lambda_2^2\lambda_3=(\bq\br)^{-(x+2)}$
for some $x\in\bZgeqo$. \par
$({\rm{irrEx3}}\mbox{-}2)$ \quad
$\lambda_1=\lambda_2=\lambda_3=1$. \par
$({\rm{irrEx3}}\mbox{-}3)$ \quad $\lambda_2=1$ and $\lambda_1\lambda_3=(\bq\br)^{-1}$. \par
$({\rm{irrEx3}}\mbox{-}4)$ \quad $\lambda_1\lambda_2=\bq^{-1}$ and $\lambda_2\lambda_3=\br^{-1}$. 
\newline\par
{\rm{(4)}} Assume that $\rkN=2$
and $\bhm\equiv\bhmp$
for some $\bhmp\in\dotXNExtra$.
Let $\bq:=\bhm(\al_2,\al_2)\in\bKtinf$ and $\bzeta:=\bhm(\al_1,\al_1)\in\bKt_3$.
Then $\dim\mclL_\bhm(\Lambda)<\infty$
if and only if one of the following conditions 
$({\rm{irrEx2}}\mbox{-}1)$-$({\rm{irrEx2}}\mbox{-}2)$ holds.
\newline\par
$({\rm{irrEx2}}\mbox{-}1)$ \quad $\Lambda\in\mbbS_2(\bhm)$ and
$\lambda_1^2\lambda_2=(\bzeta\bq^{-1})^{x+2}$
for some $x\in\bZgeqo$.  \par
$({\rm{irrEx2}}\mbox{-}2)$ \quad $\lambda_1=\lambda_2=1$. 
\newline\par
{\rm{(5)}} Assume that $\rkN=4$
and $\bhm\equiv\bhmp$
for some $\bhmp\in\dotXNExtra$.
Let $\bq:=\bhm(\al_2,\al_2)\in\bKtinf$.
Then $\dim\mclL_\bhm(\Lambda)<\infty$ if and only if one of the following conditions 
$({\rm{irrEx4}}\mbox{-}1)$-$({\rm{irrEx4}}\mbox{-}5)$ holds.
\newline\par
$({\rm{irrEx4}}\mbox{-}1)$ \quad $\Lambda\in\mbbS_1(\bhm)\cap\mbbS_2(\bhm)\cap\mbbS_4(\bhm)$ and
$\lambda_1\lambda_2^2\lambda_3^3\lambda_4=(-\bq^{-1})^{x+3}$
for some $x\in\bZgeqo$. \par 
$({\rm{irrEx4}}\mbox{-}2)$ \quad $\lambda_1=\lambda_2=\lambda_3=\lambda_4=1$.
\par 
$({\rm{irrEx4}}\mbox{-}3)$ \quad $\lambda_2=\lambda_3=1$
and $\lambda_1=\bq^{2x}$, $\lambda_4=(-\bq^{-1})^{2x+1}$ for some $x\in\bZgeqo$.
\par 
$({\rm{irrEx4}}\mbox{-}4)$ \quad $\lambda_3=1$ and there exist
$x$, $y\in\bZgeqo$ such that $\lambda_1=\bq^{2x}$, $\lambda_2=\bq^y$ and 
$\lambda_4=\bq^{-2(x+y+1)}$.
\par 
$({\rm{irrEx4}}\mbox{-}5)$ \quad There exist
$x\in\bZgeqo$ and $\fkJ_{0,{\frac x 2}}$ such that 
$\lambda_1=\bq^{x+1}$, $\lambda_2=\bq^{x-2y}$, $\lambda_3=\bq^{3(-x+2y-1)}$ and 
$\lambda_4=\bq^{-2(x+y+1)}$.
\end{theorem}
{\it{Proof.}}
We define $n$, $f\in\funcI_n$ and a map $z:\fkJ_{1.n}\to\bKt$ as follows.
Let ${\hat f}:=(f(t)|t\in\fkJ_{1,n})\in\fkI^n$ and ${\hat z}:=(z(t)|t\in\fkJ_{1,n})\in(\bKt)^n$.

If $\rkN=4$ and $\bhm\equiv\bhmp\in\dotXNSuper({\mathrm{F}}(4))$, let 
$n:=18$, \newline ${\hat f}:=(2,3,4,2,3,4,2,3,4,
1,2,3,4,1,4,3,2,1)$ and \newline  
${\hat z}:=(\bq^2,\bq^4,\bq^4,\bq^2,\bq^4,\bq^4,\bq^2,\bq^4,\bq^4,
-1,-1,-1,-1,\bq^{-6},-1,-1,-1,-1)$.

If $\rkN=3$ and $\bhm\equiv\bhmp\in\dotXNSuper({\mathrm{G}}(3))$, let 
$n:=13$, \newline ${\hat f}:=(2,3,2,3,2,3,1,2,3,1,3,2,1)$ and \newline 
${\hat z}:=(\bq^2,\bq^6,\bq^2,\bq^6,\bq^2,\bq^6, 
-1,-1,-1,-\bq^{-2},-1,-1,-1)$.

If $\rkN=3$ and $\bhm\equiv\bhmp\in\dotXNExtra$, let
$n:=7$, ${\hat f}:=(1,3,2,1,3,1,2)$ and 
${\hat z}:=(\bq,\br,-1,-1,(\bq\br)^{-1},-1,-1)$.

If $\rkN=2$ and $\bhm\equiv\bhmp\in\dotXNExtra$, let
$n:=4$, ${\hat f}:=(2,1,2,1)$ and 
${\hat z}:=(\bq,\bzeta,\bzeta\bq^{-1},\bzeta)$.

If $\rkN=4$ and $\bhm\equiv\bhmp\in\dotXNExtra$, let
$n:=15$, \newline ${\hat f}:=(1,2,1,4,3,4,2,1,4,3,1,2,4,2,1)$ and \newline 
${\hat z}:=(\bq,\bq,\bq,-\bq^{-1},
-1,-1,-1,-1,-\bq^{-1},-\bq^{-1},-1,-1,-1,-1,-1)$.

Using Lemmas~\ref{lemma:lgelmmm} and \ref{lemma:prlgelmmm-d}
(see also Figure~\ref{fig:DynkinF}, \ref{fig:DynkinG},
\ref{fig:DynkinDX}, \ref{fig:DynkinZT} or \ref{fig:DynkinRF}), we
can directly see that 
\begin{equation}\label{eqn:hatfwo}
n=|R^+(\bhm)|\quad\mbox{and}\quad\bhms_{{\hat f},n}=\lgst.
\end{equation}
We can also see that ${\hat z}(t)=\bhm_{f,t-1}(\al_{f(t)},\al_{f(t)})$
($t\in\fkJ_{1,n}$) (see also Figures~\ref{fig:DynkinF}, \ref{fig:DynkinG},
\ref{fig:DynkinDX}, \ref{fig:DynkinZT} or \ref{fig:DynkinRF}).

Let $r:=\rmMax\{\,t\in\fkJ_{1,n}\,|\,
\forall t^\prime\in\fkJ_{1,t},\,{\hat z}(t)\in \bKtinf\,\}$, and
$b:=|\{t\in\fkJ_{b+1,n}|{\hat z}(t)\in \bKtinf\}|$.
Then $b\leq 2$, and $b=2$ if and only if
$\rkN=4$ and $\bhm\equiv\bhmp\in\dotXNExtra$.
Using these facts, we can directly prove this theorem.

Let us explain more precisely how to prove the claim (5).
Note $r=4$. Since ${\hat z}(t)=-1\notin\bKtinf$ ($t\in\fkJ_{5,8}$),
by an argument similar to that for 
\eqref{eqn:IRRsb-eq-4a}, 
we see that
\begin{equation}\label{eqn:IRR-RF-1}
\mbox{$H(\bhm,\Lambda,f)\geq 8$ if and only if $\Lambda\in\mbbS_1(\bhm)\cap\mbbS_2(\bhm)\cap\mbbS_4(\bhm)$.}
\end{equation}

Assume $H(\bhm,\Lambda,f)\geq 8$. By \eqref{eqn:IRR-RF-1},
$(\lambda_1,\lambda_2,\lambda_4)
=(\bq^{l_1}, \bq^{l_2},(-\bq^{-1})^{l_4})$
for some $(l_1,l_2,l_4)\in(\bZgeqo)^3$.
Let $h_t:=h_{\bhm_{f,t-1},\Lambda_{f,t-1},f(t)}$
for $t\in\fkJ_{1,H(\bhm,\Lambda,f)}$.
By Lemma~\ref{lemma:crtypmo}, 
$h_t=l_{f(t)}$
($t\in\fkJ_{1,4}$),
$h_5=1-\delta(1,\lambda_3)$,
$h_6=1-\delta(1,(-\bq^{-1})^{h_5}\lambda_3\lambda_4)$,
$h_7=1-\delta(1,(-1)^{h_5+h_6}\lambda_2\lambda_3^2\lambda_4)$, and
$h_8=1-\delta(1,(-1)^{h_5+h_6}\bq^{h_7}\lambda_1\lambda_2\lambda_3^2\lambda_4)$.
We see that
\begin{equation}\label{eqn:Hpnine}
\begin{array}{l}
H(\bhm,\Lambda,f)\geq 9 \,\mbox{if and only if}\,\,\\
(-\bq^{-1})^c=(-\bq^{-1})^{-(h_5+h_7+h_8)}\lambda_1\lambda_2^2\lambda_3^3\lambda_4\,
\mbox{for some}\,c\in\bZgeqo.
\end{array}
\end{equation} 
We can see that
if $H(\bhm,\Lambda,f)\geq 9$, then
\begin{equation}\label{eqn:Hpten}
H(\bhm,\Lambda,f)\geq 10 \quad\Longleftrightarrow\quad
\exists c^\prime\in\bZgeqo,\,\,(-\bq^{-1})^{c^\prime}=(-\bq^{-1})^{h_5-h_6}\lambda_4.
\end{equation}
By \eqref{eqn:Hpnine} and \eqref{eqn:Hpten},
we see that $H(\bhm,\Lambda,f)\geq 9$ must be $H(\bhm,\Lambda,f)\geq 10$,
since, if $h_5=l_4=0$, then $\lambda_3=1$ and $h_6=0$.
Since ${\hat z}(t)=-1\notin\bKtinf$ ($t\in\fkJ_{11,n}$), $H(\bhm,\Lambda,f)\geq 10$
must be $H(\bhm,\Lambda,f)=n$. 
Then using Lemma~\ref{lemma:lgelmd} and \eqref{eqn:Hpnine}, by a direct argument, we can see the claim (5) holds.

This completes the proof.
\hfill $\Box$

\begin{figure}
\begin{center}
\figureF
\end{center}
\caption{Dynkin diagrams of $\bhm=\bhm_{f,0}\equiv\bhmp\in\dotXNSuper({\mathrm{F}}(4))$
with $\rkN=4$, and $\bhm_{f,u}$ with $f={\hat f}$}
\label{fig:DynkinF}
\end{figure}
\begin{figure}
\begin{center}
\figureG
\end{center}
\caption{Dynkin diagrams of $\bhm=\bhm_{f,0}\equiv\bhmp\in\dotXNSuper({\mathrm{G}}(3))$
with $\rkN=3$, and $\bhm_{f,u}$ with $f={\hat f}$}
\label{fig:DynkinG}
\end{figure}
\begin{figure}
\begin{center}
\figureDX
\end{center}
\caption{Dynkin diagrams of $\bhm=\bhm_{f,0}\equiv\bhmp\in\dotXNExtra$
with $\rkN=3$, and $\bhm_{f,u}$ with $f={\hat f}$}
\label{fig:DynkinDX}
\end{figure}
\begin{figure}
\begin{center}
\figureZT
\end{center}
\caption{Dynkin diagrams of $\bhm=\bhm_{f,0}\equiv\bhmp\in\dotXNExtra$
with $\rkN=2$, and $\bhm_{f,u}$ with $f={\hat f}$}
\label{fig:DynkinZT}
\end{figure}
\begin{figure}
\begin{center}
\figureRF
\end{center}
\caption{Dynkin diagrams of $\bhm=\bhm_{f,0}\equiv\bhmp\in\dotXNExtra$
with $\rkN=4$, and $\bhm_{f,u}$ with $f={\hat f}$}
\label{fig:DynkinRF}
\end{figure}

\begin{remark}\label{remark:last}
(1) Let ${\widehat U}(\bhm_{(b)})$ be the unital $\bK$-algebra
such that it contains $U(\bhm_{(b)})$ and the Laurent polynomial
$\bK$-algebra $X:=\bK [x_i^{\pm 1}|i\in\fkI]$ as subalgebras, the linear map
$U(\bhm_{(b)})\otimes X\to{\widehat U}(\bhm_{(b)})$
($z\otimes y\mapsto zy$) is bijective,
and $x_iK_\al x_i^{-1}=K_\al$, $x_iL_\al x_i^{-1}=L_\al$,
$x_iE_jx_i^{-1}=\bq^{\delta_{ij}}E_j$,
$x_iF_jx_i^{-1}=\bq^{-\delta_{ij}}F_j$, where recall $\bq\in \bKtinf$.
Then we have a bijective map
from $\{\,\Lambda\in\rmCh(U^0(\bhm_{(b)}))\,|\,\kp(\Lambda)\in S_{(b)}\,\}\times (\bKt)^\rkN$
to the set of equivalence classes of finite-dimensional ${\widehat U}(\bhm_{(b)})$-modules
sending $(\Lambda,\mu)$ to the equivalence class of irreducible
highest weight ${\widehat U}(\bhm_{(b)})$-modules
with highest vectors ${\widehat v}$ such
that $E_i {\widehat v}=0$, $x_i{\widehat v}=\mu_i{\widehat v}$
($i\in\fkI$) and
$K_\al L_\beta {\widehat v}=\Lambda (K_\al L_\beta ){\widehat v}$
($\al$, $\beta\in\bZPi$).
\par
(2) It seems to be not easy to recover 
Theorems~\ref{theorem:IrrRepSuperAC}, \ref{theorem:IrrRepSuperB}, 
\ref{theorem:IrrRepSuperD}, \ref{theorem:MainSec}~(1)-(2), 
from Geer's result \cite{G07} since definitions
of his and our quantum groups are not so close. \par
(3)
See \cite{RS93} for some deeper results concerning Theorem~\ref{theorem:IrrRepCartan}.
\end{remark}

\begin{remark}\label{remark:lasttwo} It is easy to be convinced that
our argument in this paper can also be applied to recover the Kac's list
(see Introduction) for
the simple Lie superalgebras.
Shu and Wang \cite[Theorem~5.3, Remark~5.4]{SW08} also recovered those for
the simple Lie superalgebras $B(m,n)$, $C(n)$
and $D(m,n)$ by using odd reflections in a way totally different from that in this paper.
Note that an intrinsic gap appears between the list for $B(m,n)$
and the one in Theorem~\ref{theorem:IrrRepSuperB}, since
$g(\Lambda)$ can be $(-\bq)^{-(2m+y)}$ for some $y\in\bN$.
\end{remark}

\vspace{1cm}
{Saeid Azam, School of Mathematics, Institute for Research in
Fundamental Sciences (IPM), P.O.Box: 19395-5746, Tehran, Iran,
and Department of Mathematics, University of Isfahan, 
P.O.Box: 81745-163, Isfahan,
Iran}, \newline {E-mail: azam@sci.ui.ac.ir} \vspace{0.5cm}

{Hiroyuki Yamane, Department of Mathematics,
Faculty of Science, University of Toyama,
 3190 Gofuku, Toyama-shi, Toyama 930-8555, Japan}, \newline
{E-mail: hiroyuki@sci.u-toyama.ac.jp}

\vspace{0.5cm}
{Malihe Yousofzadeh, \, Department of Mathematics, University of Isfahan, \, P.O.Box: 81745-163,
Isfahan, Iran, and Institute for Research in
Fundamental Sciences (IPM), P.O.Box: 19395-5746, Tehran, Iran},   \newline
{E-mail: ma.yousofzadeh@sci.ui.ac.ir}


\begin{thebibliography}{99}

\bibitem[1]{A11} I.~Angiono, {\it{On Nichols algebras of diagonal type}}, accepted in J. Reine Angew. Math.

\bibitem[2]{AAY10} N.~Andruskiewitsch, I. Angiono and H. Yamane,
{\it{On pointed Hopf superalgebras}},
Contemp. Math. 544 (2011), 123-140.


\bibitem[3]{AS98}  N.~Andruskiewitsch and H.-J.~Schneider, {\it{Lifting of quantum linear
spaces and pointed Hopf algebras of order $p^3$}}, J. Algebra 209 (1998), no. 2, 658-691.
658-691.

\bibitem[4]{AS10} \bysame, {\it{On the classification of
finite-dimensional pointed Hopf algebras}}, Annals of Mathematics Vol. 171 (2010), no. 1, 375-417.

\bibitem[5]{CH09}  M.~Cuntz and I.~Heckenberger, {\it{Weyl groupoids
with at most three objects}}, J. Pure Appl. Algebra 213 (2009), no. 6, 1112-1128.

\bibitem[6]{Dr86} V.G.~Drinfel'd, {\it{Quantum groups}}, 
Proceedings of the International Congress of Mathematicians, Vol. 1, 2 
(Berkeley, Calif., 1986), 798-820, Amer. Math. Soc., Providence, RI, 1987.

\bibitem[7]{G07} N.~Geer,
{\it{Some remarks on quantized Lie superalgebras of classical type}}, J.~Algebra 314 (2007) no. 2, 565-580.

\bibitem[8]{Hec06} I.~Heckenberger, {\it{The Weyl groupoid of a Nichols algebra of diagonal type}},
Invent. Math. 164 (2006), no. 1, 175-188.

\bibitem[9]{Hec09} \bysame, {\it{Classification of arithmetic root systems}},
 Adv. Math. 220 (2009), no. 1, 59-124.

\bibitem[10]{Hec10} \bysame, {\it{Lusztig isomorphisms for Drinfel'd doubles
of bosonizations of Nichols algebras of diagonal type}}, J. Algebra 323 (2010), no. 8, 2130-2182.

\bibitem[11]{HK07} I.~Heckenberger and S.~Kolb,
{\it{On the Bernstein-Gelfand-Gelfand resolution for Kac-Moody algebras and
quantized enveloping algebras}}, Transform.~Groups 12 (2007), no.4, 647-655.

\bibitem[12]{HY08} I.~Heckenberger and H.~Yamane, {\it{A generalization of Coxeter groups, root systems, and
Matsumoto's theorem}}, Math. Z. 259 (2008), no. 2, 255-276.



\bibitem[13]{HY10} \bysame, {\it{Drinfel'd doubles and Shapovalov determinants}},
Rev. Un. Mat. Argentina 51 (2010), no. 2, 107-146.

\bibitem[14]{Hum} J.E.~Humphreys, Reflection Groups and Coxeter Groups, Cambridge University Press, 1992

\bibitem[15]{Kac90}  V.G.~Kac, {\it{Infinite-Dimensional Lie Algebras
}}, 3rd ed., Cambridge University Press, 1990.

\bibitem[16]{Kac77}  \bysame, {\it{Lie superalgerbras}}, Adv.~in~Math., 26 (1977), no. 1,  8-96.

\bibitem[17]{Kha99} V.~Kharchenko, {\it{A quantum analogue of the Poincar\'{e}-Birkhoff-Witt
theorem}}, Algebra and Logic 38 (1999), no. 4, 259-276.

\bibitem[18]{b-Lusztig93}
G.~Lusztig, {\it{Introduction to quantum groups}}, Birkh{\"a}user, Boston, MA,
  1993.


\bibitem[19]{Lusztig90}
\bysame, {\it{On quantum groups}}, J. Algebra 131 (1990), no.~2, 466-475.

\bibitem[20]{RS93}  D.~Radford, H.-J.~Schneider, {\it{On the simple
representations of generalized quantum groups
and quantum doubles}}, J. Algebra 319 (2008), no. 9, 3689-3731.

\bibitem[21]{SW08} B.~Shu and W.~Wang,
{\it{Modular representations of the ortho-symplectic superalgebras}},
Proc. London Math. Soc., 96 (2008), no.~2, 251-271.



\bibitem[22]{Y89} H.~Yamane, {\it{A Poincar\'{e}-Birkhoff-Witt theorem for
quantized universal enveloping algebras of type
$A_N$}},
Publ. RIMS, Kyoto Univ., 25 (1989), no. 3, 503-520.

\bibitem[23]{Y94} \bysame, {\it{Quantized enveloping algebras associated to simple Lie superalgebras and their
universal $R$-matrices}}, Publ. RIMS Kyoto Univ. 30 (1994), no. 3, 15-87.

\bibitem[24]{Y99} \bysame, {\it{On defining relations of affine Lie superalgebras and affine
quantized universal enveloping superalgebras}}, Publ. RIMS Kyoto Univ. 35
(1999), no. 3, 321-390.

\bibitem[25]{Y01} \bysame, {\it{Errata to ``On defining relations of affine Lie superalgebras and
affine quantized universal enveloping superalgebras"}}, Publ. RIMS Kyoto
Univ. 37 (2001), no.~4, 615-619.


\bibitem[26]{Y06} \bysame, {\it{Representations of a $\bZ/3\bZ$-quantum group}},
Publ. RIMS, Kyoto Univ., 43 (2007), no. 1, 75-93.

\bibitem[27]{YHome} \bysame, {\it{Examples of the defining relations of the quantum affine superalgebras}},
available at http://www16.tok2.com/home/hiroyukipersonal



\end{thebibliography}
\end{document}